\documentclass[10pt]{amsart} 

\usepackage{amsmath}
\usepackage{amssymb}

\usepackage{graphicx}
\usepackage{xcolor}
\usepackage{subcaption}
\usepackage{multirow}
\usepackage{enumitem}
\usepackage{mathrsfs}
\usepackage[rightcaption]{sidecap}
\usepackage{rotating}

\newtheorem{theorem}{Theorem}[section]

\theoremstyle{definition}
\newtheorem{stencil}{Stencil}
\newtheorem{example}[theorem]{Example}
\theoremstyle{remark}

\newcommand{\bfu}{\mathbf u}
\newcommand{\bff}{\mathbf f}

\newcommand{\bfx}{\mathbf x}
\newcommand{\dx}{\mathrm d}
\newcommand{\vbar}{\overline{v}}
\newcommand{\ubar}{\overline{u}}
\newcommand{\JS}{\mathrm{JS}}
\newcommand{\M}{\mathrm M}
\newcommand{\Z}{\mathrm Z}
\newcommand{\ZR}{\mathrm{ZR}}
\newcommand{\ZL}{\mathrm{ZL}}
\newcommand{\N}{\mathcal N}
\newcommand{\Nc}{\N\!,\,}
\newcommand{\ceps}{\check{\epsilon}}
\newcommand{\teps}{\tilde{\epsilon}}
\newcommand{\vez}{\varepsilon_0}
\newcommand{\vet}{\varepsilon_2}
\newcommand{\dbp}{\mathrm{eps}} 
\newcommand{\e}{\mathrm e}
\newcommand{\cm}{\text{-}}
\newcommand{\cbeta}{\check{\beta}}
\newcommand{\tbeta}{\tilde{\beta}}
\newcommand{\tomega}{\tilde{\omega}}
\newcommand{\es}{\mathsf e}
\newcommand{\cfl}{\text{CFL}}

\newcommand{\A}{\mathscr A}
\newcommand{\B}{\mathscr B}
\newcommand{\C}{\mathscr C}
\newcommand{\D}{\mathscr D}
\newcommand{\E}{\mathscr E}

\title[A spatial-temporal weight analysis of WENO schemes]{A spatial-temporal weight analysis and novel nonlinear weights of weighted essentially non-oscillatory schemes for hyperbolic conservation laws}
\author{Xinjuan Chen}
\address{Department of Mathematics, College of Science, Jimei University, Xiamen, Fujian 361021, China}
\email{chenxinjuan@jmu.edu.cn}

\author{Jiaxi Gu}
\address{Department of Mathematics $\&$ POSTECH MINDS (Mathematical Institute for Data Science), Pohang University of Science and Technology, Pohang 37673, Korea}
\email{jiaxigu@postech.ac.kr}

\author{Jae-Hun Jung}
\address{Department of Mathematics $\&$ POSTECH MINDS (Mathematical Institute for Data Science), Pohang University of Science and Technology, Pohang 37673, Korea}
\email{jung153@postech.ac.kr}

\makeatletter
\@namedef{subjclassname@2020}{\textup{}2020 Mathematics Subject Classification}
\makeatother
\subjclass[2020]{65M08, 65M15}
\keywords{Weighted essentially non-oscillatory scheme, Finite volume, Runge-Kutta method, Nonlinear weights}

\begin{document}

\maketitle

\begin{abstract}
In this paper we analyze the weighted essentially non-oscillatory (WENO) schemes in the finite volume framework by examining the first step of the explicit third-order total variation diminishing Runge-Kutta method.
The rationale for the improved performance of the finite volume WENO-M, WENO-Z and WENO-ZR schemes over WENO-JS in the first time step is that the nonlinear weights corresponding to large errors are adjusted to increase the accuracy of numerical solutions.
Based on this analysis, we propose novel Z-type nonlinear weights of the finite volume WENO scheme for hyperbolic conservation laws.
Instead of taking the difference of the smoothness indicators for the global smoothness indicator, we employ the logarithmic function with tuners to ensure that the numerical dissipation is reduced around discontinuities while the essentially non-oscillatory property is preserved.
The proposed scheme does not necessitate substantial extra computational expenses.
Numerical examples are presented to demonstrate the capability of the proposed WENO scheme in shock capturing.
\end{abstract}

\section{Introduction} \label{sec:intro}
In this paper we are interested in the finite volume method for hyperbolic conservation laws
\begin{equation} \label{eq:hyperbolic}
\begin{aligned}
 \frac{\partial \bfu}{\partial t} + \nabla \cdot \bff(\bfu) &= 0, \\
                                             \bfu (\bfx, 0) &= \bfu_0(\bfx).
\end{aligned}
\end{equation}
The essentially non-oscillatory (ENO) and weighted ENO (WENO) schemes, designed with higher order of accuracy, give satisfactory numerical solutions to \eqref{eq:hyperbolic} with shocks and other discontinuities.
Harten et al. \cite{Harten} in 1986 introduced the ENO schemes in the finite volume framework, producing ENO transitions around the discontinuities.
In 1994, the finite volume WENO scheme was first presented by Liu et al. \cite{Liu}, where the $(r\! +\! 1)$th-order WENO scheme is constructed from the $r$th-order ENO scheme in order to increase one order of accuracy in smooth regions, while maintaining the ENO property for the sharp transitions.
Later in \cite{Jiang}, Jiang and Shu proposed the fifth-order finite difference WENO scheme with a general framework for the design of the smoothness indicators.
On the basis of the notion of the smoothness indicators, numerous finite volume WENO schemes were widely developed, including WENO schemes on structured meshes \cite{Shi,TitarevToro,QiuCF,Buchmuller}, on unstructured meshes \cite{Friedrich,HuShu,Shi,ZhangShu,TitarevTD,HuLiTang}, central WENO schemes \cite{Levy,QiuJCP,Capdeville,Kolb,Cravero}, radial basis function ENO/WENO schemes \cite{Guo,Bigoni,Hesthaven}, etc.
Recently, Zhu and Qiu \cite{ZhuJSC,ZhuSIAMJSC,ZhuJCP} extended new finite volume WENO schemes, where the linear weights are set to positive numbers artificially with their summation being one, from the structured mesh to two-dimensional triangular mesh and subsequently three-dimensional tetrahedral mesh.
The finite volume WENO schemes with adaptive order of accuracy for the unstructured mesh were designed in \cite{Balsara}.
Li, Shu and Qiu \cite{LiShuQiu} developed the finite volume multi-resolution Hermite WENO schemes on the structured mesh with both the function and its first-order derivative values.
The well-balanced finite volume WENO scheme for the Euler equations with a gravitational source term was proposed in \cite{LiWangDon}, which rectifies the dependency of the sensitivity parameter to reach a hydrostatic equilibrium state.
Chen and Wu \cite{Chen} constructed the third-order finite volume WENO scheme for special relativistic hydrodynamics on the unstructured triangular mesh, which preserves the positivity of the pressure and rest-mass density as well as the subluminal constraint on the fluid velocity.
In \cite{ZhangXiaXu}, the structure-preserving finite volume WENO schemes was presented for the shallow water equations in the arbitrary Lagrangian-Eulerian formulation.

In the finite volume framework, the reconstruction of grid values at cell boundaries from the cell averages is of essential importance to the approximations.
The WENO techniques of nonlinear weights from the finite difference version for the numerical flux can be translated to the finite volume version for the reconstruction.
Jiang and Shu \cite{Jiang} proposed the fifth-order WENO scheme (WENO-JS) in the finite difference formulation with the smoothness indicators, given as the sum of the squares of scaled $L^2$ norms for all the derivatives of the interpolation over the desired cells.
However, the nonlinear weights $\omega^\JS_k$ in \cite{Jiang}, based on the smoothness indicators, cannot satisfy the sufficient condition of fifth-order accuracy at critical points in \cite{Henrick}.
To increase the accuracy of the nonlinear weights, Henrick et al. \cite{Henrick} introduced mapping functions that are applied to the nonlinear weights $\omega^\JS_k$. 
The resulting nonlinear weights $\omega^\M_k$ meet the requirement of fifth-order accuracy and the mapped WENO scheme (WENO-M) yields sharper solutions around discontinuities. 
A different approach to devising nonlinear weights is to introduce the global smoothness indicator and weigh the linear component with the nonlinear component, which gives rise to the Z-type weights $\omega^\Z_k$ \cite{Borges}. 
The WENO-Z scheme leads to less numerical dissipation near discontinuities than WENO-M, while maintaining the ENO behavior.
It is observed that the WENO scheme produces sharper approximations around the discontinuity if the nonlinear weight(s) of the candidate substencil(s) containing the discontinuity is closer to the linear weight(s) within a certain degree. 
In \cite{Gu}, a new set of nonlinear weights denoted by $\omega^\ZR_k$ was designed by taking the $p$th root of the classical smoothness indicators and following the form of the Z-type nonlinear weights.
If $p$ is properly chosen, the nonlinear weights $\omega^\ZR_k$ satisfy the sufficient condition of fifth-order accuracy at critical points and, at the same time, the numerical dissipation is reduced around discontinuities by WENO-ZR.
In the finite volume method, we could apply the above WENO techniques in the finite difference form to the reconstruction of the grid values to yield an approximation to the flux in the conservation form.

To improve the finite volume WENO scheme, we need an detailed analysis to show why the WENO variants perform better than the original WENO-JS, where it is critical to locate specific terms accountable for the numerical dissipation.
In this paper, the WENO-JS, WENO-M, WENO-Z and WENO-ZR schemes are considered to examine the first step of the explicit third-order total variation diminishing (TVD) Runge-Kutta method within the finite volume framework.
The rationale for the improved performance of WENO-M, WENO-Z and WENO-ZR over WENO-JS is identified by investigating the first step of time integration.
Based on this analysis, we propose novel Z-type nonlinear weights for the finite volume WENO scheme, which do not require complicated techniques or sacrifice the versatility of the WENO method.
Instead of taking the difference of the smoothness indicators for the global smoothness indicator, we employ the logarithmic function with tuners to ensure that the numerical dissipation is reduced around discontinuities while the ENO property is preserved.
Besides, the proposed scheme does not necessitate substantial extra computational expenses. 

This paper is composed of the following sections. 
The next section gives a brief review of the WENO-JS, WENO-M, WENO-Z and WENO-ZR schemes, and detailed analyses with those WENO schemes after one third-order TVD Runge-Kutta step for the Riemann problem of the advection equation.
We observe that the WENO scheme, which implements with less numerical dissipation at the first time step, will continue its behavior to the final time.  
In Section \ref{sec:ZL_weight}, we introduce novel Z-type nonlinear weights $\omega^\ZL_k$ and discuss their effects on the related stencils in Section \ref{sec:FV_WENO_discontinuity}.
The two-dimensional finite volume WENO scheme is constructed in Section \ref{sec:FV_WENO_2d}.
Section \ref{sec:nr} contains several one- and two-dimensional numerical experiments for validating the new proposed scheme. 
Conclusions are presented in Section \ref{sec:conclusions}.

\section{Finite volume WENO schemes for 1D scalar conservation laws} \label{sec:FV_WENO_1d}
\subsection{WENO reconstruction from cell averages} \label{sec:WENO_reconstruction}
Suppose that we are given a one-dimensional grid with $N$ cells,
$$ 
   a = x_{1/2} < x_{3/2} < \cdots < x_{N-1/2} < x_{N+1/2} = b.
$$
For each cell $I_i = [x_{i-1/2}, \, x_{i+1/2}]$, the cell center $x_i$ and cell size $\Delta x_i$ are 
$$
   x_i = \frac{1}{2} (x_{i-1/2}+ x_{i+1/2}),~ \Delta x_i = x_{i+1/2} - x_{i-1/2},~ i = 1,2,\cdots,N.
$$ 
The maximum cell size is denoted as $\Delta x = \max_{1 \leqslant i \leqslant N} \Delta x_i$.
The cell average $\vbar_i$ of a function $v(x)$ in the $i$th cell is given as 
$$ 
   \vbar_i = \frac{1}{\Delta x_i} \int^{x_{i+1/2}}_{x_{i-1/2}} v(\xi) \, \dx \xi, \quad i = 1, \cdots, N.
$$

The classical finite volume ENO method seeks for a function $p(x)$ such that it is a $k$th-order accurate approximation to $v(x)$ inside $I_i$, that is, 
$$
   p(x) = v(x) + O(\Delta x^k), \quad x \in I_i,  
$$
and consequently the cell boundary values of $v(x)$ in $I_i$ are approximated by 
$$
   v^-_{i+1/2} = p(x_{i+1/2}), \quad  v^+_{i-1/2} = p(x_{i-1/2}), 
$$
which are at least $k$th-order accurate. 
Here the superscripts $-$ and $+$ denote the left- and right-hand limits, respectively.
For the $k$th-order reconstruction, the stencil is based on $r$ cells to the left and $s$ cells to the right including $I_i$ so that $r + s + 1 = k$.
Define $S_s(i)$ as the stencil composed of these $k$ cells, 
$$
   S_s(i) = \left\{ I_{i-r}, \cdots, I_{i+s} \right\}, \quad s=0, \cdots, k-1,
$$
and a primitive function $V(x)$ of $v(x)$ such that 
$$
   V(x) = \int_a^x v(\xi) d\xi, 
$$
where the lower limit $a$ can be chosen as $a = x_{i-k+1/2}$ \cite{Shu}.  
Then $V(x_{i-k+1/2}) = 0$ and, for $l \geqslant i-k+1$, 
$$
   V(x_{l+1/2}) = \int_{x_{i-k+1/2}}^{x_{l+1/2}} v(\xi) \, \dx \xi 
                = \sum_{j=i-k+1}^l \int_{x_{j-1/2}}^{x_{j+1/2}} v(\xi) \, \dx \xi 
                = \sum_{j=i-k+1}^l \Delta x_j \vbar_j. 
$$
There is a unique polynomial $P(x)$ of degree at most $k$, which interpolates $V(x_{i-r-1/2}), \cdots, V(x_{i+s+1/2})$ and satisfies
$$
   P(x) = V(x) + O(\Delta x^{k+1}).      
$$
Then $P'(x)$ is the approximation to $v(x)$ that we look for, and we denote this polynomial by $p(x)$, which gives the $k$th-order accurate approximation,
\begin{equation} \label{eq:reconstruction_minus_order}
 v^-_{i+1/2} = p(x_{i+1/2}) = v(x_{i+1/2}) + O(\Delta x^k).     
\end{equation}

\begin{figure}[htbp]
\centering
\includegraphics[width=\textwidth]{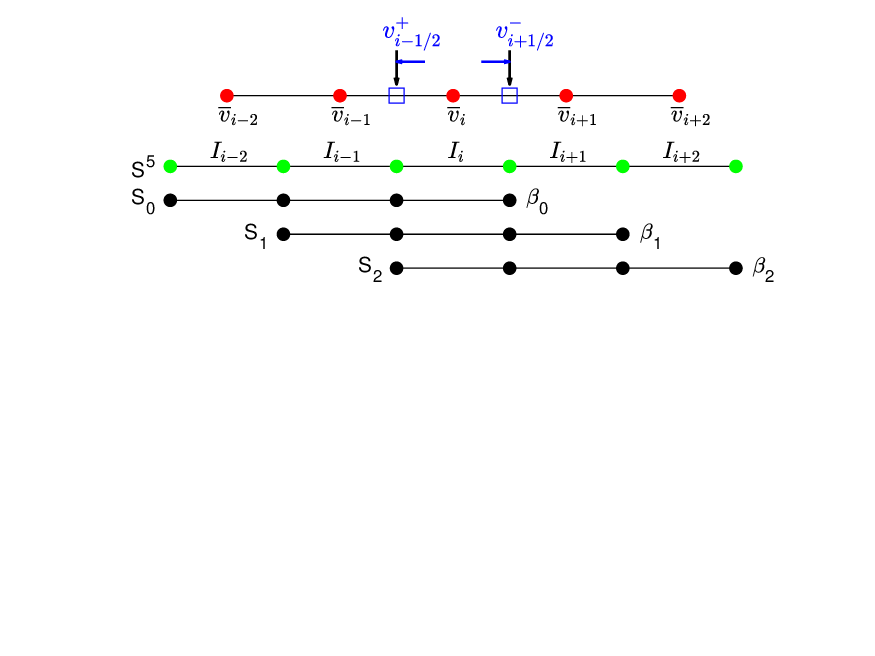}
\vspace{-2.2in}
\caption{The WENO approximation of $v(x)$ at the cell boundaries $v^-_{i+1/2}$ and $v^+_{i-1/2}$ depends on the cell averages over the stencil $S^5 = \{ I_{i-2}, \cdots, I_{i+2} \}$, as well as the substencils $S_0, S_1, S_2$.}
\label{fig:stencil}
\end{figure}

The $k$th-order ENO approximation can be upgraded to the $(2k\! -\! 1)$th-order WENO approximation.
For example, we let $k=3$ for simplicity, which is the fifth-order WENO procedure of Shu \cite{Shu} with $\Delta x = \Delta x_i$.
As shown in Figure \ref{fig:stencil}, the third-order ENO approximation in cell $I_i$ could be one of the three quadratic polynomials $p^{(0)}(x)$, $p^{(1)}(x)$ and $p^{(2)}(x)$, which are the reconstruction polynomials over stencils $S_0 = \{ I_{i-2}, I_{i-1}, I_i \}$, $S_1 = \{ I_{i-1}, I_i, I_{i+1} \}$ and $S_2 = \{ I_i, I_{i+1}, I_{i+2} \}$, respectively.
Evaluating $p^{(s)}(x)$, $s=0,1,2$ at $x = x_{i+1/2}$, gives the third-order accurate approximations as
\begin{equation} \label{eq:3rd_reconstruction_minus}
\begin{aligned}
 v^{(0)-}_{i+1/2} &= p^{(0)}(x_{i+1/2}) =   \frac{1}{3} \vbar_{i-2} - \frac{7}{6} \vbar_{i-1} + \frac{11}{6} \vbar_i, \\
 v^{(1)-}_{i+1/2} &= p^{(1)}(x_{i+1/2}) = - \frac{1}{6} \vbar_{i-1} + \frac{5}{6} \vbar_i     + \frac{1}{3}  \vbar_{i+1}, \\
 v^{(2)-}_{i+1/2} &= p^{(2)}(x_{i+1/2}) =   \frac{1}{3} \vbar_i     + \frac{5}{6} \vbar_{i+1} - \frac{1}{6}  \vbar_{i+2},
\end{aligned}
\end{equation}
and according to \eqref{eq:reconstruction_minus_order}, we have 
\begin{equation} \label{eq:3rd_reconstruction_minus_order}
 v^{(s)-}_{i+1/2} = v(x_{i+1/2}) + O(\Delta x^3), \quad s = 0, 1, 2.   
\end{equation}
Similarly, the fifth-order ENO reconstruction in the big stencil $S^5$ yields 
\begin{equation} \label{eq:5rd_reconstruction_minus}
 v^{(S^5)-}_{i+1/2} = \frac{1}{30} \vbar_{i-2} - \frac{13}{60} \vbar_{i-1} + \frac{47}{60} \vbar_i + \frac{9}{20} \vbar_{i+1} - \frac{1}{20} \vbar_{i+2},
\end{equation}
and by \eqref{eq:reconstruction_minus_order},
\begin{equation} \label{eq:5rd_reconstruction_minus_order}
 v^{(S^5)-}_{i+1/2} = v(x_{i+1/2}) + O(\Delta x^5).      
\end{equation}
Comparing \eqref{eq:3rd_reconstruction_minus} with \eqref{eq:5rd_reconstruction_minus} gives 
\begin{equation} \label{eq:reconstruction_minus_linear_weights}
 v^{(S^5)-}_{i+1/2} = \sum_{s=0}^2 d_s v^{(s)-}_{i+1/2},
\end{equation}
with $d_0 = \frac{1}{10}, d_1 = \frac{3}{5}$ and $d_2 = \frac{3}{10}$ known as the linear weights. 
Notice that $\sum_{s=0}^2 d_s = 1$ for consistency.
The WENO reconstruction is formulated as a convex combination of the three different third-order approximations in \eqref{eq:3rd_reconstruction_minus} as
\begin{equation} \label{eq:weno_appr}
 v^-_{i+1/2} = \sum_{s=0}^2 \omega_s v^{(s)-}_{i+1/2},
\end{equation}
where the nonlinear weights $\omega_s$ satisfy 
$$
   \omega_s \geqslant 0, \quad \sum_{s=0}^2 \omega_s = 1.
$$
If $v(x)$ is smooth on the stencil $S^5$, we would like to have the fifth-order approximation as \eqref{eq:5rd_reconstruction_minus_order},
$$
   v^-_{i+1/2} = v(x_{i+1/2}) + O(\Delta x^5).
$$
From \eqref{eq:3rd_reconstruction_minus_order}, \eqref{eq:5rd_reconstruction_minus_order}, \eqref{eq:reconstruction_minus_linear_weights} and \eqref{eq:weno_appr}, we have
\begin{align*}
 v^-_{i+1/2} &= \sum_{s=0}^2 d_s v^{(s)-}_{i+1/2} + \sum_{s=0}^2 ( \omega_s - d_s) v^{(s)-}_{i+1/2} \\
             &= v^{(S^5)-}_{i+1/2} + \sum_{s=0}^2 ( \omega_s - d_s ) \left[ v(x_{i+1/2}) + O(\Delta x^3) \right] \\
             &= v(x_{i+1/2}) + O(\Delta x^5) + \sum_{s=0}^2 ( \omega_s - d_s ) v(x_{i+1/2}) + \sum_{s=0}^2 ( \omega_s - d_s ) O(\Delta x^3) \\
             &= v(x_{i+1/2}) + O(\Delta x^5) + \sum_{s=0}^2 ( \omega_s - d_s ) O(\Delta x^3).
\end{align*}
It follows that 
\begin{equation} \label{eq:WENO_condition}
 \omega_s - d_s = O(\Delta x^2)
\end{equation} 
is required to attain fifth-order accuracy in smooth regions.

Now, if $v(x)$ has a discontinuity in one or more of the stencils, we would like to make the corresponding weight(s) $\omega_s$ essentially vanish. 
This could be achieved by the so-called smoothness indicators \cite{Jiang},
$$
   \beta_s = \sum_{l=1}^2 \Delta x^{2l-1} \int_{x_{i-1/2}}^{x_{i+1/2}} \left( \frac{\dx^l}{\dx x^l} p^{(s)}(x) \right)^2 \dx x, \quad s = 0,1,2,
$$
more specifically, 
\begin{equation} \label{eq:smooth_indicator}
\begin{aligned}
 \beta_0 &= \frac{13}{12} \left( \vbar_{i-2} - 2 \vbar_{i-1} + \vbar_i \right)^2 + \frac{1}{4} \left( \vbar_{i-2} - 4 \vbar_{i-1} + 3 \vbar_i \right)^2, \\
 \beta_1 &= \frac{13}{12} \left( \vbar_{i-1} - 2 \vbar_i + \vbar_{i+1} \right)^2 + \frac{1}{4} \left( \vbar_{i-1} - \vbar_{i+1} \right)^2, \\
 \beta_2 &= \frac{13}{12} \left( \vbar_i - 2 \vbar_{i+1} + \vbar_{i+2} \right)^2 + \frac{1}{4} \left( 3 \vbar_i - 4 \vbar_{i+1} + \vbar_{i+2} \right)^2.
\end{aligned}
\end{equation}
Taylor series expansion of \eqref{eq:smooth_indicator} at $x=x_i$ gives
\begin{align*}
 \beta_0 &= v'^2_i \Delta x^2 + \left( \frac{13}{12} v''^2_i - \frac{7}{12} v'_i v'''_i \right) \Delta x^4 - \left( \frac{13}{6} v''_i v'''_i - \frac{1}{2} v'_i v^{(4)}_i \right) \Delta x^5 + O(\Delta x^6), \\
 \beta_1 &= v'^2_i \Delta x^2 + \left( \frac{13}{12} v''^2_i + \frac{5}{12} v'_i v'''_i \right) \Delta x^4 + O(\Delta x^6), \\
 \beta_2 &= v'^2_i \Delta x^2 + \left( \frac{13}{12} v''^2_i - \frac{7}{12} v'_i v'''_i \right) \Delta x^4 + \left( \frac{13}{6} v''_i v'''_i - \frac{1}{2} v'_i v^{(4)}_i \right) \Delta x^5 + O(\Delta x^6),
\end{align*}
where $v'_i$, $v''_i$ and $v'''_i$ are the first, second and third derivatives of $v(x)$ at $x=x_i$, respectively.
In \cite{Jiang}, the nonlinear weights $\omega^\JS_k$ in \eqref{eq:weno_appr} were given by
\begin{equation} \label{eq:weights_JS}
 \omega^\JS_s = \frac{\alpha_s}{\sum_{r=0}^2 \alpha_r},~ \alpha_s = \frac{d_s}{(\beta_s + \epsilon)^2},~ s=0,1,2.
\end{equation} 
Here $\epsilon$ is a small positive number added to prevent the denominator from becoming zero, and is chosen as $\epsilon = 10^{-6}$ in \cite{Jiang}.
The Taylor expansion yields $\omega^\JS_s = d_s + O(\Delta x^2)$ if $v'_i \ne 0$, while $\omega^\JS_s = d_s + O(\Delta x)$ if $v'_i = 0$. 
Thus, $\omega^\JS_s$ does not satisfy the condition \eqref{eq:WENO_condition} at the critical points. 
In \cite{Henrick}, Henrick et al. introduced the mapping functions 
$$
   g_s(\omega) = \frac{ (d_s + d_s^2 - 3 d_s \omega + \omega^2)\omega}{d_s^2 + (1 - 2d_s)\omega}, \quad s=0,1,2,
$$
satisfying the following conditions:
\begin{itemize}
\item every function $g_s$ is monotonically increasing on $[0,1]$ with a finite slope;
\item $g_s(0) = 0$, $g_s(1) = 1$, $g_s(d_s) = d_s$;
\item $g_s(\omega) > \omega$ on $(0, d_s)$ and $g_s(\omega) < \omega$ on $(d_s, 1)$, 
\end{itemize}
and defined the mapped nonlinear weights $\omega^\M_s$ as
\begin{equation} \label{eq:weights_M}
 \omega^\M_s = \frac{\alpha^*_s}{\sum_{r=0}^2 \alpha^*_r},~ \alpha^*_s = g_s(\omega^\JS_s),~ s=0,1,2,
\end{equation}
with $\omega^\JS_s$ in \eqref{eq:weights_JS} and $\epsilon = 10^{-40}$.
Taylor series expansion shows that the weights $\omega^\M_s$ satisfy $\omega^\M_s = d_s + O(\Delta x^3)$ regardless of the value of $v'_i$, which agrees with the condition \eqref{eq:WENO_condition}.
In \cite{Borges}, the Z-type nonlinear weights $\omega^\Z_s$ were defined as 
\begin{equation} \label{eq:weights_Z}
 \omega^\Z_s = \frac{\alpha_s}{\sum_{r=0}^2 \alpha_r},~ \alpha_s = d_s \left( 1 + \frac{\tau_5}{\beta_s + \epsilon} \right),~ s=0,1,2,
\end{equation}
with the global smoothness indicator $\tau_5 = |\beta_0 - \beta_2|$ and $\epsilon = 10^{-40}$.
It can be shown that $\omega^\Z_s = d_s + O(\Delta x^3)$ for $v'_i \ne 0$ whereas $\omega^\Z_s = d_s + O(\Delta x)$ for $v'_i = 0$, implying that the weights $\omega^\Z_s$ satisfy the condition \eqref{eq:WENO_condition} only at non-critical points.
In order to achieve fifth-order accuracy at critical points, we introduced the new Z-type nonlinear weights $\omega^\ZR_s$ \cite{Gu},
\begin{equation} \label{eq:weights_ZR}
 \omega^\ZR_s = \frac{\alpha_s}{\sum_{r=0}^2 \alpha_r},~ \alpha_s = d_s \left( 1 + \left( \frac{\tau^\ZR_5}{\sqrt[p]{\beta_s} + \epsilon} \right)^p \right),~ s=0,1,2,
\end{equation}
with the global smoothness indicator $\tau^\ZR_5 = | \sqrt[p]{\beta_0} - \sqrt[p]{\beta_2} |$.
Note that for $p=1$, the nonlinear weights $\omega^\ZR_s$ in \eqref{eq:weights_ZR} coincide with $\omega^\Z_s$ in \eqref{eq:weights_Z}. 
It was shown in \cite{Gu} that $\omega^\ZR_s = d_s + O(\Delta x^{3p})$ if $v'_i \ne 0$, while $\omega^\ZR_s = d_s + O(\Delta x^p)$ if $v'_i = 0$.
By \eqref{eq:WENO_condition} we see that the nonlinear weights with $p \geqslant 2$ reach fifth-order accuracy if there is a first-order critical point at $x = x_i$.

\subsection{Finite volume WENO scheme for 1D scalar case} \label{sec:WENO_schemes}
Consider the one-dimensional scalar conservation law:
\begin{equation} \label{eq:1d_scalar_hyperbolic}
 \frac{\partial u}{\partial t} + \frac{\partial f(u)}{\partial x} = 0,
\end{equation}
with $u = u(x,t)$.
We integrate \eqref{eq:1d_scalar_hyperbolic} over the interval $I_i$ to obtain
\begin{equation} \label{eq:1d_scalar_hyperbolic_integral}
 \frac{\dx \ubar(x_i,t)}{\dx t}  = - \frac{f(u(x_{i+1/2},t)) - f(u(x_{i-1/2},t))}{\Delta x_i},
\end{equation}
with the cell average
$$
   \ubar(x_i,t) = \frac{1}{\Delta x_i} \int_{x_{i-1/2}}^{x_{i+1/2}} u(\xi, t) \, \dx \xi.
$$
We approximate \eqref{eq:1d_scalar_hyperbolic_integral} by the following conservative method of the form
\begin{equation} \label{eq:1d_scalar_hyperbolic_FV}
 \frac{\dx \ubar_i(t)}{\dx t}  = - \frac{\hat{f}_{i+1/2} - \hat{f}_{i-1/2}}{\Delta x_i},
\end{equation}
where $\ubar_i(t)$ is the numerical approximation to $\ubar(x_i,t)$, and the numerical flux $\hat{f}_{i+1/2}$ is given by
\begin{equation} \label{eq:1d_scalar_hyperbolic_FV_flux}
   \hat{f}_{i+1/2} = h(u^-_{i+1/2}, u^+_{i+1/2}),
\end{equation}
with $u^-_{i+1/2}$ from the WENO reconstruction in \eqref{sec:WENO_reconstruction} and $u^+_{i+1/2}$ in the reverse order.
For the numerical flux, we use the Lax-Friedrichs flux $h$,
\begin{equation} \label{eq:LF_flux}
 h(a,b) = {1 \over 2} [f(a) + f(b) - \alpha (b - a)], 
\end{equation}
with $\alpha = \max_{u} |f'(u)|$ being a constant taken over the relevant range of $u$. 

\subsection{Spatial fourth-order accuracy in smooth regions} \label{sec:1d_FV_WENO_smooth}
When the condition \eqref{eq:WENO_condition} is satisfied in smooth regions, we have $u^\pm_{i+1/2} = u(x_{i+1/2}) + O(\Delta x^5)$.
Taylor series expansion yields
\begin{align*}
 f(u^{\pm}_{i+1/2}) &= f \left( u(x_{i+1/2}) \right) + f' \left( u(x_{i+1/2}) \right) \cdot O(\Delta x^5) + O(\Delta x^{10}) \\
                    &= f \left( u(x_{i+1/2}) \right) + O(\Delta x^5),
\end{align*}
where the time variable $t$ is dropped for the notational convenience as the time is fixed. 
Then
\begin{align*}
 \hat{f}_{i+1/2} &= {1 \over 2} [f(u^-_{i+1/2}) + f(u^+_{i+1/2}) - \alpha (u^+_{i+1/2} - u^-_{i+1/2})] \\
                 &= {1 \over 2} [2 f \left( u(x_{i+1/2}) \right) + O(\Delta x^5) - \alpha \left( u(x_{i+1/2})-u(x_{i+1/2}) + O(\Delta x^5) \right)] \\
                 & = f \left( u(x_{i+1/2}) \right) + O(\Delta x^5).
\end{align*}
Similarly,
$$
   \hat{f}_{i-1/2} = f \left( u(x_{i-1/2}) \right) + O(\Delta x^5).
$$
From these equations we obtain
\begin{align*}
 \frac{\hat{f}_{i+1/2} - \hat{f}_{i-1/2}}{\Delta x_i} = {} & \frac{f \left( u(x_{i+1/2}) \right) + O(\Delta x^5) - f \left( u(x_{i-1/2}) \right) - O(\Delta x^5)}{\Delta x_i} \\
                                                      = {} & \frac{f \left( u(x_{i+1/2}) \right) - f \left( u(x_{i-1/2}) \right)}{\Delta x_i} + O(\Delta x^4),
\end{align*} 
and hence fourth-order accuracy is achieved.

\subsection{Dissipation around discontinuities} \label{sec:FV_WENO_discontinuity}
We now analyze the behavior of the numerical solutions by the WENO-JS, WENO-M, WENO-Z and WENO-ZR schemes introduced in Sections \ref{sec:WENO_reconstruction} and \ref{sec:WENO_schemes} after a single time step using the third-order TVD Runge-Kutta method.
Our discussion will focus on the Riemann problem of the advection equation. 
In every stage, we start with the computation of nonlinear weights based on the cell averages, which leads to the construction of the numerical flux. 
The approximations at this stage are computed by the conservative form \eqref{eq:1d_scalar_hyperbolic_FV} with the numerical fluxes. 
Through this analysis, we can see the effect of each WENO scheme on the numerical solutions.

The Riemann problem of the advection equation has the form
\begin{equation} \label{eq:weight_comparison}
\begin{aligned}
 u_t + u_x &= 0, \\
    u(x,0) &= \left\{ 
               \begin{array}{ll} 
                u_L, & x < 0, \\ 
                u_R, & x \geqslant 0,
               \end{array} 
              \right.
\end{aligned}
\end{equation}
where the initial condition is a step function containing a single jump discontinuity at $x=0$ with the difference $\delta = u_L - u_R$. 
Let $I_0 = [x_{-1/2}, \, x_{1/2}] = [0, \Delta x]$ and $x_0 = \Delta x/2$.
Then the rest of the cells are defined accordingly.
At $t=0$, the jump discontinuity lies at the boundary between the cells $I_{-1}$ and $I_0$, and 
$$
   \ubar^0_j = \left\{ 
               \begin{array}{ll} 
                u_L, & j < 0, \\ 
                u_R, & j \geqslant 0.
               \end{array} 
              \right.
$$
As the CFL condition requires that the Courant number $\nu = \frac{\Delta t}{\Delta x} \leqslant 1$, we set $0 < \nu \leqslant 0.5$.
At $t = \Delta t$, the discontinuity moves to $x = \Delta t$ in the cell $I_0$.
It is easy to compute the cell average over the cell $I_0$ after the time step $\Delta t$,
$$
   \ubar(x_0,\Delta t) = \frac{1}{\Delta x} \int_0^{\Delta x} u(\xi, \Delta t) \, \dx \xi 
                       = \frac{1}{\Delta x} \left[ u_L \Delta t + u_R (\Delta x - \Delta t) \right]
                       = u_R + \nu \delta. 
$$
So we have 
\begin{equation} \label{eq:exact_sol}
   \ubar^1_j = \left\{ 
               \begin{array}{ll} 
                u_L, & j<0, \\ 
                u_R + \nu \delta, & j=0, \\
                u_R, & j>0.
               \end{array} 
              \right.
\end{equation} 
Set $\ceps = 10^{-12}$ for WENO-JS and $\teps = 10^{-40}$ for the rest WENO schemes.
We assume that $p \geqslant 1$, $\delta > 0$ and $\sqrt[p]{\delta} \gg \ceps$.

To begin, we list some formulas that will be used in the following weight analysis, based on \eqref{eq:weights_JS}, \eqref{eq:weights_Z} and \eqref{eq:weights_ZR}.
When $\beta^\N_0 = \beta^\N_1 = 0$, $\N \in \{ \JS, \M, \Z, \ZR \}$, we have 
\begin{subequations} \label{eq:weights0}
 \begin{equation} \label{eq:nw_JS0}
  \begin{aligned}
  \cbeta^\JS_2 &= \left( \frac{\beta^\JS_2}{\ceps} + 1 \right)^2, \\
  \omega^\JS_0 &= \frac{d_0}{d_0 + d_1 + {d_2 \over \cbeta^\JS_2}},~           
  \omega^\JS_1  = \frac{d_1}{d_0 + d_1 + {d_2 \over \cbeta^\JS_2}},~            
  \omega^\JS_2  = \frac{d_2}{d_0 \cbeta^\JS_2 + d_1 \cbeta^\JS_2 + d_2};   
  \end{aligned}  
 \end{equation}    
 \begin{equation} \label{eq:nw_JSM0}
  \begin{aligned}
  \tbeta^\M_2   &= \left( \frac{\beta^\M_2}{\teps} + 1 \right)^2, \\
  \tomega^\JS_0 &= \frac{d_0}{d_0 + d_1 + {d_2 \over \tbeta^\M_2}},~          
  \tomega^\JS_1  = \frac{d_1}{d_0 + d_1 + {d_2 \over \tbeta^\M_2}},~         
  \tomega^\JS_2  = \frac{d_2}{d_0 \tbeta^\M_2 + d_1 \tbeta^\M_2 + d_2};   
  \end{aligned}
 \end{equation}
 \begin{equation} \label{eq:nw_Z0}
  \begin{aligned}
  \tbeta^\Z_2 &= \frac{1 + \frac{\beta^\Z_2}{\teps}}{1 + \frac{\beta^\Z_2}{\beta^\Z_2+\teps}}, \\
  \omega^\Z_0 &= \frac{d_0}{d_0 + d_1 + {d_2 \over \tbeta^\Z_2}},~
  \omega^\Z_1  = \frac{d_1}{d_0 + d_1 + {d_2 \over \tbeta^\Z_2}},~
  \omega^\Z_2  = \frac{d_2}{d_0 \tbeta^\Z_2 + d_1 \tbeta^\Z_2 + d_2};  
  \end{aligned}
 \end{equation}
 \begin{equation} \label{eq:nw_ZR0}
  \begin{aligned}
  \tbeta^\ZR_2 &= \frac{1 + \frac{\beta^\ZR_2}{\teps^p}}{1 + \frac{\beta^\ZR_2}{\left( \sqrt[p]{\beta^\ZR_2}+\teps \right)^p}}, \\
  \omega^\ZR_0 &= \frac{d_0}{d_0 + d_1 + {d_2 \over \tbeta^\ZR_2}},~
  \omega^\ZR_1  = \frac{d_1}{d_0 + d_1 + {d_2 \over \tbeta^\ZR_2}},~
  \omega^\ZR_2  = \frac{d_2}{d_0 \tbeta^\ZR_2 + d_1 \tbeta^\ZR_2 + d_2}.  
  \end{aligned}
 \end{equation}
\end{subequations}
Otherwise, we have
\begin{subequations} \label{eq:weights1}
 \begin{equation} \label{eq:nw_JS1}
  \begin{aligned}
  \cbeta^\JS_1 &= \left( \frac{\beta^\JS_1+\ceps}{\beta^\JS_0+\ceps} \right)^2,~\cbeta^\JS_2 = \left( \frac{\beta^\JS_2+\ceps}{\beta^\JS_0+\ceps} \right)^2, \\
  \omega^\JS_0 &= \frac{d_0}{d_0 + {d_1 \over \cbeta^\JS_1} + {d_2 \over \cbeta^\JS_2}},~
  \omega^\JS_1  = \frac{d_1}{d_0 \cbeta^\JS_1 + d_1 +  d_2 {\cbeta^\JS_1 \over \cbeta^\JS_2}},~
  \omega^\JS_2  = \frac{d_2}{d_0 \cbeta^\JS_2 + d_1{ \cbeta^\JS_2 \over \cbeta^\JS_1} + d_2};  
  \end{aligned}
 \end{equation}    
 \begin{equation} \label{eq:nw_JSM1}
  \begin{aligned}
  \tbeta^\M_1   &= \left( \frac{\beta^\M_1+\teps}{\beta^\M_0+\teps} \right)^2,~\tbeta^\M_2 = \left( \frac{\beta^\M_2+\teps}{\beta^\M_0+\teps} \right)^2, \\
  \tomega^\JS_0 &= \frac{d_0}{d_0 + { d_1 \over \tbeta^\M_1} + { d_2 \over \tbeta^\M_2}},~
  \tomega^\JS_1  = \frac{d_1}{d_0 \tbeta^\M_1 + d_1 + d_2 {\tbeta^\M_1 \over \tbeta^\M_2}},~
  \tomega^\JS_2  = \frac{d_2}{d_0 \tbeta^\M_2 + d_1 {\tbeta^\M_2 \over \tbeta^\M_1} + d_2};  
  \end{aligned}
 \end{equation}
 \begin{equation} \label{eq:nw_Z1}
  \begin{aligned}
  \tbeta^\Z_1 &= \frac{1+\frac{\tau^\Z_5}{\beta^\Z_0+\teps}}{1+\frac{\tau^\Z_5}{\beta^\Z_1+\teps}},~\tbeta^\Z_2 = \frac{1+\frac{\tau^\Z_5}{\beta^\Z_0+\teps}}{1+\frac{\tau^\Z_5}{\beta^\Z_2+\teps}}, \\
  \omega^\Z_0 &= \frac{d_0}{d_0 + {d_1 \over \tbeta^\Z_1} + {d_2 \over \tbeta^\Z_2}},~
  \omega^\Z_1  = \frac{d_1}{d_0 \tbeta^\Z_1 + d_1 + d_2{ \tbeta^\Z_1 \over \tbeta^\Z_2}},~
  \omega^\Z_2  = \frac{d_2}{d_0 \tbeta^\Z_2 + d_1 {\tbeta^\Z_2 \over \tbeta^\Z_1} + d_2};  
  \end{aligned}
 \end{equation}
 \begin{equation} \label{eq:nw_ZR1}
  \begin{aligned}
  \tbeta^\ZR_1 &= \frac{1+\left( \frac{\tau^\ZR_5}{\sqrt[p]{\beta^\ZR_0}+\teps} \right)^p}{1+\left( \frac{\tau^\ZR_5(p)}{\sqrt[p]{\beta^\ZR_1}+\teps} \right)^p},~\tbeta^\ZR_2 = \frac{1+\left( \frac{\tau^\ZR_5}{\sqrt[p]{\beta^\ZR_0}+\teps} \right)^p}{1+\left( \frac{\tau^\ZR_5}{\sqrt[p]{\beta^\ZR_2}+\teps} \right)^p}, \\
  \omega^\ZR_0 &= \frac{d_0}{d_0 + {d_1 \over \tbeta^\ZR_1} + {d_2 \over \tbeta^\ZR_2}},~
  \omega^\ZR_1  = \frac{d_1}{d_0 \tbeta^\ZR_1 + d_1 + d_2 {\tbeta^\ZR_1 \over \tbeta^\ZR_2}},~
  \omega^\ZR_2  = \frac{d_2}{d_0 \tbeta^\ZR_2 + d_1 {\tbeta^\ZR_2 \over \tbeta^\ZR_1} + d_2}.  
  \end{aligned}
 \end{equation}
\end{subequations}

\subsubsection{First stage of the TVD Runge-Kutta method} \label{sec:Euler}
Consider the first stage of the optimal third-order TVD Runge-Kutta method \cite{Gottlieb, ShuOsherI}, which is simply the Euler method,
\begin{equation} \label{eq:RK1st} 
 u^{(1)} = u^n + \Delta t L(u^n),
\end{equation}
where $L$ is the numerical spatial operator approximating the spatial derivative in \eqref{eq:weight_comparison}. 
As the discontinuity lies at $x=0$, we deal with the following 6 stencils:
\begin{enumerate}[label=\arabic*., font=\itshape]
\item $\{ I_{j-4}, I_{j-3}, I_{j-2}, I_{j-1}, I_j \},~j \leqslant -1$;
\item $\{ I_{-4}, I_{-3}, I_{-2}, I_{-1}, I_{0} \}$; \label{item:stt1}
\item $\{ I_{-3}, I_{-2}, I_{-1}, I_{0}, I_{1} \}$;
\item $\{ I_{-2}, I_{-1}, I_{0}, I_{1}, I_{2} \}$;
\item $\{ I_{-1}, I_{0}, I_{1}, I_{2}, I_{3} \}$;    \label{item:end1}
\item $\{ I_j, I_{j+1}, I_{j+2}, I_{j+3}, I_{j+4} \},~j \geqslant 0$,
\end{enumerate}
where the second to fifth stencils contain the discontinuity.

The superscript $\N \in \{ \JS, \M, \Z, \ZR \}$ for the smoothness indicator $\beta^\N_s$ is dropped for all the WENO schemes share the same smoothness indicators in this stage.

\begin{stencil} \label{case:1}
$\{ I_{j-4}, I_{j-3}, I_{j-2}, I_{j-1}, I_j \},~j \leqslant -1$ \\
The smoothness indicators satisfy $\beta_s = 0$ for $s = 0,1,2$. 
From \eqref{eq:weights_JS}, \eqref{eq:weights_Z} and \eqref{eq:weights_ZR}, it is easy to have
$$
   \omega^\N_s = d_s,~ \N \in \{ \JS, \Z, \ZR \},~ s = 0,1,2,
$$
which corresponds to $x=-0.03, -0.02$ in Table \ref{tab:weight_comparison_first}.
However, the WENO-M nonlinear weights $\omega^\M_0$ and $\omega^\M_2$ slightly deviate from the respective linear weights $d_0$ and $d_2$ due to the use of the mapping functions with the double precision calculations.
\end{stencil}

\begin{stencil} \label{case:2}
$\{ I_{-4}, I_{-3}, I_{-2}, I_{-1}, I_{0} \}$ \\
The smoothness indicators are $\beta_0 = \beta_1 = 0$ and $\beta_2 = \frac{4}{3} \delta^2$. 
Using \eqref{eq:weights0}, the properties of $g_s(x)$ and the fact that $d_0 < \omega^\JS_0, \tomega^\JS_0 \lessapprox \frac{1}{7}$, we have   
$$
   d_0 < \alpha^*_0 = g_0 \left( \tomega^\JS_0 \right) \lessapprox g_0 \left( \frac{1}{7} \right) = \frac{143}{1421} < \omega^\JS_0,
$$
and similarly, 
$$
   d_1 < \alpha^*_1 = g_1 \left( \tomega^\JS_1 \right) \lessapprox g_1 \left( \frac{6}{7} \right) = \frac{372}{539} < \omega^\JS_1,
$$
since $d_1 < \omega^\JS_1, \tomega^\JS_1 \lessapprox \frac{6}{7}$.
As $\tomega^\JS_2 \gtrapprox 0$, we have $\alpha^*_2 = g_2 \left( \tomega^\JS_2 \right) \gtrapprox 0$.
By \eqref{eq:weights_M}, 
$$ 
   \omega^\M_0 \approx \frac{1573}{12361},~ \omega^\M_1 \approx \frac{10788}{12361},~ \omega^\M_2 \approx 0. 
$$
Furthermore, with $0 \lessapprox \omega^\JS_2, \tomega^\JS_2 \ll d_2$, we get 
$$
   \tomega^\JS_2 < \alpha^*_2 = g_2 \left( \tomega^\JS_2 \right) 
 = \tomega^\JS_2 \frac{d_2 + d_2^2 - 3 d_2 \tomega^\JS_2 + \left( \tomega^\JS_2 \right)^2}{d_2^2 + \tomega^\JS_2 (1 - 2d_2)}
 < \tomega^\JS_2 \frac{d_2 + d_2^2 + \left( \tomega^\JS_2 \right)^2}{d_2^2} < 6 \tomega^\JS_2 , 
$$
and
$$
   \frac{\omega^\JS_2}{\tomega^\JS_2} = \frac{d_0 \tbeta_2 + d_1 \tbeta_2 + d_2}{d_0 \cbeta_2 + d_1 \cbeta_2 + d_2}
 = \left( \frac{\ceps}{\teps} \right)^2 \frac{d_0 \left( \beta_2+\teps \right)^2 + d_1 \left( \beta_2+\teps \right)^2 + d_2 \teps^2}{d_0 \left( \beta_2+\ceps \right)^2 + d_1 \left( \beta_2+\ceps \right)^2 + d_2 \ceps^2} \approx 10^{56}, 
$$
which leads to 
$$
   \alpha^*_2 \ll (d_0 + d_1) \omega^\JS_2 < (\alpha^*_0 + \alpha^*_1 + \alpha^*_2) \omega^\JS_2,
$$
and thus
$$
   \omega^\M_2 = \frac{\alpha^*_2}{\alpha^*_0 + \alpha^*_1 + \alpha^*_2} \ll \omega^\JS_2.
$$
Note that $d_0 < \omega^\Z_0, \omega^\ZR_0 \lessapprox \frac{1}{7}$, $d_1 < \omega^\Z_1, \omega^\ZR_1 \lessapprox \frac{6}{7}$ and $\omega^\Z_2, \omega^\ZR_2 \gtrapprox 0$.
Also the estimates
$$
   \cbeta^\JS_2 = \left( \frac{\beta_2}{\ceps} + 1 \right)^2 \approx 10^{24} \beta_2^2,~
   \tbeta^\Z_2  = \frac{1 + \frac{\beta_2}{\teps}}{2 - \frac{\teps}{\beta_2+\teps}} \approx \frac{10^{40}}{2} \beta_2,~
   \tbeta^\ZR_2 = \frac{1 + \frac{\beta_2}{\teps^p}}{1 + \left( 1 - \frac{\teps}{\sqrt[p]{\beta_2}+\teps} \right)^p} \approx \frac{10^{40p}}{2} \beta_2,
$$
give $\cbeta^\JS_2 \ll \tbeta^\Z_2 \ll \tbeta^\ZR_2$, and hence $\omega^\JS_s < \omega^\Z_s < \omega^\ZR_s,~s=0,1$ and $\omega^\ZR_2 \ll \omega^\Z_2 \ll \omega^\JS_2$.
In summary, we have
\begin{gather*}
 \omega^\N_0 \lessapprox \frac{1}{7},~\omega^\N_1 \lessapprox \frac{6}{7},~\omega^\N_2 \gtrapprox 0,~\N \in \{ \JS, \Z, \ZR \}, \\
 \omega^\M_0 \approx \frac{1573}{12361},~ \omega^\M_1 \approx \frac{10788}{12361},~ \omega^\M_2 \gtrapprox 0, \\
 \omega^\M_2 \ll \omega^\JS_2,~ \omega^\ZR_2 \ll \omega^\Z_2 \ll \omega^\JS_2, 
\end{gather*}
corresponding to $x=-0.01$ in Table \ref{tab:weight_comparison_first}.
\end{stencil}

\begin{stencil} \label{case:3}
$\{ I_{-3}, I_{-2}, I_{-1}, I_{0}, I_{1} \}$ \\
Now the smoothness indicators satisfy $\beta_0 = 0$ and $\beta_1 = \frac{4}{3} \delta^2$, $\beta_2 = \frac{10}{3} \delta^2$.
By \eqref{eq:weights1},
we have $0 \lessapprox \omega^\JS_1, \tomega^\JS_1 \ll d_1$, and hence 
$$
   \alpha^*_1 = g_1 \left( \tomega^\JS_1 \right) 
 = \tomega^\JS_1 \frac{\left( d_1 - \tomega^\JS_1 \right)^2 + d_1 \left( 1 - \tomega^\JS_1 \right)}{\left( d_1 - \tomega^\JS_1 \right)^2 + \tomega^\JS_1 \left( 1 - \tomega^\JS_1 \right)} 
 \approx \tomega^\JS_1 \frac{d_1^2 + d_1}{d_1^2} = \frac{8}{3} \tomega^\JS_1,
$$
and
$$
   \frac{\omega^\JS_1}{\tomega^\JS_1} = \frac{d_0 \tbeta_1 + d_1 + d_2 \tbeta_1 / \tbeta_2}{d_0 \cbeta_1 + d_1 + d_2 \cbeta_1 / \cbeta_2}
 = \left( \frac{\ceps}{\teps} \right)^2 \frac{d_0 \left( \beta_1+\teps \right)^2 + d_1 \teps^2 + d_2 \left( \frac{\beta_1+\teps}{\beta_2+\teps} \teps \right)^2}{d_0 \left( \beta_1+\ceps \right)^2 + d_1 \ceps^2 + d_2 \left( \frac{\beta_1+\ceps}{\beta_2+\ceps} \ceps \right)^2} \approx 10^{56}.
$$
Similarly, since $0 \lessapprox \omega^\JS_2, \tomega^\JS_2 \ll d_2$, we have
$$
   \alpha^*_2 = g_2 \left( \tomega^\JS_2 \right) 
 = \tomega^\JS_2 \frac{\left( d_2 - \tomega^\JS_2 \right)^2 + d_2 \left( 1 - \tomega^\JS_2 \right)}{\left( d_2 - \tomega^\JS_2 \right)^2 + \tomega^\JS_2 \left( 1 - \tomega^\JS_2 \right)} 
 \approx \tomega^\JS_2 \frac{d_2^2 + d_2}{d_2^2} = \frac{13}{3} \tomega^\JS_2,
$$
and
$$
   \frac{\omega^\JS_2}{\tomega^\JS_2} = \frac{d_0 \tbeta_2 + d_1 \tbeta_2 / \tbeta_1 + d_2}{d_0 \cbeta_2 + d_1 \cbeta_2 / \cbeta_1 + d_2}
 = \left( \frac{\ceps}{\teps} \right)^2 \frac{d_0 \left( \beta_2+\teps \right)^2 + d_1 \left( \frac{\beta_2+\teps}{\beta_1+\teps} \teps \right)^2 + d_2 \teps^2}{d_0 \left( \beta_2+\ceps \right)^2 + d_1 \left( \frac{\beta_2+\ceps}{\beta_1+\ceps} \ceps \right)^2 + d_2 \ceps^2} \approx 10^{56}.
$$
Also, $\alpha^*_s \ll (\alpha^*_0 + \alpha^*_1 + \alpha^*_2) \omega^\JS_s$ for $s=1,2$ implies that $\omega^\M_s \ll \omega^\JS_s$.
Note that $\omega^\Z_0, \omega^\ZR_0 \lessapprox 1$ and $\omega^\Z_s, \omega^\ZR_s \gtrapprox 0,~s=1,2$.
With the following calculations
\begin{gather*}
 \cbeta^\JS_1 = \left( \frac{\beta_1}{\ceps} + 1 \right)^2 \approx 10^{24} \beta_1^2,~
 \tbeta^\Z_1  = \frac{1 + \frac{\beta_2}{\teps}}{1 + \frac{\beta_2}{\beta_1+\teps}} \approx \frac{2}{7} 10^{40}\beta_2,~
 \tbeta^\ZR_1 = \frac{1 + \frac{\beta_2}{\teps^p}}{1 + \frac{\beta_2}{\left( \sqrt[p]{\beta_1}+\teps \right)^p}} \approx \frac{2}{7} 10^{40p}\beta_2, \\
 \cbeta^\JS_2 = \left( \frac{\beta_2}{\ceps} + 1 \right)^2 \approx 10^{24} \beta_2^2,~
 \tbeta^\Z_2  = \frac{1 + \frac{\beta_2}{\teps}}{1 + \frac{\beta_2}{\beta_2+\teps}} \approx \frac{10^{40}}{2} \beta_2,~
 \tbeta^\ZR_2 = \frac{1+\frac{\beta_2}{\teps^p}}{1+\frac{\beta_2}{\left( \sqrt[p]{\beta_2}+\teps \right)^p}} \approx \frac{10^{40p}}{2} \beta_2, \\
 \frac{\cbeta^\JS_1}{\cbeta^\JS_2} = \left( \frac{\beta_1+\ceps}{\beta_2+\ceps} \right)^2 \approx \frac{4}{25},~
 \frac{\tbeta^\Z_1}{\tbeta^\Z_2} = \frac{1 + \frac{\beta_2}{\beta_2+\teps}}{1 + \frac{\beta_2}{\beta_1+\teps}} \approx \frac{4}{7},
 \frac{\tbeta^\ZR_1}{\tbeta^\ZR_2} = \frac{1+\frac{\beta_2}{\left( \sqrt[p]{\beta_2}+\teps \right)^p}}{1+\frac{\beta_2}{\left( \sqrt[p]{\beta_1}+\teps \right)^p}} \approx \frac{4}{7},
\end{gather*}
we have $\omega^\JS_0 < \omega^\Z_0 < \omega^\ZR_0$ and $\omega^\ZR_s \ll \omega^\Z_s \ll \omega^\JS_s,~s=1,2$.
So we obtain
\begin{gather*}
 \omega^\N_0 \lessapprox 1,~~ \omega^\N_s \gtrapprox 0,~s=1,2,~~ \N \in \{ \JS, \M, \Z, \ZR \}, \\
 \omega^\M_s \ll \omega^\JS_s,~ \omega^\ZR_s \ll \omega^\Z_s \ll \omega^\JS_s,~ s = 1,2,
\end{gather*}
in correspondence to $x=0$ in Table \ref{tab:weight_comparison_first}.
\end{stencil}

\begin{stencil} \label{case:4}
$\{ I_{-2}, I_{-1}, I_{0}, I_{1}, I_{2} \}$ \\
Applying the techniques in Stencil \ref{case:3} to this stencil gives
\begin{gather*}
 \omega^\N_s \gtrapprox 0,~s=0,1,~~ \omega^\N_2 \lessapprox 1,~~ \N \in \{ \JS, \M, \Z, \ZR \}, \\
 \omega^\M_s \ll \omega^\JS_s,~ \omega^\ZR_s \ll \omega^\Z_s \ll \omega^\JS_s,~ s = 0,1,
\end{gather*}
which corresponds to $x=0.01$ in Table \ref{tab:weight_comparison_first}.
\end{stencil}

\begin{stencil} \label{case:5}
$\{ I_{-1}, I_{0}, I_{1}, I_{2}, I_{3} \}$ \\
The similar analysis as it is used in Stencil \ref{case:2} results in
\begin{gather*}
 \omega^\N_0 \gtrapprox 0,~\omega^\N_1 \lessapprox \frac{2}{3},~\omega^\N_2 \lessapprox \frac{1}{3},~\N \in \{ \JS, \Z, \ZR \}, \\
 \omega^\M_0 \gtrapprox 0,~ \omega^\M_1 \approx \frac{6164}{9241},~ \omega^\M_2 \approx \frac{3077}{9241}, \\
 \omega^\M_0 \ll \omega^\JS_0,~ \omega^\ZR_0 \ll \omega^\Z_0 \ll \omega^\JS_0, 
\end{gather*}
which is associated with $x=0.02$ in Table \ref{tab:weight_comparison_first}.
\end{stencil}

\begin{stencil} \label{case:6}
$\{ I_j, I_{j+1}, I_{j+2}, I_{j+3}, I_{j+4} \},~j \geqslant 0$ \\
We have the same conclusion at $x=0.03, 0.04$ in Table \ref{tab:weight_comparison_first}, as is done in Stencil \ref{case:1} for the smoothness indicators $\beta_0 = \beta_1 = \beta_2 = 0$.
\end{stencil}

All the nonlinear weights $\omega_s$ computed around the discontinuity for the Riemann problem \eqref{eq:weight_comparison} are presented in Table \ref{tab:weight_comparison_first}, with $u_L=1, u_R=0$ and $\nu = 0.5$.

Since $\alpha = 1$, by \eqref{eq:1d_scalar_hyperbolic_FV_flux} and \eqref{eq:LF_flux}, we have $\hat{f}_{i+1/2} = u^-_{i+1/2}$.
And based on \eqref{eq:3rd_reconstruction_minus} and \eqref{eq:weno_appr}, we obtain the following numerical fluxes:
\begin{enumerate}[label=\arabic*., font=\itshape]
\item $\hat{f}^{\Nc 1}_{j-5/2} = u_L,~j \leqslant 0$;
\item $\hat{f}^{\Nc 1}_{-3/2} = u_L + \frac{1}{6} \omega^\N_{2,-3/2} \delta$;
\item $\hat{f}^{\Nc 1}_{-1/2} = u_L - \frac{1}{3} \A^{\Nc 1}_{-1/2} \delta$;
\item $\hat{f}^{\Nc 1}_{1/2} = u_R - \frac{1}{6} \B^{\Nc 1}_{1/2} \delta$;
\item $\hat{f}^{\Nc 1}_{3/2} = u_R + \frac{1}{3} \omega^\N_{0,3/2} \delta$;
\item $\hat{f}^{\Nc 1}_{j+5/2} = u_R,~j \geqslant 0$,
\end{enumerate}
with $\N \in \{ \JS, \M, \Z, \ZR \}$ and $\A, \B$ given in Appendix.
However, due to round-off errors, the numerical fluxes computed with double precision do not fully agree with the analyzed ones above, as it is shown at $x = -0.01, 0$ in Table \ref{tab:numerical_flux_first}.

\begin{SCfigure}[0.5][htbp]
\includegraphics[width=0.5\textwidth]{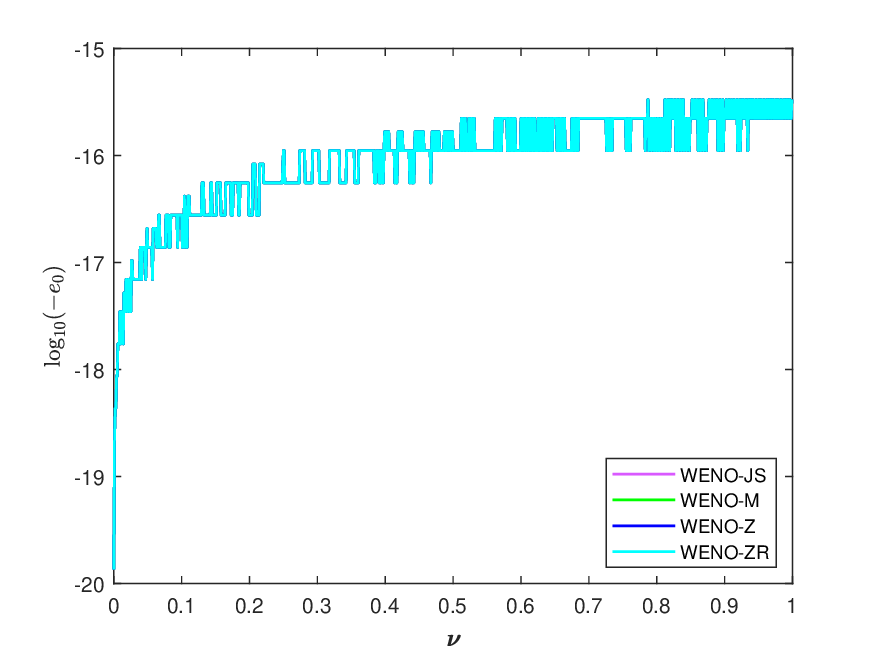}
\caption{The error $e^{\Nc 1}_0$ for cell $I_0$. Notice that the errors of WENO-JS, WENO-M, WENO-Z and WENO-ZR are almost the same.}
\label{fig:error1}
\end{SCfigure}

Let $\ubar^{\Nc k}_j$ denote the $k$th-stage approximation to $\ubar^1_j$ \eqref{eq:exact_sol} by WENO-$\N$.
We use \eqref{eq:RK1st} for the time integration at the first stage with the updating formula  
$$
   \ubar^{\Nc 1}_j = \ubar^0_j - \nu \left( \hat{f}^{\Nc 1}_{j+1/2} - \hat{f}^{\Nc 1}_{j-1/2} \right),
$$ 
then we have
\begin{equation} \label{eq:numerical_solution_first}
\begin{aligned}
 \ubar^{\Nc 1}_j &= \ubar^1_j,~j \leqslant -3, \\
 \ubar^{\Nc 1}_j &= \ubar^1_j + e^{\Nc 1}_j,~j = -2, \cdots, 2, \\
 \ubar^{\Nc 1}_j &= \ubar^1_j,~j \geqslant 3,
\end{aligned}
\end{equation}
with the errors $e^{\Nc 1}_j$ given in Appendix.
It is easy to see that all the errors in this stage are small and close to 0.
We could also check that $e^{\Nc 1}_{-2},~e^{\Nc 1}_{1} < 0$ and $e^{\Nc 1}_{-1},~e^{\Nc 1}_{2} > 0$, while whether $e^{\Nc 1}_0$ is positive or negative cannot be determined by the formula. 
Table \ref{tab:numerical_solution_first} and Figure \ref{fig:error1} show that $e^{\Nc 1}_0 < 0$.
With many experiments tested, we notice that $e^{\Nc 1}_{-2}= 0$ for $0 < \nu \leqslant 1$, and either $e^{\Nc 1}_{-1}=0$ if $0 < \nu \leqslant 0.5$ or $| e^{\Nc 1}_{-1} |$ is on the order of $\dbp =$ 2.220e-16 if $0.5 < \nu \leqslant 1$.
So the formulas \eqref{eq:numerical_solution_first} agree in part with the results of the approximations after the first stage of one TVD Runge-Kutta time step displayed in Table \ref{tab:numerical_solution_first} because of the double precision.
Therefore, $\ubar^{\Nc 1}_j = 1,~j = -2, -1$ corresponding to $x = -0.015, -0.005$ with $e^{\Nc 1}_j = 0$, and $\ubar^{\Nc 1}_{0} = 0.5-\dbp/2$ at $x = 0.005$, for $\N \in \{ \JS, \M, \Z, \ZR \}$.

\begin{sidewaystable}
\renewcommand{\arraystretch}{1.2}
\centering

\vspace{4.6in}
\begin{tabular}{ccccccccc} 
\hline
$x$ & -0.03 & -0.02 & -0.01 & 0 & 0.01 & 0.02 & 0.03 & 0.04 \\ 
\hline 
$\omega^\JS_0$ & 0.1        & 0.1        & 0.142857 & 1 & 3.000e-26  & 6.250e-26  & 0.1        & 0.1 \\  
$\omega^\M_0$  & 0.1+$\vez$ & 0.1+$\vez$ & 0.127255 & 1 & 3.300e-81  & 7.626e-81  & 0.1+$\vez$ & 0.1+$\vez$ \\  
$\omega^\Z_0$  & 0.1        & 0.1        & 0.142857 & 1 & 2.000e-41  & 1.667e-41  & 0.1        & 0.1 \\
$\omega^\ZR_0$ & 0.1        & 0.1        & 0.142857 & 1 & 2.000e-121 & 1.667e-121 & 0.1        & 0.1 \\
\hline 
$\omega^\JS_1$ & 0.6 & 0.6 & 0.857143 & 3.375e-24  & 1.125e-24  & 0.666667 & 0.6 & 0.6 \\
$\omega^\M_1$  & 0.6 & 0.6 & 0.872745 & 9.000e-80  & 3.000e-80  & 0.667027 & 0.6 & 0.6 \\  
$\omega^\Z_1$  & 0.6 & 0.6 & 0.857143 & 6.300e-40  & 2.100e-40  & 0.666667 & 0.6 & 0.6 \\
$\omega^\ZR_1$ & 0.6 & 0.6 & 0.857143 & 6.300e-120 & 2.100e-120 & 0.666667 & 0.6 & 0.6 \\
\hline 
$\omega^\JS_2$ & 0.3        & 0.3        & 2.411e-25  & 2.700e-13  & 1 & 0.333333 & 0.3        & 0.3 \\  
$\omega^\M_2$  & 0.3+$\vet$ & 0.3+$\vet$ & 1.321e-80  & 1.170e-80  & 1 & 0.332973 & 0.3+$\vet$ & 0.3+$\vet$ \\  
$\omega^\Z_2$  & 0.3        & 0.3        & 6.429e-41  & 1.800e-40  & 1 & 0.333333 & 0.3        & 0.3 \\
$\omega^\ZR_2$ & 0.3        & 0.3        & 6.429e-121 & 1.800e-120 & 1 & 0.333333 & 0.3        & 0.3 \\    
\hline
\end{tabular}
\caption{The nonlinear weights $\omega^\N_s$ near the discontinuity at $x=0$, with the linear weights $(d_0, d_1, d_2) = (0.1, 0.6, 0.3)$, for WENO-JS, WENO-M, WENO-Z and WENO-ZR at the first stage \eqref{eq:RK1st} of the single TVD Runge-Kutta time step for the Riemann problem \eqref{eq:weight_comparison} with $\vez = \cm 2.776 \e \cm 17$ and $\vet = 5.551\e \cm 17$.}
\label{tab:weight_comparison_first}

\vspace{0.2in}

\begin{tabular}{ccccccccc} 
\hline
$x$ & -0.03 & -0.02 & -0.01 & 0 & 0.01 & 0.02 & 0.03 & 0.04 \\
\hline
$\hat{f}^\JS_j$ & 1 & 1 & 1 & 1-2.220e-16 & -2.125e-25  & 2.083e-26  & 0 & 0 \\  
$\hat{f}^\M_j$  & 1 & 1 & 1 & 1-2.220e-16 & -7.750e-81  & 2.542e-81  & 0 & 0 \\  
$\hat{f}^\Z_j$  & 1 & 1 & 1 & 1-2.220e-16 & -5.167e-41  & 5.556e-42  & 0 & 0 \\
$\hat{f}^\ZR_j$ & 1 & 1 & 1 & 1-2.220e-16 & -5.167e-121 & 5.556e-122 & 0 & 0 \\ 
\hline
\end{tabular}
\caption{The numerical fluxes $\hat{f}^\N_j$ near the discontinuity at $x=0$ for WENO-JS, WENO-M, WENO-Z and WENO-ZR at the first stage \eqref{eq:RK1st} of the single TVD Runge-Kutta time step for the Riemann problem \eqref{eq:weight_comparison}.}
\label{tab:numerical_flux_first}

\vspace{0.2in}

\begin{tabular}{cccccccc} 
\hline 
$x$ & -0.025 & -0.015 & -0.005 & 0.005 & 0.015 & 0.025 & 0.035  \\ 
\hline
$\ubar^\JS_j$ & 1 & 1 & 1 & 0.5-1.110e-16 & -1.167e-25  & 1.042e-26  & 0 \\  
$\ubar^\M_j$  & 1 & 1 & 1 & 0.5-1.110e-16 & -5.146e-81  & 1.271e-81  & 0 \\  
$\ubar^\Z_j$  & 1 & 1 & 1 & 0.5-1.110e-16 & -2.861e-41  & 2.778e-42  & 0 \\
$\ubar^\ZR_j$ & 1 & 1 & 1 & 0.5-1.110e-16 & -2.861e-121 & 2.778e-122 & 0 \\ 
$\ubar_j$     & 1 & 1 & 1 & 0.5           & 0           & 0          & 0 \\ 
\hline
\end{tabular}
\caption{The numerical solutions $\ubar^\N_j$, compared with the exact solution $\ubar_j$, near the discontinuity at $x=0$ for WENO-JS, WENO-M, WENO-Z and WENO-ZR at the first stage \eqref{eq:RK1st} of the single TVD Runge-Kutta time step for the Riemann problem \eqref{eq:weight_comparison}.}
\label{tab:numerical_solution_first}

\end{sidewaystable}

\subsubsection{Second stage of the TVD Runge-Kutta method} \label{sec:RK2nd}
Next we move to the second stage,
\begin{equation} \label{eq:RK2nd} 
 u^{(2)} = \frac{3}{4} u^n + \frac{1}{4} u^{(1)} + \frac{1}{4} \Delta t L \left( u^{(1)} \right),
\end{equation}
and compute the numerical solutions by WENO schemes.
Based on the observations in Section \ref{sec:Euler}, 9 stencils should be considered:
\begin{enumerate}[label=\arabic*., font=\itshape]
\setcounter{enumi}{6}
\item $\{ I_{j-4}, I_{j-3}, I_{j-2}, I_{j-1}, I_j \},~j \leqslant -1$;
\item $\{ I_{-4}, I_{-3}, I_{-2}, I_{-1}, I_{0} \}$; \label{item:stt2}
\item $\{ I_{-3}, I_{-2}, I_{-1}, I_{0}, I_{1} \}$;
\item $\{ I_{-2}, I_{-1}, I_{0}, I_{1}, I_{2} \}$;
\item $\{ I_{-1}, I_{0}, I_{1}, I_{2}, I_{3} \}$;    
\item $\{ I_{0}, I_{1}, I_{2}, I_{3}, I_{4} \}$;     \label{item:end2}
\item $\{ I_{1}, I_{2}, I_{3}, I_{4}, I_{5} \}$;
\item $\{ I_{2}, I_{3}, I_{4}, I_{5}, I_{6} \}$;
\item $\{ I_j, I_{j+1}, I_{j+2}, I_{j+3}, I_{j+4} \},~j \geqslant 3$,
\end{enumerate}
among which Stencils \ref{case:8} - \ref{case:12} contain the jump discontinuity. 
We also assume from those observations that
\begin{itemize} 
\item $| e^{\Nc 1}_0 |$ is on the order of $\dbp$ for $\N \in \{ \JS, \M, \Z, \ZR \}$, \\
      $| e^{\JS,\,1}_j | \ll \ceps,~| e^{\M,\,1}_j |, | e^{\ZR,\,1}_j | \ll \teps$ and $| e^{\Z,\,1}_j |$ is on the order of $\teps$ for $j = 1,2$;
\item $\delta \gg e^{\Nc 1}_j,~\N \in \{ \JS, \M, \Z, \ZR \},~j = 0,1,2$.
\end{itemize}
In addition, $e^{\Nc 1}_0,~e^{\Nc 1}_1 < 0$ and $e^{\Nc 1}_2 > 0$ hold from Section \ref{sec:Euler}.

\begin{stencil} \label{case:7}
$\{ I_{j-4}, I_{j-3}, I_{j-2}, I_{j-1}, I_j \},~j \leqslant -1$ \\
The same conclusion as in Stencil \ref{case:1}, can be made corresponding to $x=-0.03, -0.02$ in Table \ref{tab:weight_comparison_second}.
\end{stencil}

\begin{stencil} \label{case:8}
$\{ I_{-4}, I_{-3}, I_{-2}, I_{-1}, I_{0} \}$ \\
With $\beta^\N_0 = \beta^\N_1 = 0,~\beta^\N_2 = \frac{4}{3} \left[ ( 1-\nu ) \delta - e^\N_0 \right]^2$ for $\N \in \{ \JS, \M, \Z, \ZR \}$, we obtain the same results as in Stencil \ref{case:2}, which corresponds to $x=-0.01$ in Table \ref{tab:weight_comparison_second}.
\end{stencil}

\begin{stencil} \label{case:9}
$\{ I_{-3}, I_{-2}, I_{-1}, I_{0}, I_{1} \}$ \\
The smoothness indicators in this stencil are
\begin{align*}
 \beta^\N_0 &= 0,~\beta^\N_1 = \frac{4}{3} \left[ ( 1-\nu ) \delta - e^\N_0 \right]^2, \\
 \beta^\N_2 &= \frac{13}{12} \left[ ( 1-2\nu ) \delta - 2e^\N_0 + e^\N_1 \right]^2 + 
               \frac{1}{4}   \left[ ( 3-4\nu ) \delta - 4e^\N_0 + e^\N_1 \right]^2
\end{align*} 
for $\N \in \{ \JS, \M, \Z, \ZR \}$.
Similar to the analysis in Stencil \ref{case:3}, we have
\begin{gather*}
 \omega^\N_0 \lessapprox 1,~~ \omega^\N_s \gtrapprox 0,~s=1,2,~~ \N \in \{ \JS, \M, \Z, \ZR \}, \\
 \omega^\Z_s \gg \omega^\ZR_s,~s= 1,2,
\end{gather*}
which is associated with $x=0$ in Table \ref{tab:weight_comparison_second}.
\end{stencil}

\begin{stencil} \label{case:10}
$\{ I_{-2}, I_{-1}, I_{0}, I_{1}, I_{2} \}$ \\
For $\N \in \{ \JS, \M, \Z, \ZR\}$, the smoothness indicators satisfy 
\begin{align*} 
 \beta^\N_0 &= \frac{10}{3} \left[ ( 1-\nu ) \delta - e^\N_0 \right]^2,~\beta^\N_1 = \frac{13}{12} \left[ ( 1-2\nu ) \delta - 2e^\N_0 + e^\N_1 \right]^2 + \frac{1}{4} \left( \delta - e^\N_1 \right)^2,  \\
 \beta^\N_2 &= \frac{13}{12} \left(  \nu \delta +  e^\N_0 - 2e^\N_1 + e^\N_2 \right)^2 + 
               \frac{1}{4}   \left( 3\nu \delta + 3e^\N_0 - 4e^\N_1 + e^\N_2 \right)^2.
\end{align*}
Then we can estimate 
$$
   \beta^\N_0 \approx \frac{10}{3}(1-\nu)^2 \delta^2,~
   \beta^\N_1 \approx \frac{13}{12}(1- 2 \nu)^2 \delta^2 + \frac{1}{4} \delta^2,~
   \beta^\N_2 \approx \frac{10}{3} \nu^2 \delta^2.
$$
Figure \ref{fig:weight12} plots the graphs of $\omega^\N_s,~s = 0,1,2$ versus $\nu$.
Unlike the nonlinear weights in previous stencils, the weights $\omega^\N_s$ in this stencil vary with $\nu$ in the unit interval.
These weights are important components of the errors $e^{\Nc 2}_0$ \eqref{eq:e2b0} and $e^{\Nc 2}_1$ \eqref{eq:e2b1}.
For $\nu$ fixed in the range $0 < \nu \leqslant 0.5$, the value of $\omega_0$ from WENO-ZR is the largest, and WENO-JS gives rise to the smallest value of $\omega_0$ with the values of $\omega_0$ by WENO-M and WENO-Z in between, which we will use latter in this subsection. 
\begin{figure}[htbp]
\centering
\includegraphics[width=0.325\textwidth]{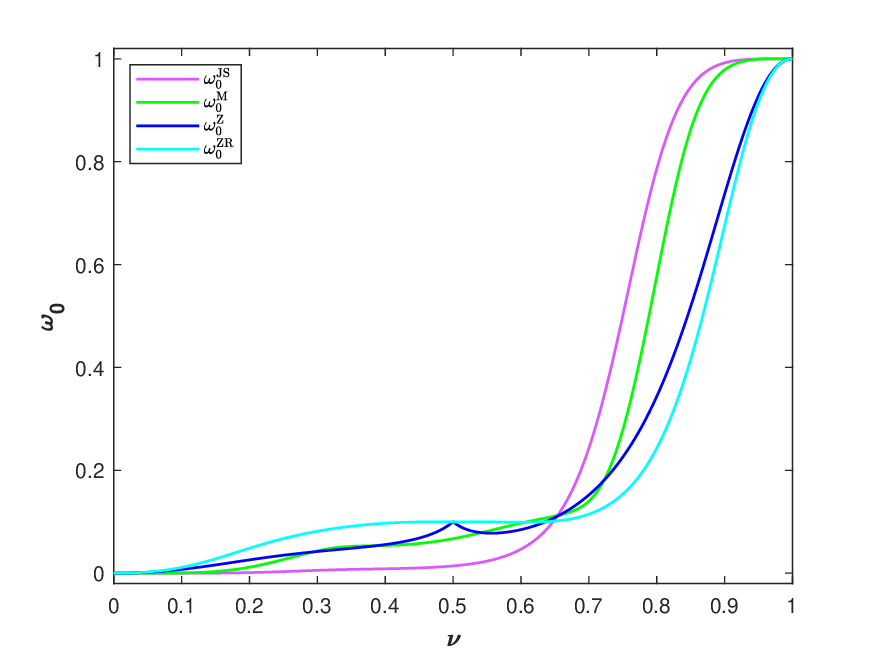}
\includegraphics[width=0.325\textwidth]{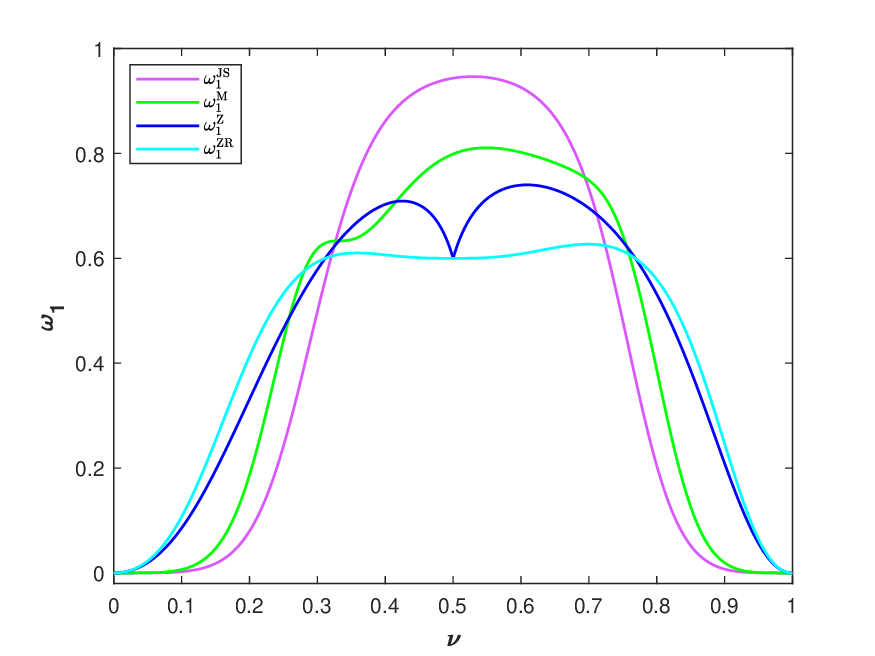}
\includegraphics[width=0.325\textwidth]{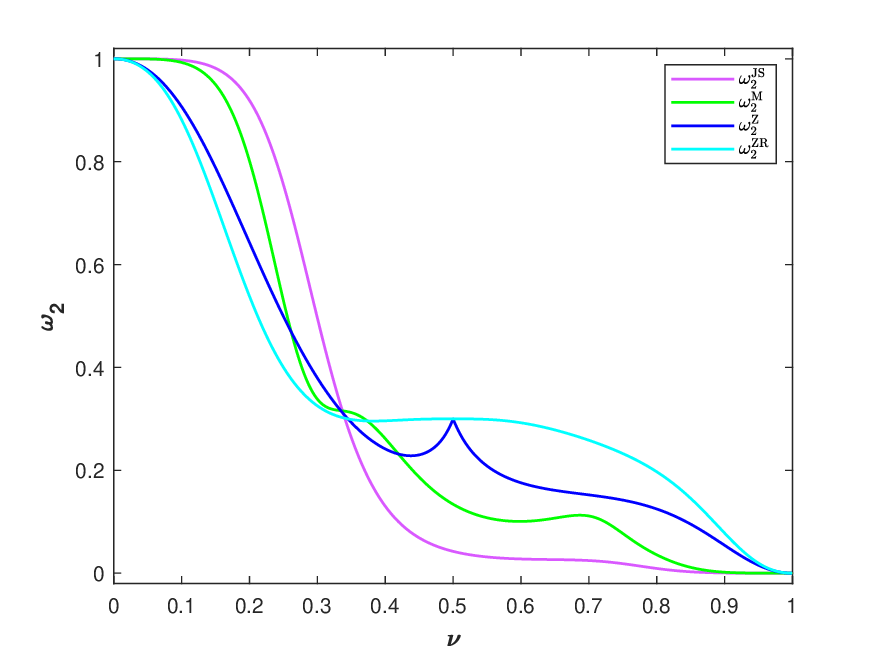}
\vspace{-0.3cm}
\caption{The nonlinear weights $\omega^\N_0$ (left), $\omega^\N_1$ (middle) and $\omega^\N_2$ (right) in Stencil \ref{case:10} by different WENO techniques.}
\label{fig:weight12}
\end{figure}
\end{stencil}

\begin{stencil} \label{case:11}
$\{ I_{-1}, I_{0}, I_{1}, I_{2}, I_{3} \}$ \\
Using the smoothness indicators 
\begin{align*}
 \beta^\N_0 &= \frac{13}{12} \left[ ( 1-2\nu ) \delta - 2e^\N_0 + e^\N_1 \right]^2 + 
               \frac{1}{4}   \left( ( 1-4\nu ) \delta - 4e^\N_0 + 3e^\N_1 \right]^2, \\
 \beta^\N_1 &= \frac{13}{12} \left[ \nu \delta + e^\N_0 - 2e^\N_1 + e^\N_2 \right]^2 + 
               \frac{1}{4}   \left[ \nu \delta + e^\N_0 - e^\N_2 \right]^2, \\
 \beta^\N_2 &= \frac{13}{12} \left( e^\N_1-2e^\N_2 \right)^2 + \frac{1}{4} \left( 3e^\N_1-4e^\N_2 \right)^2,
\end{align*}
for $\N \in \{ \JS, \M, \Z, \ZR \}$, we obtain $\beta^\N_0, \beta^\N_1 \gg \beta^\N_2$.
A similar analysis as in Stencil \ref{case:3} gives 
\begin{gather*}
 \omega^\N_s \gtrapprox 0,~s=0,1,~~ \omega^\N_2 \lessapprox 1,~~ \N \in \{ \JS, \M, \Z, \ZR \}, \\
 \omega^\Z_s \gg \omega^\ZR_s,~s = 0,1,
\end{gather*}
in correspondence to $x=0.02$ in Table \ref{tab:weight_comparison_second}.
\end{stencil}

\begin{stencil} \label{case:12}
$\{ I_{0}, I_{1}, I_{2}, I_{3}, I_{4} \}$ \\
Note that for $\N \in \{ \JS, \M, \Z, \ZR \}$,
\begin{align*}
 \beta^\N_0 &= \frac{13}{12} \left( \nu \delta + e^\N_0 - 2e^\N_1 + e^\N_2 \right)^2 + 
               \frac{1}{4}   \left( \nu \delta + e^\N_0 - 4e^\N_1 + 3e^\N_2 \right)^2, \\
 \beta^\N_1 &= \frac{13}{12} \left( e^\N_1-2e^\N_2 \right)^2 + \frac{1}{4} \left( e^\N_1 \right)^2,~
 \beta^\N_2  = \frac{10}{3} \left( e^\N_2 \right)^2,
\end{align*} 
where $\beta^\N_0 \gg \beta^\N_1, \beta^\N_2$.
Applying the techniques in Stencil \ref{case:2} would lead to
\begin{gather*}
 \omega^\N_0 \gtrapprox 0,~\omega^\N_1 \approx \frac{2}{3},~\omega^\N_2 \approx \frac{1}{3},~\N \in \{ \JS, \Z, \ZR \}, \\
 \omega^\M_0 \gtrapprox 0,~\omega^\M_1 \approx \frac{6164}{9241},~\omega^\M_2 \approx \frac{3077}{9241}, \\
 \omega^\Z_0 \gg \omega^\ZR_0,
\end{gather*}
which corresponds to $x=0.03$ in Table \ref{tab:weight_comparison_second}.
\end{stencil}

\begin{stencil} \label{case:13}
$\{ I_{1}, I_{2}, I_{3}, I_{4}, I_{5} \}$ \\
The smoothness indicators satisfy 
\begin{align*}
 \beta^\N_0 &= \frac{13}{12} \left( e^\N_1-2e^\N_2 \right)^2 + \frac{1}{4} \left( e^\N_1-4e^\N_2 \right)^2, \\
 \beta^\N_1 &= \frac{4}{3} \left( e^\N_2 \right)^2,~\beta^\N_2  = 0,
\end{align*} 
for $\N \in \{ \JS, \M, \Z, \ZR \}$.
It is easy to see from \eqref{eq:nw_JS1} that $\cbeta^\JS_1, \cbeta^\JS_2, \cbeta^\JS_2 / \cbeta^\JS_1 \lessapprox 1$, so that $d_1/\cbeta^\JS_1 \gtrapprox d_1$, $d_2/\cbeta^\JS_2 \gtrapprox d_2$ and $d_0 + d_1/\cbeta^\JS_1 + d_2/\cbeta^\JS_2 \gtrapprox 1$, and consequently 
$$ 
   \omega^\JS_0 \lessapprox d_0,~ \omega^\JS_1 \approx d_1,~ \omega^\JS_2 \gtrapprox d_2.
$$
A similar analysis shows that
\begin{gather*}
 \omega^\N_0 \lessapprox d_0,~\omega^\N_1 \approx d_1,~\omega^\N_2 \gtrapprox d_2,~\N \in \{ \Z, \ZR \}, \\
 \tomega^\JS_0 \lessapprox d_0,~\tomega^\JS_1 \approx d_1,~\tomega^\JS_2 \gtrapprox d_2, 
\end{gather*}
and it is expected that $\omega^\M_s \approx d_s,~s = 0,1,2$.
Therefore, we have
\begin{gather*}
 \omega^\N_s \approx d_s,~ s = 0,1,2,~ \N \in \{ \JS, \M, \Z, \ZR \},
\end{gather*}
at $x=0.04$ in Table \ref{tab:weight_comparison_second}.
\end{stencil}

\begin{stencil} \label{case:14}
$\{ I_{2}, I_{3}, I_{4}, I_{5}, I_{6} \}$ \\
With the smoothness indicators $\beta^\N_0 = \frac{4}{3} \left( e^\N_2 \right)^2$ and $\beta^\N_1 = \beta^\N_2 = 0$, 
\eqref{eq:nw_JS1} gives $\cbeta^\JS_1, \cbeta^\JS_2 \lessapprox 1$, $\cbeta^\JS_1 / \cbeta^\JS_2 = 1$, resulting in $\omega^\JS_0 \lessapprox d_0$ and $\omega^\JS_s \gtrapprox d_s,~s=1,2$.
Applying similar analysis yields
\begin{gather*}
 \omega^\N_0 \lessapprox d_0,~~ \omega^\N_s \gtrapprox d_s,~s=1,2,~~ \N \in \{ \Z, \ZR \}, \\
 \tomega^\JS_0 \lessapprox d_0,~~ \tomega^\JS_s \gtrapprox d_s,~s=1,2.
\end{gather*}
Hence $\omega^\M_s \approx d_s,~s = 0,1,2$.
In summary, we have
\begin{gather*}
 \omega^\N_s \approx d_s,~ s = 0,1,2,~ \N \in \{ \JS, \M, \Z, \ZR \},
\end{gather*}
corresponding to $x=0.05$ in Table \ref{tab:weight_comparison_second}.
\end{stencil}

\begin{stencil} \label{case:15}
$\{ I_j, I_{j+1}, I_{j+2}, I_{j+3}, I_{j+4} \},~j \geqslant 3$ \\
Each nonlinear weight is equal to its corresponding linear one at $x=0.06, 0.07$ in Table \ref{tab:weight_comparison_second} as in Stencil \ref{case:1}.
\end{stencil}

Based on \eqref{eq:3rd_reconstruction_minus}, \eqref{eq:weno_appr} and the assumptions given in the beginning of this subsection, we have the following numerical fluxes,
\begin{enumerate}[label=\arabic*., font=\itshape]
\setcounter{enumi}{6}
\item $\hat{f}^{\Nc 2}_{j-1/2} = u_L,~~j \leqslant -2$;
\item $\hat{f}^{\Nc 2}_{-3/2} = u_L + \frac{1}{6} \omega^{\Nc 2}_{2,-3/2} (1-\nu) \delta - \frac{1}{6} \omega^{\Nc 2}_{2,-3/2} e^{\Nc 1}_0$;
\item $\hat{f}^{\Nc 2}_{-1/2} = u_L + \frac{1}{6} \omega^{\Nc 2}_{2,-1/2} \delta - \frac{1}{6} \C^{\Nc 2}_{-1/2} (1-\nu) \delta + \frac{1}{6} \left( \C^{\Nc 2}_{-1/2} e^{\Nc 1}_0 - \omega^{\Nc 2}_{2,-1/2} e^{\Nc 1}_1 \right)$; 
\item $\hat{f}^{\Nc 2}_{1/2} = u_L - \frac{1}{3} \A^{\Nc 2}_{1/2} \delta - \frac{1}{6} \D^{\Nc 2}_{1/2} (1-\nu) \delta + \frac{1}{6} \left( \D^{\Nc 2}_{1/2} e^{\Nc 1}_0 + \C^{\Nc 2}_{1/2} e^{\Nc 1}_1 - \omega^{\Nc 2}_{2,1/2} e^{\Nc 1}_2 \right)$;
\item $\hat{f}^{\Nc 2}_{3/2} = u_R - \frac{1}{6} \B^{\Nc 2}_{3/2} \delta + \frac{1}{6} \E^{\Nc 2}_{3/2} (1-\nu) \delta + \frac{1}{6} \left( - \E^{\Nc 2}_{3/2} e^{\Nc 1}_0 + \D^{\Nc 2}_{3/2} e^{\Nc 1}_1 + \C^{\Nc 2}_{3/2} e^{\Nc 1}_2 \right)$;
\item $\hat{f}^{\Nc 2}_{5/2} = u_R + \frac{1}{3} \omega^{\Nc 2}_{0,5/2} \nu \delta + \frac{1}{6} \left( 2 \omega^{\Nc 2}_{0,5/2} e^{\Nc 1}_0 - \E^{\Nc 2}_{5/2} e^{\Nc 1}_1 + \D^{\Nc 2}_{5/2} e^{\Nc 1}_2 \right)$;
\item $\hat{f}^{\Nc 2}_{7/2} = u_R + \frac{1}{6} \left( 2 \omega^{\Nc 2}_{0,7/2} e^{\Nc 1}_1 - \E^{\Nc 2}_{7/2} e^{\Nc 1}_2 \right)$;
\item $\hat{f}^{\Nc 2}_{9/2} = u_R + \frac{1}{3} \omega^{\Nc 2}_{0,9/2} e^{\Nc 1}_2$;
\item $\hat{f}^{\Nc 2}_{j+1/2} = u_R,~j \geqslant 5$;
\end{enumerate} 
for $\N \in \{ \JS, \M, \Z, \ZR \}$ and $\A, \B, \C, \D, \E$ given in Appendix.
Because of round-off errors, the numerical fluxes with double precision in Table \ref{tab:numerical_flux_second} are not consistent with the analyzed ones above at $x = -0.01, 0$, as is the case in Section \ref{sec:Euler}.

Using the updating formula  
$$
   \ubar^{\Nc 2}_j = \frac{3}{4} \ubar^0_j + \frac{1}{4} \ubar^{\Nc 1}_j - \frac{1}{4} \nu \left( \hat{f}^{\Nc 2}_{j+1/2} - \hat{f}^{\Nc 2}_{j-1/2} \right),
$$
from \eqref{eq:RK2nd} for the time integration at the second stage, the jump discontinuity remains unchanged at $x = \Delta t$, and we have the updated approximations
\begin{equation} \label{eq:numerical_solution_second}
\begin{aligned}
 \ubar^{\Nc 2}_j &= \ubar^1_j,~j \leqslant -3, \\
 \ubar^{\Nc 2}_j &= \ubar^1_j + e^{\Nc 2}_j,~j = -2, \cdots, 5, \\
 \ubar^{\Nc 2}_j &= \ubar^1_j,~j \geqslant 6,
\end{aligned}
\end{equation}
with the errors $e^{\Nc 2}_j$ given in Appendix.
The formulas \eqref{eq:numerical_solution_second} agree with the approximations after the second stage in Table \ref{tab:numerical_solution_second} except at $x = -0.015, -0.005$, corresponding to $\ubar^{\Nc 2}_j = 1,~j = -2, -1$ with $e^{\Nc 2}_j = 0$ for $\N \in \{ \JS, \M, \Z, \ZR \}$, due to the round-off errors.

It is worth noting from Table \ref{tab:numerical_solution_second} that $e^{\Nc 2}_0$ and $e^{\Nc 2}_1$, corresponding to $x = 0.005, 0.015$, are large compared to the rest of errors. 
From Appendix, we find that $-\frac{3}{4} \nu \delta + \frac{1}{24} \left[ 2 \A^{\Nc 2}_{1/2} + \D^{\Nc 2}_{1/2} (1-\nu) \right]\nu \delta$ dominates $e^{\Nc 2}_0$ whereas $\frac{1}{4} \nu \delta - \frac{1}{24} \left[ 2 \A^{\Nc 2}_{1/2} + \D^{\Nc 2}_{1/2} (1-\nu) \right] \nu \delta$ is the main contribution to $e^{\Nc 2}_1$.  
Then
\begin{align}
 e^{\Nc 2}_0 &\approx -\frac{3}{4} \nu \delta + \frac{\nu \delta}{24} \left[ 2 \A^{\Nc 2}_{1/2} + \D^{\Nc 2}_{1/2} (1-\nu) \right] \label{eq:e2b0} \\
             &= - \frac{\nu \delta}{24} \left[ (7+11\nu) + (4-6\nu) \omega^{\Nc 2}_{1,1/2} + (5-9\nu) \omega^{\Nc 2}_{2,1/2} \right], \nonumber \\
 e^{\Nc 2}_1 &\approx \frac{1}{4} \nu \delta - \frac{\nu \delta}{24} \left[ 2 \A^{\Nc 2}_{1/2} + \D^{\Nc 2}_{1/2} (1-\nu) \right] \label{eq:e2b1} \\ 
             &= \frac{\nu \delta}{24} \left[ 2\nu - (5-9\nu) \omega^{\Nc 2}_{0,1/2} - (1-3\nu) \omega^{\Nc 2}_{1,1/2} \right], \nonumber \\
 e^{\Nc 2}_0 - e^{\Nc 2}_1 &\approx -\nu \delta + \frac{\nu \delta}{12} \left[ 2 \A^{\Nc 2}_{1/2} + \D^{\Nc 2}_{1/2} (1-\nu) \right]. \label{eq:e2b01}
\end{align}

\begin{figure}[htbp]
\centering
\includegraphics[width=0.325\textwidth]{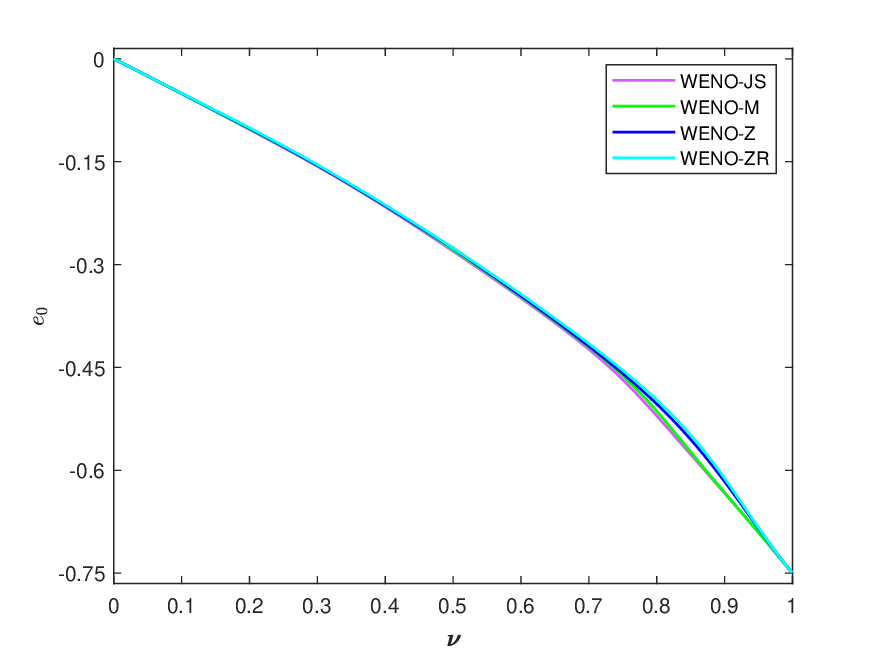}
\includegraphics[width=0.325\textwidth]{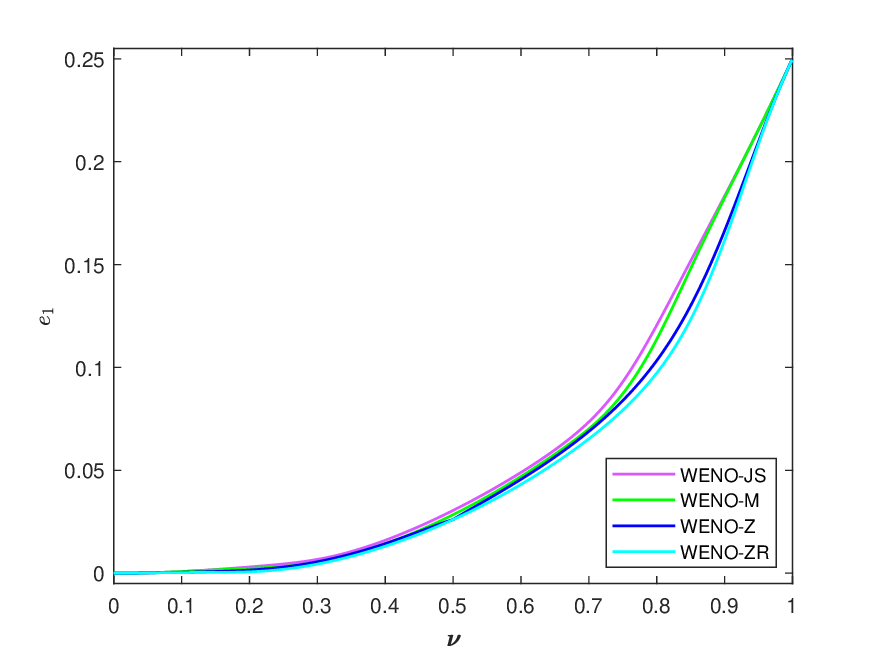}
\includegraphics[width=0.325\textwidth]{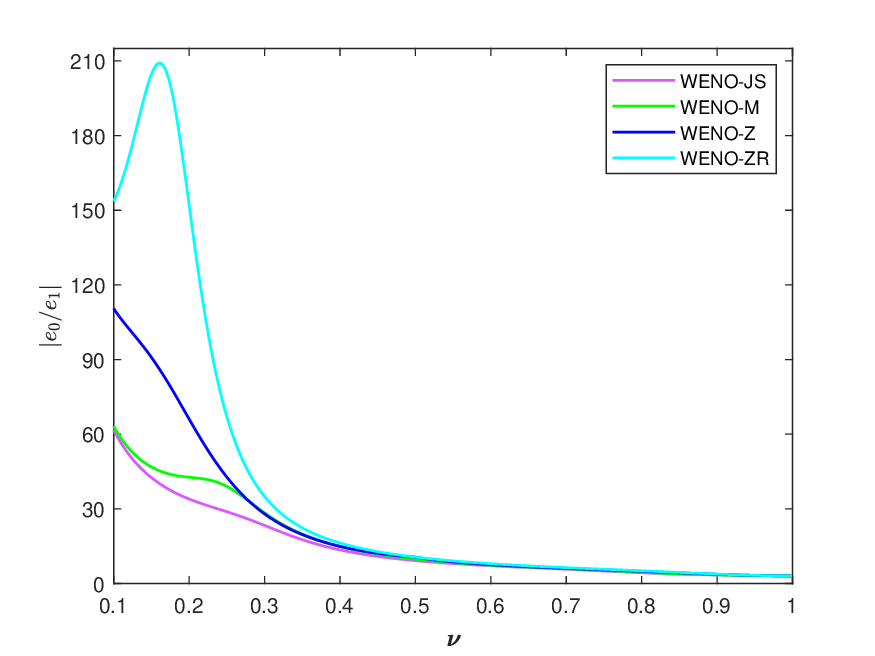}
\vspace{-0.3cm}
\caption{$e^{\Nc 2}_0$ (left), $e^{\Nc 2}_1$ (middle) and $\left| e^{\Nc 2}_0/e^{\Nc 2}_1 \right|$ (right).}
\label{fig:error2}
\end{figure}

It is obvious from \eqref{eq:e2b0} that $e^{\Nc 2}_0 < 0$ for $0 < \nu \leqslant 0.5$, which could be verified in Figure \ref{fig:error2}.
Then our goal is to reduce the values of both $\omega^{\Nc 2}_{1,1/2}$ and $\omega^{\Nc 2}_{2,1/2}$ corresponding to Stencil \ref{case:10}.
Roughly speaking, the larger value of $\omega^{\Nc 2}_{0,1/2}$ in Figure \ref{fig:weight12} is what we want.

From Figure \ref{fig:error2}, we have $e^{\Nc 2}_1 > 0$ .
When $0 < \delta < \frac{1}{3}$, the larger $\omega^{\Nc 2}_{0,1/2}$ and $\omega^{\Nc 2}_{1,1/2}$ are desired to suppress $e^{\Nc 2}_1$.
For $\frac{1}{3} < \delta < 0.5$, we want the larger $\omega^{\Nc 2}_{0,1/2}$ and the smaller $\omega^{\Nc 2}_{1,1/2}$ to reduce $e^{\Nc 2}_1$.
This gives one way to reduce the error $e^{\Nc 2}_1$, at least increasing the value of $\omega^{\Nc 2}_{0,1/2}$.

So at the second stage, we wish to increase the value of $\omega^{\Nc 2}_{0,1/2}$ in order to reduce the errors $e^{\Nc 2}_0$ and $e^{\Nc 2}_1$.
Figure \ref{fig:weight12} shows that WENO-ZR produces larger value of $\omega^{\ZR,\,2}_{0,1/2}$ for $0 < \nu \leqslant 0.5$, which results in better numerical solutions $\ubar^{\ZR,\,2}_0$ and $\ubar^{\ZR,\,2}_1$ with smaller errors $e^{\ZR,\,2}_0$ and $e^{\ZR,\,2}_1$ in Figure \ref{fig:error2}.

\begin{sidewaystable}
\renewcommand{\arraystretch}{1.2}
\centering

\vspace{4.6in}
\begin{tabular}{cccccccccccc} 
\hline
$x$ & -0.03 & -0.02 & -0.01 & 0 & 0.01 & 0.02 & 0.03 & 0.04 & 0.05 & 0.06 & 0.07 \\ 
\hline
$\omega^\JS_0$ & 0.1        & 0.1        & 0.142857 & 1 & 0.014151      & 5.333e-24  & 1.000e-24  & 0.1        & 0.1        & 0.1        & 0.1 \\
$\omega^\M_0$  & 0.1+$\vez$ & 0.1+$\vez$ & 0.127255 & 1 & 0.067167      & 5.867e-79  & 1.220e-79  & 0.1+$\vez$ & 0.1+$\vez$ & 0.1+$\vez$ & 0.1+$\vez$ \\    
$\omega^\Z_0$  & 0.1        & 0.1        & 0.142857 & 1 & 0.1-1.388e-16 & 2.667e-40  & 6.667e-41  & 0.1        & 0.1        & 0.1        & 0.1 \\
$\omega^\ZR_0$ & 0.1        & 0.1        & 0.142857 & 1 & 0.1           & 2.667e-120 & 6.667e-121 & 0.1        & 0.1        & 0.1        & 0.1 \\
\hline
$\omega^\JS_1$ & 0.6 & 0.6 & 0.857143 & 5.400e-23  & 0.943396      & 1.800e-23  & 0.666667 & 0.6 & 0.6 & 0.6 & 0.6 \\  
$\omega^\M_1$  & 0.6 & 0.6 & 0.872745 & 1.440e-78  & 0.798816      & 4.800e-79  & 0.667027 & 0.6 & 0.6 & 0.6 & 0.6 \\   
$\omega^\Z_1$  & 0.6 & 0.6 & 0.857143 & 4.200e-39  & 0.6+6.661e-16 & 1.400e-39  & 0.666667 & 0.6 & 0.6 & 0.6 & 0.6 \\
$\omega^\ZR_1$ & 0.6 & 0.6 & 0.857143 & 4.200e-119 & 0.6           & 1.400e-119 & 0.666667 & 0.6 & 0.6 & 0.6 & 0.6 \\
\hline	
$\omega^\JS_2$ & 0.3        & 0.3        & 3.857e-24  & 4.800e-23  & 0.042453      & 1 & 0.333333 & 0.3        & 0.3        & 0.3        & 0.3 \\   
$\omega^\M_2$  & 0.3+$\vet$ & 0.3+$\vet$ & 2.114e-79  & 2.080e-78  & 0.134016      & 1 & 0.332973 & 0.3+$\vet$ & 0.3+$\vet$ & 0.3+$\vet$ & 0.3+$\vet$ \\    
$\omega^\Z_2$  & 0.3        & 0.3        & 2.571e-40  & 2.400e-39  & 0.3-3.886e-16 & 1 & 0.333333 & 0.3        & 0.3        & 0.3        & 0.3 \\
$\omega^\ZR_2$ & 0.3        & 0.3        & 2.571e-120 & 2.400e-119 & 0.3           & 1 & 0.333333 & 0.3        & 0.3        & 0.3        & 0.3 \\   
\hline
\end{tabular}
\caption{The nonlinear weights $\omega^\N_s$ near the discontinuity at $x = 0.005$, with the linear weights $(d_0, d_1, d_2) = (0.1, 0.6, 0.3)$, for WENO-JS, WENO-M, WENO-Z and WENO-ZR at the second stage \eqref{eq:RK2nd} of the single TVD Runge-Kutta time step for the Riemann problem \eqref{eq:weight_comparison} with $\vez = \cm 2.776 \e \cm 17$ and $\vet = 5.551\e \cm 17$.}
\label{tab:weight_comparison_second}

\vspace{0.2in}

\begin{tabular}{cccccccccccc} 
\hline
$x$ & -0.03 & -0.02 & -0.01 & 0 & 0.01 & 0.02 & 0.03 & 0.04 & 0.05 & 0.06 & 0.07 \\ 
\hline	
$\hat{f}^\JS_j$ & 1 & 1 & 1 & 1-2.220e-16 & 0.244104 & -2.864e-24  & 1.866e-25  & -6.146e-27  & 3.472e-28  & 0 & 0 \\  
$\hat{f}^\M_j$  & 1 & 1 & 1 & 1-2.220e-16 & 0.227637 & -1.873e-79  & 2.176e-80  & -4.469e-82  & 4.237e-83  & 0 & 0 \\  
$\hat{f}^\Z_j$  & 1 & 1 & 1 & 1-2.220e-16 & 0.208333 & -1.906e-40  & 1.614e-41  & -1.556e-42  & 9.259e-44  & 0 & 0 \\
$\hat{f}^\ZR_j$ & 1 & 1 & 1 & 1-2.220e-16 & 0.208333 & -1.906e-120 & 1.614e-121 & -1.556e-122 & 9.259e-124 & 0 & 0 \\ 
\hline
\end{tabular}
\caption{The numerical fluxes $\hat{f}^\N_j$ near the discontinuity at $x = 0.005$ for WENO-JS, WENO-M, WENO-Z and WENO-ZR at the second stage \eqref{eq:RK2nd} of the single TVD Runge-Kutta time step for the Riemann problem \eqref{eq:weight_comparison}.}
\label{tab:numerical_flux_second}

\vspace{0.2in}

\begin{tabular}{ccccccccccc} 
\hline
$x$ & -0.025 & -0.015 & -0.005 & 0.005 & 0.015 & 0.025 & 0.035 & 0.045 & 0.055 & 0.065 \\ 
\hline	
$\ubar^\JS_j$ & 1 & 1 & 1 & 0.219487 & 0.030513 & -3.787e-25  & 2.409e-26  & -8.116e-28  & 4.340e-29  & 0 \\  
$\ubar^\M_j$  & 1 & 1 & 1 & 0.221545 & 0.028455 & -2.582e-80  & 2.775e-81  & -6.116e-83  & 5.296e-84  & 0 \\  
$\ubar^\Z_j$  & 1 & 1 & 1 & 0.223958 & 0.026042 & -2.514e-41  & 2.212e-42  & -2.060e-43  & 1.157e-44  & 0 \\
$\ubar^\ZR_j$ & 1 & 1 & 1 & 0.223958 & 0.026042 & -2.514e-121 & 2.212e-122 & -2.060e-123 & 1.157e-124 & 0 \\
$\ubar_j$     & 1 & 1 & 1 & 0.5      & 0        & 0           & 0          & 0           & 0          & 0 \\  
\hline
\end{tabular}
\caption{The numerical solutions $\ubar^\N_j$, compared with the exact solution $\ubar_j$, near the discontinuity at $x = 0.005$ for WENO-JS, WENO-M, WENO-Z and WENO-ZR at the second stage \eqref{eq:RK2nd} of the single TVD Runge-Kutta time step for the Riemann problem \eqref{eq:weight_comparison}.}
\label{tab:numerical_solution_second}

\end{sidewaystable}

\subsubsection{Third stage of the third-order TVD Runge-Kutta method} \label{sec:RK3rd}
We will go through the numerical solutions in the third stage to distinguish the large error terms and reveal the main source of error.
Consider the final stage of the TVD Runge-Kutta time step,
\begin{equation} \label{eq:RK3rd}
 u^{n+1} = \frac{1}{3} u^n + \frac{2}{3} u^{(2)} + \frac{2}{3} \Delta t L \left( u^{(2)} \right).
\end{equation} 
From the observation in Section \ref{sec:RK2nd}, 12 stencils are to be studied:
\begin{enumerate}[label=\arabic*., font=\itshape]
\setcounter{enumi}{15}
\item $\{ I_{j-4}, I_{j-3}, I_{j-2}, I_{j-1}, I_j \},~j \leqslant -1$;
\item $\{ I_{-4}, I_{-3}, I_{-2}, I_{-1}, I_0 \}$; \label{item:stt3}
\item $\{ I_{-3}, I_{-2}, I_{-1}, I_0, I_1 \}$; 
\item $\{ I_{-2}, I_{-1}, I_0, I_1, I_2 \}$;
\item $\{ I_{-1}, I_0, I_1, I_2, I_3 \}$;
\item $\{ I_0, I_1, I_2, I_3, I_4 \}$;             \label{item:end3}
\item $\{ I_1, I_2, I_3, I_4, I_5 \}$;
\item $\{ I_2, I_3, I_4, I_5, I_6 \}$;
\item $\{ I_3, I_4, I_5, I_6, I_7 \}$;
\item $\{ I_4, I_5, I_6, I_7, I_8 \}$;
\item $\{ I_5, I_6, I_7, I_8, I_9 \}$;
\item $\{ I_j, I_{j+1}, I_{j+2}, I_{j+3}, I_{j+4} \},~j \geqslant 6$,
\end{enumerate}
where Stencils \ref{case:17} - \ref{case:21} contain the discontinuity.
Again, the analysis performed on each stencil will use the assumptions listed as below, based on the observation in Section \ref{sec:RK2nd}:
\begin{itemize}
\item $e^{\Nc 2}_0 < 0,~e^{\Nc 2}_1 > 0,~| e^{\Nc 2}_0/e^{\Nc 2}_1 | > 10$ for $\N \in \{ \JS, \M, \Z, \ZR \}$, \\
      $| e^{\JS,\,2}_j | > | e^{\M,\,2}_j |,~| e^{\JS,\,2}_j | > | e^{\Z,\,2}_j | > | e^{\ZR,\,2}_j |$ for $j=0,1$;
\item $e^{\Nc 2}_2 < 0,~e^{\Nc 2}_3 > 0,~e^{\Nc 2}_4 < 0,~e^{\Nc 2}_5 > 0$, \\
      $| e^{\JS,\,2}_j | \ll \ceps,~| e^{\M,\,2}_j |, | e^{\ZR,\,2}_j | \ll \teps$ and $| e^{\Z,\,2}_j |$ is on the order of $\teps$ for $j = 2, 3, 4, 5$.
\end{itemize} 

\begin{stencil} \label{case:16}
$\{ I_{j-4}, I_{j-3}, I_{j-2}, I_{j-1}, I_j \},~j \leqslant -1$ \\
We arrive at the same conclusion corresponding to $x=-0.03, -0.02$ in Table \ref{tab:weight_comparison_third} as in Stencil \ref{case:1}.
\end{stencil}

\begin{stencil} \label{case:17}
$\{ I_{-4}, I_{-3}, I_{-2}, I_{-1}, I_0 \}$ \\
As we have $\beta^\N_0 = \beta^\N_1 = 0,~\beta^\N_2 = \frac{4}{3} \left[ ( 1-\nu) \delta - e^\N_0 \right]^2$ for $\N \in \{ \JS, \M, \Z, \ZR \}$, similar analysis as in Stencil \ref{case:2} gives
\begin{gather*}
 \omega^\N_0 \approx \frac{1}{7},~\omega^\N_1 \approx \frac{6}{7},~\omega^\N_2 \gtrapprox 0,~\N \in \{ \JS, \Z, \ZR \}, \\
 \omega^\M_0 \approx \frac{1573}{12361},~ \omega^\M_1 \approx \frac{10788}{12361},~ \omega^\M_2 \gtrapprox 0, \\
 \omega^\Z_2 \gg \omega^\ZR_2,
\end{gather*}
corresponding to $x=-0.01$ in Table \ref{tab:weight_comparison_third}.
\end{stencil}

\begin{stencil} \label{case:18}
$\{ I_{-3}, I_{-2}, I_{-1}, I_0, I_1 \}$ \\
With the smoothness indicators
\begin{align*}
 \beta^\N_0 &= 0,~\beta^\N_1 = \frac{4}{3} \left[ (1-\nu) \delta - e^\N_0 \right]^2, \\
 \beta^\N_2 &= \frac{13}{12} \left[ (1-2\nu) \delta - 2e^\N_0 + e^\N_1 \right]^2 + 
               \frac{1}{4}   \left[ (3-4\nu) \delta - 4e^\N_0 + e^\N_1 \right]^2
\end{align*} 
for $\N \in \{ \JS, \M, \Z, \ZR \}$, the results that corresponds to $x=0$ in Table \ref{tab:weight_comparison_third}, are consistent with the ones in Stencil \ref{case:3}.
\end{stencil}

\begin{stencil} \label{case:19}
$\{ I_{-2}, I_{-1}, I_0, I_1, I_2 \}$ \\
The smoothness indicators satisfy
\begin{align*}
 \beta^\N_0 &= \frac{10}{3} \left[ (1-\nu) \delta - e^\N_0 \right]^2,~
 \beta^\N_1 = \frac{13}{12} \left[ (1-2\nu) \delta - 2e^\N_0 + e^\N_1 \right]^2 + 
               \frac{1}{4}   \left( \delta - e^\N_1 \right)^2, \\
 \beta^\N_2 &= \frac{13}{12} \left( \nu \delta + e^\N_0 - 2e^\N_1 + e^\N_2 \right)^2 + 
               \frac{1}{4}   \left( 3\nu \delta + 3e^\N_0 - 4e^\N_1 + e^\N_2 \right)^2,
\end{align*} 
for $\N \in \{\JS, \M, \Z, \ZR \}$.
With some algebra, we have
\begin{align*}
 \beta^\N_0 &= \frac{13}{12} \left[ (1-\nu) \delta - e^\N_0 \right]^2 +
               \frac{9}{4}   \left( \delta - \nu \delta - e^\N_0 \right)^2, \\
 \beta^\N_1 &= \frac{13}{12} \left[ (1-\nu) \delta - e^\N_0 - (\nu \delta + e^\N_0 - e^\N_1) \right]^2 + 
               \frac{1}{4}   \left( \delta - e^\N_1 \right)^2, \\
 \beta^\N_2 &\approx \frac{13}{12} \left[ (\nu \delta + e^\N_0 - e^\N_1) - e^\N_1 \right]^2 + 
                     \frac{1}{4}   \left[ 3(\nu \delta + e^\N_0 - e^\N_1) - e^\N_1 \right]^2.
\end{align*}
From \eqref{eq:e2b0}, \eqref{eq:e2b1} and \eqref{eq:e2b01}, $\beta^\N_0 > \beta^\N_1 > \beta^\N_2$.
Then, by \eqref{eq:nw_JS1}, \eqref{eq:nw_Z1} and \eqref{eq:nw_ZR1}, $\omega^\N_0 < d_0$ and $\omega^\N_2 > d_2$.
The graphs of $\omega^\N_s,~s = 0,1,2$ versus $\nu$ are shown in Figure \ref{fig:weight25}, from which
$$
   \omega^\JS_s < \omega^\M_s < \omega^\Z_s < \omega^\ZR_s < d_s,~s = 0,1,~~ \omega^\JS_2 > \omega^\M_2 > \omega^\Z_2 > \omega^\ZR_2 > d_2.
$$
The weights $\omega^\N_s$, combined with previous errors $e^{\Nc 2}_0$ and $e^{\Nc 2}_1$, will appear as important components in the errors $\es^{\Nc 1}_0$ \eqref{eq:e1b0} and $\es^{\Nc 1}_1$ \eqref{eq:e1b1} below.
\begin{figure}[htbp]
\centering
\includegraphics[width=0.325\textwidth]{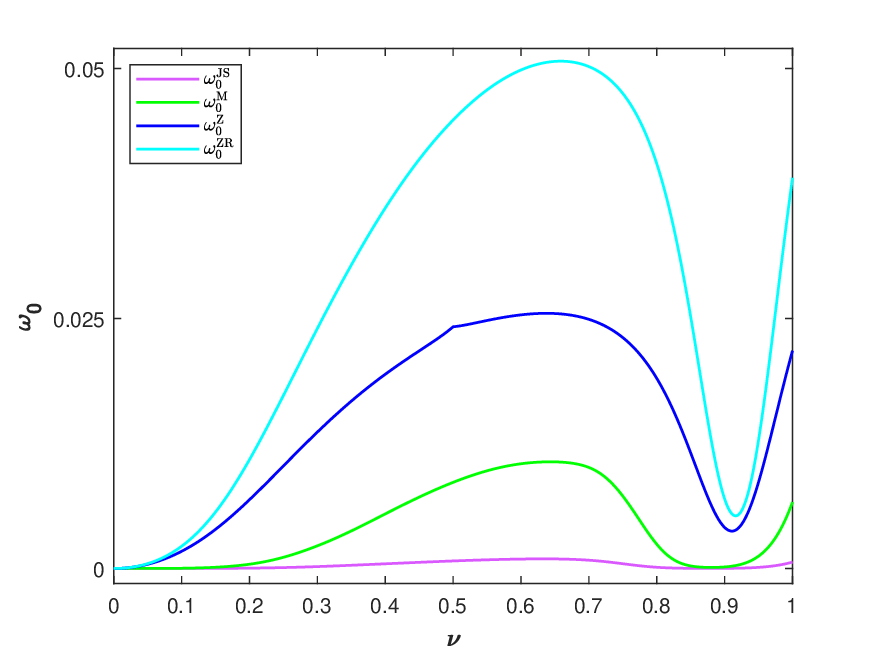}
\includegraphics[width=0.325\textwidth]{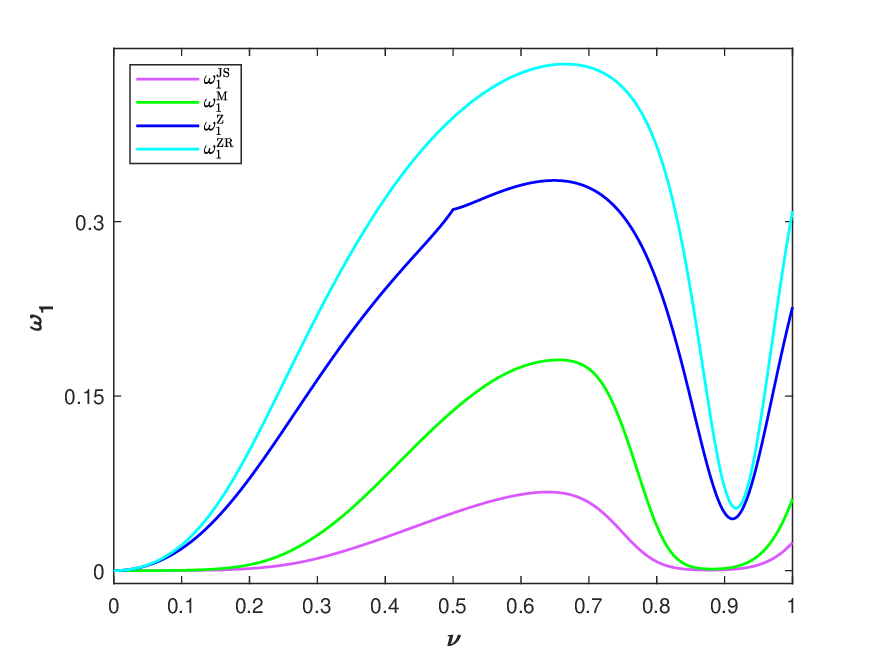}
\includegraphics[width=0.325\textwidth]{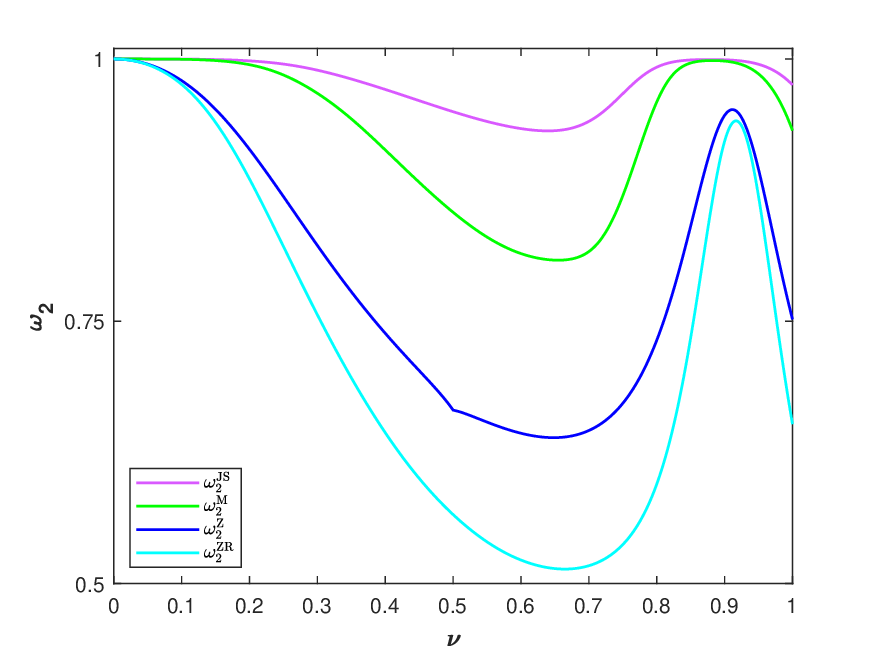}
\vspace{-0.3cm}
\caption{The nonlinear weights $\omega^\N_0$ (left), $\omega^\N_1$ (middle) and $\omega^\N_2$ (right) in Stencil \ref{case:19} by different WENO techniques.}
\label{fig:weight25}
\end{figure}
\end{stencil}

\begin{stencil} \label{case:20}
$\{ I_{-1}, I_0, I_1, I_2, I_3 \}$ \\
The smoothness indicators are
\begin{align*}
 \beta^\N_0 &= \frac{13}{12} \left[ (1-2\nu) \delta - 2e^\N_0 + e^\N_1 \right]^2 + 
               \frac{1}{4}   \left[ (1-4\nu) \delta - 4e^\N_0 + 3e^\N_1 \right]^2, \\
 \beta^\N_1 &= \frac{13}{12} \left( \nu \delta + e^\N_0 - 2e^\N_1 + e^\N_2 \right)^2 + 
               \frac{1}{4}   \left( \nu \delta + e^\N_0 - e^\N_2 \right)^2, \\
 \beta^\N_2 &= \frac{13}{12} \left( e^\N_1 - 2e^\N_2 + e^\N_3 \right)^2 + 
               \frac{1}{4}   \left( 3e^\N_1 - 4e^\N_2 + e^\N_3 \right)^2,
\end{align*} 
for $\N \in \{\JS, \M, \Z, \ZR \}$. 
Since
\begin{align*}
 \beta^\N_0 &= \frac{13}{12} \left[ \delta - e^\N_1 + 2 (\nu \delta + e^\N_0 - e^\N_1) \right]^2 +
               \frac{1}{4}   \left[ \delta - e^\N_1 + 4 (\nu \delta + e^\N_0 - e^\N_1) \right]^2, \\
 \beta^\N_1 &\approx \frac{13}{12} \left[ (\nu \delta + e^\N_0 - e^\N_1) - e^\N_1 \right]^2 + 
                     \frac{1}{4}   \left( \nu \delta + e^\N_0 \right)^2, \\
 \beta^\N_2 &\approx \frac{13}{12} \left( e^\N_1 \right)^2 + \frac{1}{4} \left( 3e^\N_1 \right)^2,
\end{align*}
we have $\beta^\N_0 > \beta^\N_1, \beta^\N_2$ by \eqref{eq:e2b0}, \eqref{eq:e2b1} and \eqref{eq:e2b01}.
Thus $\omega^\N_0 < d_0$ by \eqref{eq:nw_JS1}, \eqref{eq:nw_Z1} and \eqref{eq:nw_ZR1}.
Figure \ref{fig:weight26} shows the graphs of $\omega^\N_s,~s = 0,1,2$ versus $\nu$.
The weights $\omega^\N_s$, as well as previous errors $e^{\Nc 2}_0$ and $e^{\Nc 2}_1$, are important components in the errors $\es^{\Nc 1}_1$ \eqref{eq:e1b1} and $\es^{\Nc 1}_2$ \eqref{eq:e1b2} below.
If we keep $\nu \in [0.4, 0.5]$ fixed, then $\omega^\JS_0 < \omega^\M_0 < \omega^\Z_0 < \omega^\ZR_0$ as shown in Figure \ref{fig:weight26}.  
\begin{figure}[htbp]
\centering
\includegraphics[width=0.325\textwidth]{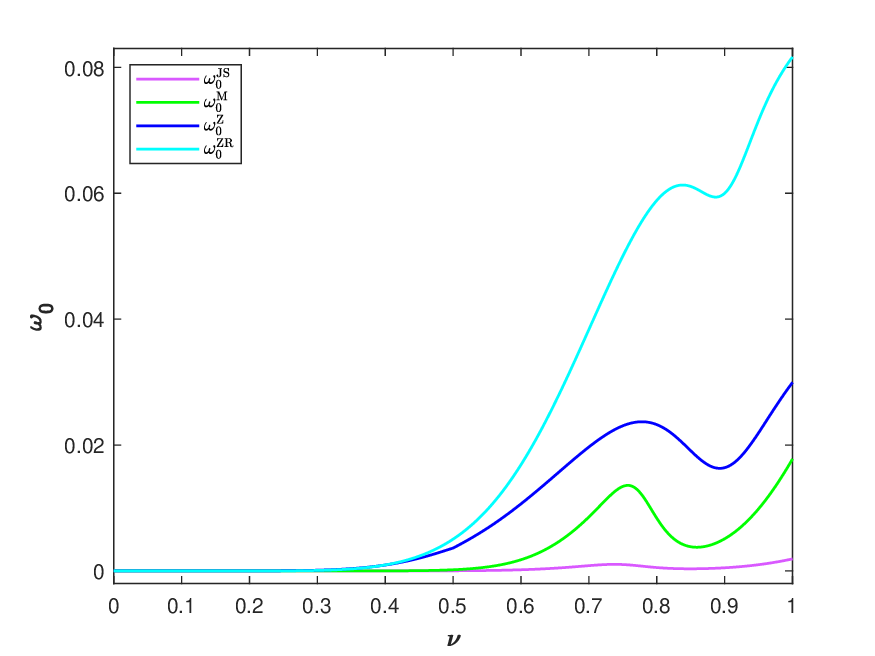}
\includegraphics[width=0.325\textwidth]{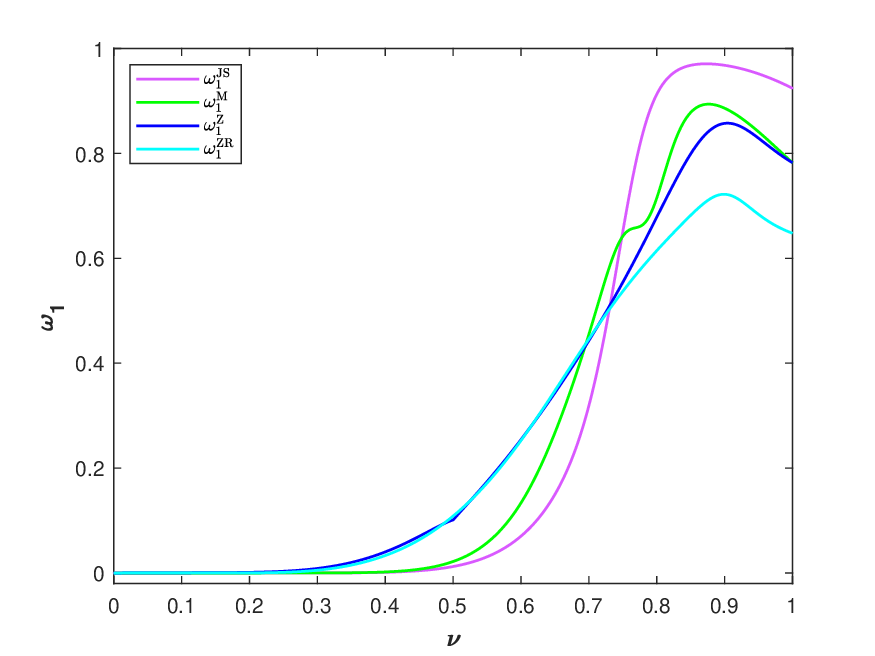}
\includegraphics[width=0.325\textwidth]{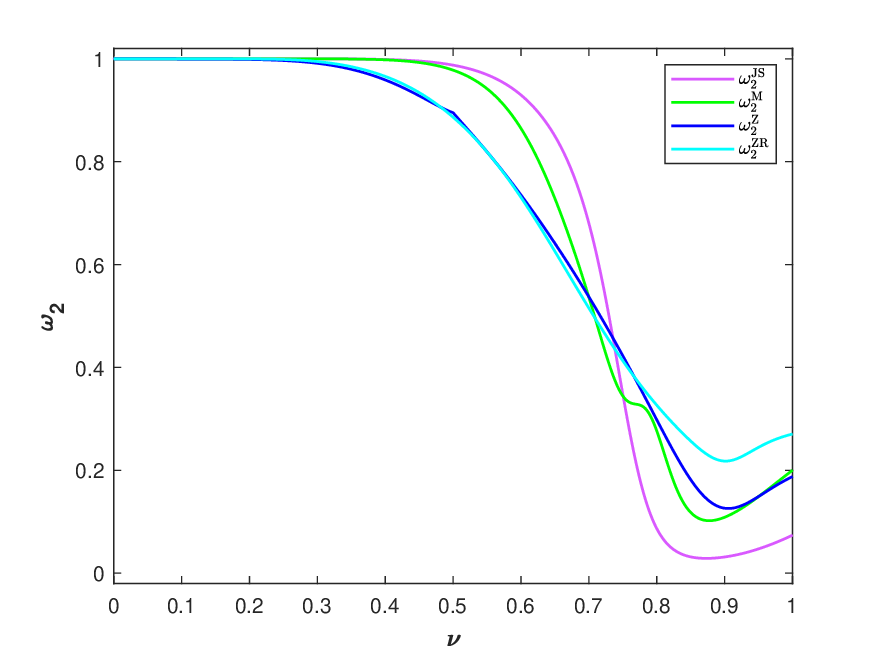}
\vspace{-0.3cm}
\caption{The nonlinear weights $\omega^\N_0$ (left), $\omega^\N_1$ (middle) and $\omega^\N_2$ (right) in Stencil \ref{case:20} by different WENO techniques.}
\label{fig:weight26}
\end{figure}
\end{stencil}

\begin{stencil} \label{case:21}
$\{ I_0, I_1, I_2, I_3, I_4 \}$ \\
From
\begin{align*}
 \beta^\N_0 &= \frac{13}{12} \left( \nu \delta + e^\N_0 - 2e^\N_1 + e^\N_2 \right)^2 + 
               \frac{1}{4}   \left( \nu \delta + e^\N_0 - 4e^\N_1 + 3e^\N_2 \right)^2, \\ 
 \beta^\N_1 &= \frac{13}{12} \left( e^\N_1 - 2e^\N_2 + e^\N_3 \right)^2 + 
               \frac{1}{4}   \left( e^\N_1 - e^\N_3 \right)^2, \\
 \beta^\N_2 &= \frac{13}{12} \left( e^\N_2 - 2e^\N_3 + e^\N_4 \right)^2 + 
               \frac{1}{4}   \left( 3e^\N_2 - 4e^\N_3 + e^\N_4 \right)^2.
\end{align*} 
for $\N \in \{ \JS, \M, \Z, \ZR \}$, it is easy to check that $\beta^\N_0, \beta^\N_1 \gg \beta^\N_2$.
Then we have the same conclusion as in Stencil \ref{case:11}, corresponding to $x=0.03$ in Table \ref{tab:weight_comparison_third}.
\end{stencil}

\begin{stencil} \label{case:22}
$\{ I_1, I_2, I_3, I_4, I_5 \}$ \\
Note that for $\N \in \{ \JS, \M, \Z, \ZR \}$,
\begin{align*}
 \beta^\N_0 &= \frac{13}{12} \left( e^\N_1 - 2e^\N_2 + e^\N_3 \right)^2 + 
               \frac{1}{4}   \left( e^\N_1 - 4e^\N_2 + 3e^\N_3 \right)^2, \\
 \beta^\N_1 &= \frac{13}{12} \left( e^\N_2 - 2e^\N_3 + e^\N_4 \right)^2 + 
               \frac{1}{4}   \left( e^\N_2 - e^\N_4 \right)^2, \\
 \beta^\N_2 &= \frac{13}{12} \left( e^\N_3 - 2e^\N_4 + e^\N_5 \right)^2 + 
               \frac{1}{4}   \left( 3e^\N_3 - 4e^\N_4 + e^\N_5 \right)^2.
\end{align*} 
The results in correspondence to $x=0.04$ in Table \ref{tab:weight_comparison_third}, are in agreement with the ones in Stencil \ref{case:12} .
\end{stencil}

\begin{stencil} \label{case:23}
$\{ I_2, I_3, I_4, I_5, I_6 \}$ \\
With the smoothness indicators 
\begin{align*}
 \beta^\N_0 &= \frac{13}{12} \left( e^\N_2 - 2e^\N_3 + e^\N_4 \right)^2 + 
               \frac{1}{4}   \left( e^\N_2 - 4e^\N_3 + 3e^\N_4 \right)^2, \\
 \beta^\N_1 &= \frac{13}{12} \left( e^\N_3 - 2e^\N_4 + e^\N_5 \right)^2 + 
               \frac{1}{4}   \left( e^\N_3 - e^\N_5 \right)^2, \\
 \beta^\N_2 &= \frac{13}{12} \left( e^\N_4 - 2e^\N_5 \right)^2 + \frac{1}{4} \left( 3e^\N_4 - 4e^\N_5 \right)^2,
\end{align*} 
for $\N \in \{ \JS, \M, \Z, \ZR \}$, we have $\cbeta^\JS_1, \cbeta^\JS_2 \approx 1$, and hence
$$
   \omega^\JS_s \approx d_s,~s=0,1,2.
$$
A similar analysis shows that
\begin{gather*}
 \omega^\N_s \approx d_s,~s = 0,1,2,~\N \in \{ \M, \Z, \ZR \},
\end{gather*}
at $x=0.05$ in Table \ref{tab:weight_comparison_third}.
\end{stencil}

\begin{stencil} \label{case:24}
$\{ I_3, I_4, I_5, I_6, I_7 \}$ \\
Given
\begin{align*}
 \beta^\N_0 &= \frac{13}{12} \left( e^\N_3 - 2e^\N_4 + e^\N_5 \right)^2 + 
               \frac{1}{4}   \left( e^\N_3 - 4e^\N_4 + 3e^\N_5 \right)^2, \\
 \beta^\N_1 &= \frac{13}{12} \left( e^\N_4 - 2e^\N_5 \right)^2 + \frac{1}{4} \left( e^\N_4 \right)^2,~
 \beta^\N_2  = \frac{10}{3} \left( e^\N_5 \right)^2
\end{align*} 
for $\N \in \{ \JS, \M, \Z, \ZR \}$, we have the same conclusion as in Stencil \ref{case:23}, which corresponds to $x=0.06$ in Table \ref{tab:weight_comparison_third}.
\end{stencil}

\begin{stencil} \label{case:25}
$\{ I_4, I_5, I_6, I_7, I_8 \}$ \\
The smoothness indicators are
\begin{align*}
 \beta^\N_0 &= \frac{13}{12} \left( e^\N_4 - 2e^\N_5 \right)^2 + \frac{1}{4} \left( e^\N_4-4e^\N_5 \right)^2, \\
 \beta^\N_1 &= \frac{4}{3} \left( e^\N_5 \right)^2,~\beta^\N_2 = 0,
\end{align*} 
for $\N \in \{ \JS, \M, \Z, \ZR \}$. 
The results, which is associated with $x=0.07$ in Table \ref{tab:weight_comparison_third}, agree with the ones in Stencil \ref{case:13}.
\end{stencil}

\begin{stencil} \label{case:26}
$\{ I_5, I_6, I_7, I_8, I_9 \}$ \\
Since the smoothness indicators satisfy $\beta^\N_0 = \frac{4}{3} \left( e^\N_5 \right)^2$ and $\beta^\N_1 = \beta^\N_2 = 0$ for $\N \in \{ \JS, \M, \Z, \ZR \}$, we obtain the same conclusion as in Stencil \ref{case:14}, corresponding to $x=0.08$ in Table \ref{tab:weight_comparison_third}.
\end{stencil}

\begin{stencil} \label{case:27}
$\{ I_j, I_{j+1}, I_{j+2}, I_{j+3}, I_{j+4} \},~j \geqslant 6$ \\
We get the same conclusion at $x=0.09, 0.1$ in Table \ref{tab:weight_comparison_third} as in Stencil \ref{case:1}.
\end{stencil}

\begin{sidewaystable}
\renewcommand{\arraystretch}{1.2}
\centering
\vspace{4.6in}
\begin{tabular}{cccccccc} 
\hline
$x$ & -0.03 & -0.02 & -0.01 & 0 & 0.01 & 0.02 & 0.03 \\ 
\hline 
$\omega^\JS_0$ & 0.1        & 0.1        & 0.142857 & 1 & 0.000755 & 2.080e-5 & 3.811e-22 \\  
$\omega^\M_0$  & 0.1+$\vez$ & 0.1+$\vez$ & 0.127255 & 1 & 0.008614 & 1.817e-4 & 3.522e-77  \\  
$\omega^\Z_0$  & 0.1        & 0.1        & 0.142857 & 1 & 0.024174 & 3.631e-3 & 1.873e-39  \\
$\omega^\ZR_0$ & 0.1        & 0.1        & 0.142857 & 1 & 0.044847 & 5.050e-3 & 1.873e-119 \\  
\hline 
$\omega^\JS_1$ & 0.6 & 0.6 & 0.857143 & 9.094e-24  & 0.049509 & 0.012352 & 1.298e-18  \\  
$\omega^\M_1$  & 0.6 & 0.6 & 0.872745 & 2.451e-79  & 0.137537 & 0.021935 & 4.576e-74  \\  
$\omega^\Z_1$  & 0.6 & 0.6 & 0.857143 & 1.148e-39  & 0.310561 & 0.101230 & 2.268e-37  \\
$\omega^\ZR_1$ & 0.6 & 0.6 & 0.857143 & 1.148e-119 & 0.389471 & 0.108576 & 2.268e-117 \\ 
\hline 
$\omega^\JS_2$ & 0.3        & 0.3        & 6.496e-25  & 1.269e-24  & 0.949737 & 0.987628 & 1 \\  
$\omega^\M_2$  & 0.3+$\vet$ & 0.3+$\vet$ & 3.597e-80 & 5.638e-80 & 0.853848 & 0.977883 & 1 \\  
$\omega^\Z_2$  & 0.3        & 0.3        & 1.067e-40  & 4.009e-40  & 0.665265 & 0.895138 & 1 \\
$\omega^\ZR_2$ & 0.3        & 0.3        & 1.067e-120 & 4.009e-120 & 0.565682 & 0.886373 & 1 \\
\hline
$x$ & 0.04 & 0.05 & 0.06 & 0.07 & 0.08 & 0.09 & 0.1 \\ 
\hline 
$\omega^\JS_0$ & 7.210e-20  & 0.1        & 0.1        & 0.1        & 0.1        & 0.1        & 0.1 \\  
$\omega^\M_0$  & 1.163e-74  & 0.1+$\vez$ & 0.1+$\vez$ & 0.1+$\vez$ & 0.1+$\vez$ & 0.1+$\vez$ & 0.1+$\vez$ \\  
$\omega^\Z_0$  & 2.458e-38  & 0.1        & 0.1        & 0.1        & 0.1        & 0.1        & 0.1 \\
$\omega^\ZR_0$ & 2.458e-118 & 0.1        & 0.1        & 0.1        & 0.1        & 0.1        & 0.1 \\ 
\hline 
$\omega^\JS_1$ & 0.666667 & 0.6 & 0.6 & 0.6 & 0.6 & 0.6 & 0.6 \\  
$\omega^\M_1$  & 0.667027 & 0.6 & 0.6 & 0.6 & 0.6 & 0.6 & 0.6 \\  
$\omega^\Z_1$  & 0.666667 & 0.6 & 0.6 & 0.6 & 0.6 & 0.6 & 0.6 \\
$\omega^\ZR_1$ & 0.666667 & 0.6 & 0.6 & 0.6 & 0.6 & 0.6 & 0.6 \\ 
\hline 
$\omega^\JS_2$ & 0.333333 & 0.3        & 0.3        & 0.3        & 0.3        & 0.3        & 0.3 \\  
$\omega^\M_2$  & 0.332973 & 0.3+$\vet$ & 0.3+$\vet$ & 0.3+$\vet$ & 0.3+$\vet$ & 0.3+$\vet$ & 0.3+$\vet$ \\  
$\omega^\Z_2$  & 0.333333 & 0.3        & 0.3        & 0.3        & 0.3        & 0.3        & 0.3 \\
$\omega^\ZR_2$ & 0.333333 & 0.3        & 0.3        & 0.3        & 0.3        & 0.3        & 0.3 \\ 
\hline
\end{tabular}
\caption{The nonlinear weights $\omega^\N_s$ near the discontinuity at $x = \Delta t$, with the linear weights $(d_0, d_1, d_2) = (0.1, 0.6, 0.3)$, for WENO-JS, WENO-M, WENO-Z and WENO-ZR at the third stage \eqref{eq:RK3rd} of the single TVD Runge-Kutta time step for the Riemann problem \eqref{eq:weight_comparison} with $\vez = \cm 2.776 \e \cm 17$ and $\vet = 5.551\e \cm 17$.}
\label{tab:weight_comparison_third}
\end{sidewaystable}

By \eqref{eq:3rd_reconstruction_minus} and \eqref{eq:weno_appr}, we have the numerical fluxes:
\begin{enumerate}[label=\arabic*., font=\itshape]
\setcounter{enumi}{15}
\item $\hat{f}^{\Nc 3}_{j-1/2} = u_L,~j \leqslant -2$; 
\item $\hat{f}^{\Nc 3}_{-3/2} = u_L + \frac{1}{6} \omega^{\Nc 3}_{2,-3/2} (1-\nu) \delta - \frac{1}{6} \omega^{\Nc 3}_{2,-3/2} e^{\Nc 2}_0$;
\item $\hat{f}^{\Nc 3}_{-1/2} = u_L + \frac{1}{6} \omega^{\Nc 3}_{2,-1/2} \delta - \frac{1}{6} \C^{\Nc 3}_{-1/2} (1-\nu) \delta + \frac{1}{6} \left( \C^{\Nc 3}_{-1/2} e^{\Nc 2}_0 - \omega^{\Nc 3}_{2,-1/2} e^{\Nc 2}_1 \right)$;
\item $\hat{f}^{\Nc 3}_{1/2} = u_L - \frac{1}{3} \A^{\Nc 3}_{1/2} \delta - \frac{1}{6} \D^{\Nc 3}_{1/2} (1-\nu) \delta + \frac{1}{6} \left( \D^{\Nc 3}_{1/2} e^{\Nc 2}_0 + \C^{\Nc 3}_{1/2} e^{\Nc 2}_1 - \omega^{\Nc 3}_{2,1/2} e^{\Nc 2}_2 \right)$;
\item $\hat{f}^{\Nc 3}_{3/2} = u_R - \frac{1}{6} \B^{\Nc 3}_{3/2} \delta + \frac{1}{6} \E^{\Nc 3}_{3/2} (1-\nu) \delta$ \\
      \hspace*{45pt} $+ \frac{1}{6} \left( - \E^{\Nc 3}_{3/2} e^{\Nc 2}_0 + \D^{\Nc 3}_{3/2} e^{\Nc 2}_1 + \C^{\Nc 3}_{3/2} e^{\Nc 2}_2 - \omega^{\Nc 3}_{2,3/2} e^{\Nc 2}_3 \right)$; 
\item $\hat{f}^{\Nc 3}_{5/2} = u_R + \frac{1}{3} \omega^{\Nc 3}_{0,5/2} \nu \delta$ \\
      \hspace*{45pt} $+ \frac{1}{6} \left( 2 \omega^{\Nc 3}_{0,5/2} e^{\Nc 2}_0 - \E^{\Nc 3}_{5/2} e^{\Nc 2}_1 + \D^{\Nc 3}_{5/2} e^{\Nc 2}_2 + \C^{\Nc 3}_{5/2} e^{\Nc 2}_3 - \omega^{\Nc 3}_{2,5/2} e^{\Nc 2}_4 \right)$;
\item $\hat{f}^{\Nc 3}_{7/2} = u_R + \frac{1}{6} \left( 2 \omega^{\Nc 3}_{0,7/2} e^{\Nc 2}_1 - \E^{\Nc 3}_{7/2} e^{\Nc 2}_2 + \D^{\Nc 3}_{7/2} e^{\Nc 2}_3 + \C^{\Nc 3}_{7/2} e^{\Nc 2}_4 - \omega^{\Nc 3}_{2,7/2} e^{\Nc 2}_5 \right)$;
\item $\hat{f}^{\Nc 3}_{9/2} = u_R + \frac{1}{6} \left( 2 \omega^{\Nc 3}_{0,9/2} e^{\Nc 2}_2 - \E^{\Nc 3}_{9/2} e^{\Nc 2}_3 + \D^{\Nc 3}_{9/2} e^{\Nc 2}_4 + \C^{\Nc 3}_{9/2} e^{\Nc 2}_5 \right)$;
\item $\hat{f}^{\Nc 3}_{11/2} = u_R + \frac{1}{6} \left( 2 \omega^{\Nc 3}_{0,11/2} e^{\Nc 2}_3 - \E^{\Nc 3}_{11/2} e^{\Nc 2}_4 + \D^{\Nc 3}_{13/2} e^{\Nc 2}_5 \right)$;
\item $\hat{f}^{\Nc 3}_{13/2} = u_R + \frac{1}{6} \left( 2 \omega^{\Nc 3}_{0,13/2} e^{\Nc 2}_4 - \E^{\Nc 3}_{13/2} e^{\Nc 2}_5 \right)$;
\item $\hat{f}^{\Nc 3}_{15/2} = u_R + \frac{1}{3} \omega^{\Nc 3}_{0,15/2} e^{\Nc 2}_5$;
\item $\hat{f}^{\Nc 3}_{j+1/2} = u_R,~j \geqslant 8$,
\end{enumerate}
for  $\N \in \{ \JS, \M, \Z, \ZR \}$ and $\A, \B, \C, \D, \E$ given in Appendix.
Similar to the computation of the numerical flux in Sections \ref{sec:Euler} and \ref{sec:RK2nd}, the numerical fluxes computed with round-off errors in Table \ref{tab:numerical_flux_third} do not agree with the analyzed ones above at $x = -0.01, 0$.

If we use \eqref{eq:RK3rd} for the time integration at the third stage, then the discontinuity will stay at $x = \Delta t$, and the  numerical solutions after one TVD Runge-Kutta time step are as follows:
\begin{equation} \label{eq:numerical_solution_third}
\begin{aligned}
 \ubar^\N_j &= \ubar^1_j,~j \leqslant -3, \\
 \ubar^\N_j &= \ubar^1_j + \es^{\Nc 1}_j,~j = -2, \cdots, 8, \\
 \ubar^\N_j &= \ubar^1_j,~j \geqslant 9,
\end{aligned}
\end{equation}
where the errors $\es^{\Nc 1}_j$ given in Appendix.
It is expected to mach the formulas \eqref{eq:numerical_solution_third} with the numerical solutions in Table \ref{tab:numerical_solution_second} except at $x = -0.015, -0.005$, corresponding to $\ubar^\N_j = 1,~j = -2, -1$ with $\es^{\Nc 1}_j = 0$ for $\N \in \{ \JS, \M, \Z, \ZR \}$, as is the case in Sections \ref{sec:Euler} and \ref{sec:RK2nd}.
Also note that $\es^{\Nc 1}_3, \cdots, \es^{\Nc 1}_8$ are small and less than $\dbp$ according to the formula.

Now we take a look at the main source of errors $\es^{\Nc 1}_0, \es^{\Nc 1}_1$ and $\es^{\Nc 1}_2$:
\begin{align}
 \es^{\Nc 1}_0 \approx & \frac{1}{3} \nu \delta + \frac{1}{9} \left( \B^{\Nc 3}_{1/2} - \D^{\Nc 3}_{1/2} \nu \right) \nu \delta + \frac{\nu}{9} \left[ \left( - \D^{\Nc 3}_{1/2} + 6/\nu \right) e^{\Nc 2}_0 - \C^{\Nc 3}_{1/2} e^{\Nc 2}_1 \right] \label{eq:e1b0} \\
               = & \frac{1}{9} \left[ \left( 8 - 11 \nu \right) \nu \delta + \left( 6 - 11 \nu \right) e^{\Nc 2}_0 \right] +
                   \frac{\nu}{9} \left[ \left( 6\nu-4 \right) \delta + 6e^{\Nc 2}_0 - 2e^{\Nc 2}_1 \right] \omega^{\Nc 3}_{1,1/2}+ \nonumber \\
                 & \frac{\nu}{9} \left[ \left( 9\nu-5 \right) \delta + 9e^{\Nc 2}_0 - 5e^{\Nc 2}_1 \right] \omega^{\Nc 3}_{2,1/2}, \nonumber \\
 \es^{\Nc 1}_1 \approx & \frac{1}{9} \left[ - \B^{\Nc 3}_{1/2} + \left( \D^{\Nc 3}_{1/2} + \omega^{\Nc 3}_{1,3/2} \right) \nu \right] \nu \delta \label{eq:e1b1}  \\
                       & + \frac{\nu}{9} \left[ \left( \D^{\Nc 3}_{1/2} + \omega^{\Nc 3}_{1,3/2} \right) e^{\Nc 2}_0 + \left( \C^{\Nc 3}_{1/2} - 5 \omega^{\Nc 3}_{1,3/2} - 2 \omega^{\Nc 3}_{2,3/2} + 6/\nu \right) e^{\Nc 2}_1 \right] \nonumber \\
               = & \frac{1}{9} \left[ 2 \nu \left( \nu \delta + e^{\Nc 2}_0 \right) + \left( 5\nu+6 \right) e^{\Nc 2}_1 \right] +
                   \frac{\nu}{9} \left[ \left( 9\nu-5 \right) \delta + 9e^{\Nc 2}_0 - 5e^{\Nc 2}_1 \right] \omega^{\Nc 3}_{0,1/2}+ \nonumber \\
                 & \frac{1}{9} \left[ \left( 3\nu-1 \right) \nu \delta + 3 \nu e^{\Nc 2}_0 - 3 \nu e^{\Nc 2}_1 \right] \omega^{\Nc 3}_{1,1/2} +
                   \frac{\nu}{9} \left( \nu \delta + e^{\Nc 2}_0 - 5 e^{\Nc 2}_1 \right) \omega^{\Nc 3}_{1,3/2} - \frac{2}{9} \nu e^{\Nc 2}_1 \omega^{\Nc 3}_{2,3/2}, \nonumber \\                         
 \es^{\Nc 1}_2 \approx & - \frac{1}{9} \omega^{\Nc 3}_{1,3/2} \nu^2 \delta + \frac{\nu}{9} \left[ - \omega^{\Nc 3}_{1,3/2} e^{\Nc 2}_0 + \left( 5 \omega^{\Nc 3}_{1,3/2} + 2 \omega^{\Nc 3}_{2,3/2} \right) e^{\Nc 2}_1 \right] \label{eq:e1b2} \\
               = & \frac{\nu}{9} \left( - \nu \delta - e^{\Nc 2}_0 + 5 e^{\Nc 2}_1 \right) \omega^{\Nc 3}_{1,3/2} + \frac{2}{9} \nu e^{\Nc 2}_1 \omega^{\Nc 3}_{2,3/2} \nonumber.
\end{align}

From a multitude of numerical experiments, we find
\begin{gather*}
 \frac{1}{9} \left[ \left( 3 \nu-1 \right) \nu \delta + 3 \nu e^{\Nc 2}_0 - 3 \nu e^{\Nc 2}_1 \right] < 0,~0 < \delta \leqslant 1,\\
 \frac{\nu}{9} \left( \nu \delta + e^{\Nc 2}_0 - 5 e^{\Nc 2}_1 \right) > 0,~ 0 < \delta \leqslant 0.5.
\end{gather*}

\begin{figure}[htbp]
\centering
\includegraphics[width=0.325\textwidth]{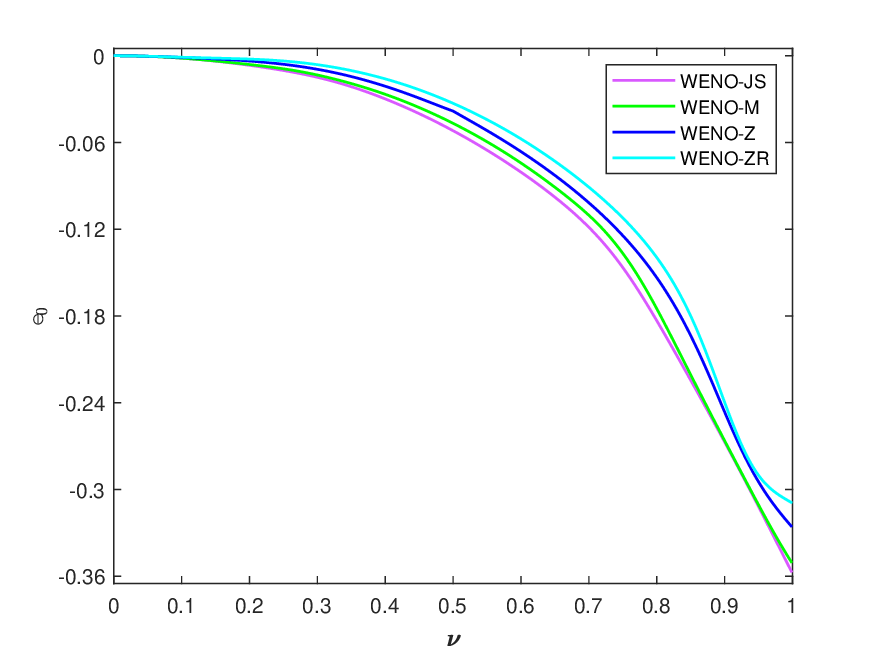}
\includegraphics[width=0.325\textwidth]{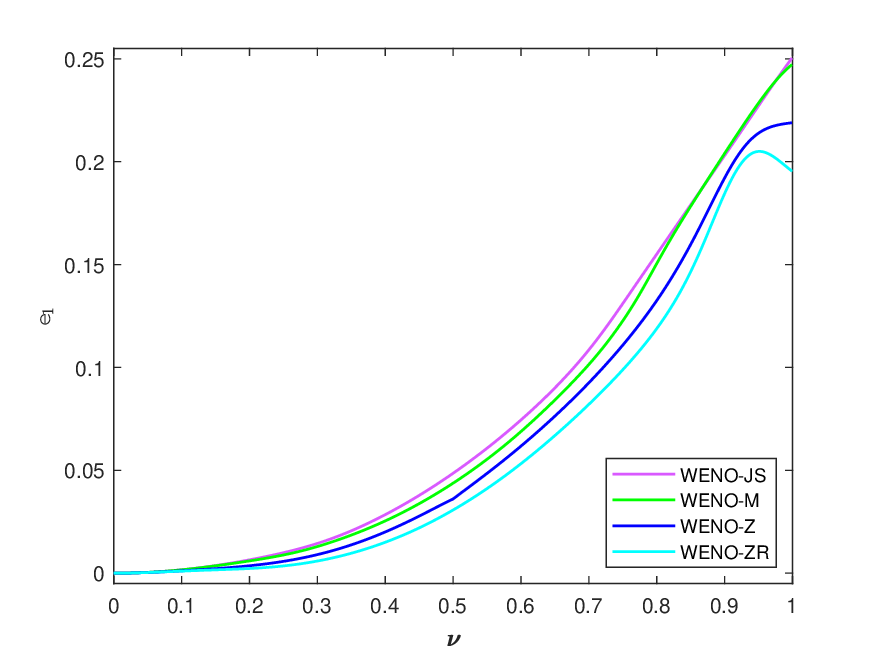}
\includegraphics[width=0.325\textwidth]{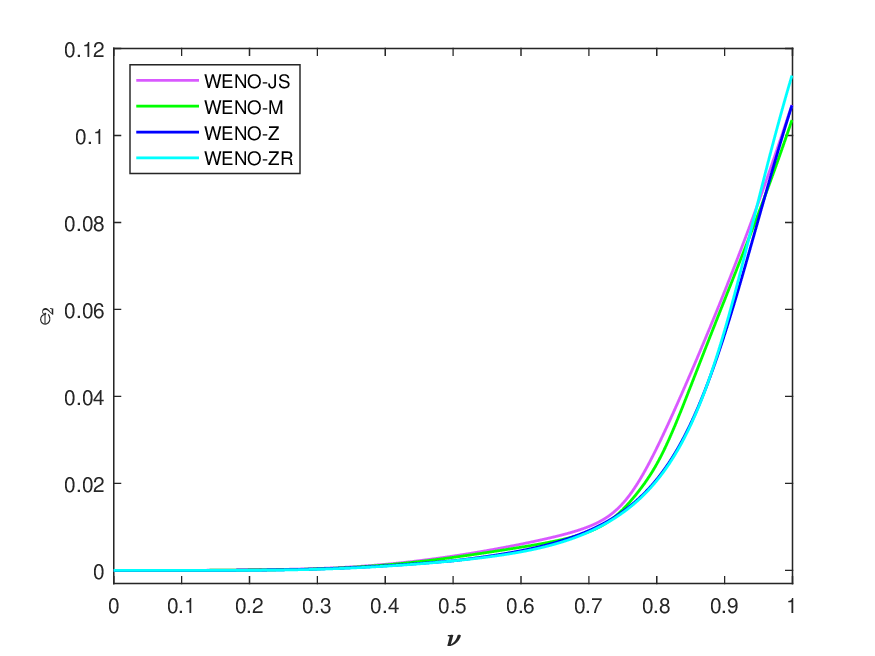}
\vspace{-0.3cm}
\caption{$\es^\N_0$ (left), $\es^\N_1$ (middle) and $\es^\N_2$ (right).}
\label{fig:error3}
\end{figure}

It is easy to see $\es^{\Nc 1}_0 < 0$ from Figure \ref{fig:error3}. 
Then we aim at reducing the values of both $\omega^{\Nc 3}_{1,1/2}$ and $\omega^{\Nc 3}_{2,1/2}$ corresponding to Stencil \ref{case:19}.
In other words, we want the larger value of $\omega^{\Nc 3}_{0,1/2}$ in Figure \ref{fig:weight25}.

Figure \ref{fig:error3} also shows $\es^{\Nc 1}_1 > 0$ and $\es^{\Nc 1}_2 > 0$.
Then we could increase the value of $\omega^{\Nc 3}_{0,1/2}$ in order to decrease the error $\es^{\Nc 1}_1$.

Thus at the final stage, increasing the value of $\omega^{\Nc 3}_{0,1/2}$ gives the smaller errors $\es^{\Nc 1}_0$ and $\es^{\Nc 1}_1$. 
As shown in Figure \ref{fig:weight25}, WENO-ZR yields larger $\omega^{\ZR,\,3}_{0,1/2}$ than WENO-JS, WENO-M and WENO-Z for $0 < \nu \leqslant 0.5$, leading to the better numerical solutions $\ubar^\ZR_0$ and $\ubar^\ZR_1$ in Table \ref{tab:numerical_solution_third} with smaller errors $\es^{\ZR,\,1}_0$ and $\es^{\ZR,\,1}_1$ in Figure \ref{fig:error3} and less numerical dissipation around the jump discontinuity.

\begin{sidewaystable}
\renewcommand{\arraystretch}{1.2}
\centering

\vspace{5.0in}
\begin{tabular}{cccccccc} 
\hline
$x$ & -0.03 & -0.02 & -0.01 & 0 & 0.01 & 0.02 & 0.03  \\ 
\hline
$\hat{f}^\JS_j$ & 1 & 1 & 1 & 1-2.220e-16 & 0.094617 & 0.009910 & -6.586e-21  \\  
$\hat{f}^\M_j$  & 1 & 1 & 1 & 1-2.220e-16 & 0.083396 & 0.009008 & -2.156e-76  \\  
$\hat{f}^\Z_j$  & 1 & 1 & 1 & 1-2.220e-16 & 0.062778 & 0.006624 & -9.080e-40  \\
$\hat{f}^\ZR_j$ & 1 & 1 & 1 & 1-2.220e-16 & 0.046704 & 0.006603 & -9.080e-120 \\
\hline
$x$ & 0.04 & 0.05 & 0.06 & 0.07 & 0.08 & 0.09 & 0.1 \\ 
\hline
$\hat{f}^\JS_j$ & 7.334e-22  & -1.846e-26  & 1.013e-27  & -3.646e-29  & 1.447e-30  & 0 & 0 \\  
$\hat{f}^\M_j$  & 1.103e-76  & -1.507e-81  & 1.099e-82  & -3.186e-84  & 1.765e-85  & 0 & 0 \\  
$\hat{f}^\Z_j$  & 2.175e-40  & -1.474e-42  & 1.274e-43  & -9.375e-45  & 3.858e-46  & 0 & 0 \\
$\hat{f}^\ZR_j$ & 2.175e-120 & -1.474e-122 & 1.274e-123 & -9.375e-125 & 3.858e-126 & 0 & 0 \\
\hline
\end{tabular}
\caption{The numerical fluxes $\hat{f}^\N_j$ near the discontinuity at $x = 0.005$ for WENO-JS, WENO-M, WENO-Z and WENO-ZR at the third stage \eqref{eq:RK3rd} of the single TVD Runge-Kutta time stage for the Riemann problem \eqref{eq:weight_comparison}.}
\label{tab:numerical_flux_third}

\vspace{0.2in}

\begin{tabular}{cccccccc} 
\hline
$x$ & -0.025 & -0.015 & -0.005 & 0.005 & 0.015 & 0.025 & 0.035 \\ 
\hline
$\ubar^\JS_j$ & 1 & 1 & 1 & 0.448119 & 0.048578 & 0.003303 & -2.440e-21  \\  
$\ubar^\M_j$  & 1 & 1 & 1 & 0.453231 & 0.043766 & 0.003003 & -1.086e-76  \\  
$\ubar^\Z_j$  & 1 & 1 & 1 & 0.461713 & 0.036079 & 0.002208 & -3.737e-40  \\
$\ubar^\ZR_j$ & 1 & 1 & 1 & 0.467071 & 0.030728 & 0.002201 & -3.737e-120 \\
$\ubar_j$     & 1 & 1 & 1 & 0.5      & 0        & 0        & 0 \\ 
\hline
$x$ & 0.045 & 0.055 & 0.065 & 0.075 & 0.085 & 0.095 \\
\hline
$\ubar^\JS_j$ & 2.445e-22  & -6.461e-27  & 3.498e-28  & -1.264e-29  & 4.823e-31  & 0 \\  
$\ubar^\M_j$  & 3.678e-77  & -5.356e-82  & 3.770e-83  & -1.121e-84  & 5.885e-86  & 0 \\  
$\ubar^\Z_j$  & 7.285e-41  & -5.260e-43  & 4.561e-44  & -3.254e-45  & 1.286e-46  & 0 \\
$\ubar^\ZR_j$ & 7.285e-121 & -5.260e-123 & 4.561e-124 & -3.254e-125 & 1.286e-126 & 0 \\
$\ubar_j$     & 0          & 0           & 0          & 0           & 0          & 0 \\ 
\hline
\end{tabular}
\caption{The numerical solutions $\ubar^\N_j$, compared with the exact solution $\ubar_j$, near the discontinuity at $x = 0.005$ for WENO-JS, WENO-M, WENO-Z and WENO-ZR at the third stage \eqref{eq:RK3rd} of the single TVD Runge-Kutta time step for the Riemann problem \eqref{eq:weight_comparison}.}
\label{tab:numerical_solution_third}

\vspace{0.2in}

\begin{tabular}{ccccccccc} 
\hline
$x$ & 0.9655 & 0.975 & 0.985 & 0.995 & 1.005 & 1.015 & 1.025 & 1.035 \\ 
\hline
$\ubar^\JS_j$ & 0.957444 & 0.900244 & 0.781778 & 0.602513 & 0.399953 & 0.219345 & 0.098083 & 0.041337 \\  
$\ubar^\M_j$  & 0.978247 & 0.932979 & 0.816058 & 0.618327 & 0.384776 & 0.183072 & 0.063658 & 0.021682 \\  
$\ubar^\Z_j$  & 0.988074 & 0.947798 & 0.830036 & 0.625016 & 0.381100 & 0.170407 & 0.047171 & 0.009916 \\
$\ubar^\ZR_j$ & 0.990215 & 0.952712 & 0.835129 & 0.627611 & 0.379859 & 0.165692 & 0.041284 & 0.007268 \\
$\ubar_j$     & 1        & 1        & 1        & 1        & 0        & 0        & 0        & 0 \\ 
\hline
\end{tabular}
\caption{The numerical solutions $\ubar^\N_j$ at $T=1$, compared with the exact solution $\ubar_j$, near the discontinuity at $x = 1$ for WENO-JS, WENO-M, WENO-Z and WENO-ZR for the Riemann problem \eqref{eq:weight_comparison}.}
\label{tab:numerical_solution_final}

\end{sidewaystable}

We also provide the values of numerical solutions near the discontinuity at the location of $x=1$, at time $T=1$, in Table \ref{tab:numerical_solution_final}.
The errors corresponding to WENO-JS, WENO-M, WENO-Z and WENO-ZR gradually decrease, as is seen from Table \ref{tab:numerical_solution_third}.
In other words, the WENO scheme with more accurate approximation at the first TVD Runge-Kutta time step, will continue its behavior to the final time.

To conclude this subsection, we present the nonlinear weights that can be adjusted to improve the accuracy:
\begin{itemize}
\item In the second stage, the increase of $\omega^{\Nc 2}_{0,1/2}$ corresponding to Stencil \ref{case:10} is an approach to reducing the errors $e^{\Nc 2}_0$ and $e^{\Nc 2}_1$.
\item In the third stage, we can increase $\omega^{\Nc 3}_{0,1/2}$ corresponding to Stencil \ref{case:19} so that the errors $\es^{\Nc 1}_0$ and $\es^{\Nc 1}_1$ decrease.
\end{itemize}

\section{New Z-type nonlinear weights $\omega^{\ZL}_k$} \label{sec:ZL_weight}
We have seen where the main source of error is after the first TVD Runge-Kutta time step in Section \ref{sec:FV_WENO_discontinuity}.
The WENO-ZR scheme gives better numerical solutions as the nonlinear weights $\omega^\ZR_s$ suppress the error around the jump discontinuity. 
However, computing $\omega^\ZR_s$ with the $p$th root of the smoothness indicators $\beta_s$ takes longer time in \cite{Gu} than the other weights $\omega^\N_s,~\N \in \{ \JS, \M, \ZR \}$.
To improve the efficiency, we will propose a novel nonlinear weight formula in this section and explain how the error terms in Section \ref{sec:FV_WENO_discontinuity} are used to verify that the proposed one enhances the performance.

To take all that are discussed in Section \ref{sec:FV_WENO_discontinuity} into account, we introduce a novel global smoothness indicator $\tau^\ZL_5$ on the stencil $S^5$ as
\begin{equation} \label{eq:global_smooth_indicator_ZL}
 \tau^\ZL_5 = \frac{1}{p} \left| \ln \frac{1+\beta_0}{1+\beta_2} \right|. 
\end{equation}
Hence the nonlinear weights denoted by $\omega^\ZL_s$ are defined as
\begin{equation} \label{eq:weights_ZL}
 \omega^\ZL_s = \frac{\alpha_s}{\sum_{r=0}^2 \alpha_r},~ \alpha_s = d_s \left( 1 + \left( \frac{\tau^\ZL_5}{\beta_s + \epsilon} \right)^q \right),~ s=0,1,2,
\end{equation}
with $p>0,\, q \geqslant 1$. 
Note that, in the limit of $p \to \infty$, 
\begin{equation} \label{eq:tau5_limit_infty} 
 \lim_{p \to \infty} \tau^{\ZL}_5  = 0.
\end{equation}
Then we have from \eqref{eq:tau5_limit_infty} 
$$
   \lim_{p \to \infty} \alpha_s = \lim_{p \to \infty} d_s \left( 1 + \left( \frac{\tau^\ZL_5}{\beta_s + \epsilon} \right)^q \right) = d_s,
$$
and hence
\begin{equation} \label{eq:weights_ZL_infty}
 \lim_{p \to \infty} \omega^{\ZL}_s = \lim_{p \to \infty} \frac{\alpha_s}{\sum^2_{r=0} \alpha_r} = d_s,
\end{equation}
where the nonlinear weights return to the linear ones in the stencil.
The parameter $q$ has the following effect on maintaining fifth-order accuracy in smooth regions.
Using Taylor series expansions, we obtain $\omega^\ZL_s = d_s + O(\Delta x^{3q})$ for $v'_i \ne 0$, whereas $\omega^\ZR_s = d_s + O(\Delta x^q)$ for $v'_i = 0$.
The condition \eqref{eq:WENO_condition} can be satisfied provided $q \geqslant 2$.
We will see that introducing the parameters $p$ and $q$ allows more flexibility to reduce numerical dissipation around discontinuities. 

\subsection{Dissipation around discontinuities} \label{sec:FV_WENO_discontinuity_ZL} 
When $\beta^\ZL_0 = \beta^\ZL_1 = 0$, we have 
\begin{equation} \label{eq:nw_ZL0}
 \omega^\ZL_0 = \frac{d_0}{d_0 + d_1 + {d_2 \over \tbeta^\ZL_2}},~           
 \omega^\ZL_1 = \frac{d_1}{d_0 + d_1 + {d_2 \over \tbeta^\ZL_2}},~            
 \omega^\ZL_2 = \frac{d_2}{d_0 \tbeta^\ZL_2 + d_1 \tbeta^\ZL_2 + d_2},
\end{equation} 
where
$$
   \tbeta_2 = \tbeta_2(p,q) = \frac{1 + \left[ \frac{\ln (1+\beta_2)}{p \teps} \right]^q}{1 + \left[ \frac{\ln (1+\beta_2)}{p(\beta_2+\teps)} \right]^q}
            = \frac{(\beta_2+\teps)^q}{\teps^q} \frac{(p \teps)^q + \left[ \ln (1+\beta_2) \right]^q}{p^q (\beta_2+\teps)^q + \left[ \ln (1+\beta_2) \right]^q}
$$
is a bivariate function. 
Otherwise, we have
\begin{equation} \label{eq:nw_ZL1}  
 \omega^\ZL_0 = \frac{d_0}{d_0 + {d_1 \over \tbeta^\ZL_1} + {d_2 \over \tbeta^\ZL_2}},~
 \omega^\ZL_1 = \frac{d_1}{d_0 \tbeta^\ZL_1 + d_1 +  d_2 {\tbeta^\ZL_1 \over \tbeta^\ZL_2}},~
 \omega^\ZL_2 = \frac{d_2}{d_0 \tbeta^\ZL_2 + d_1 {\tbeta^\ZL_2 \over \tbeta^\ZL_1} + d_2}, 
\end{equation} 
where 
\begin{align*}
 \tbeta_1 &= \tbeta_1(p,q)= \frac{1 + \left[ \frac{\Delta}{p (\beta_0+\teps)} \right]^q}{1 + \left[ \frac{\Delta}{p (\beta_1+\teps)} \right]^q}
           = \frac{(\beta_1+\teps)^q}{(\beta_0+\teps)^q} \frac{p^q (\beta_0+\teps)^q + \Delta^q}{p^q (\beta_1+\teps)^q + \Delta^q}, \\
 \tbeta_2 &= \tbeta_2(p,q)= \frac{1 + \left[ \frac{\Delta}{p (\beta_0+\teps)} \right]^q}{1 + \left[ \frac{\Delta}{p (\beta_2+\teps)} \right]^q}
           = \frac{(\beta_2+\teps)^q}{(\beta_0+\teps)^q} \frac{p^q (\beta_0+\teps)^q + \Delta^q}{p^q (\beta_2+\teps)^q + \Delta^q}, \\
 \frac{\tbeta_1}{\tbeta_2} &= \frac{1 + \left[ \frac{\Delta}{p(\beta_2+\teps)} \right]^q}{1 + \left[ \frac{\Delta}{p(\beta_1+\teps)} \right]^q}                 
\end{align*}
are bivariate functions with $\Delta = | \ln (1+\beta_0) - \ln (1+\beta_2) |$.
The function 
$$ B(p,q) = \frac{1 + \left[ \frac{\Delta}{p (\beta_r+\teps)} \right]^q}{1 + \left[ \frac{\Delta}{p (\beta_s+\teps)} \right]^q} $$
has the following properties:
\begin{itemize}
\item $B$ is increasing in $p$ if $\beta_r > \beta_s$ while $B$ is decreasing in $p$ if $\beta_r < \beta_s$;
\item $B$ is decreasing in $q$ if $\beta_r > \beta_s$ whereas $B$ is increasing in $p$ if $\beta_r < \beta_s$. 
\end{itemize}

We omit the stencils where the analysis of nonlinear weights is similar to their counterparts in Section \ref{sec:FV_WENO_discontinuity} and focus on the stencils, of which the nonlinear weights affect the numerical solutions.

\setcounter{stencil}{1}
\begin{stencil}
$\{ I_{-4}, I_{-3}, I_{-2}, I_{-1}, I_{0} \}$ \\
With $\beta_0 = \beta_1 = 0$ and $\beta_2 = \frac{4}{3} \delta^2$, we have $\tau^\ZL_5 = \frac{1}{p} \ln (1+\beta_2)$.
Applying similar analysis as in Stencil \ref{case:2}, Section \ref{sec:Euler} gives
$$
   \omega^\ZL_s > d_s,~s = 0,1,~~0 \lessapprox \omega^\ZL_2 \ll d_2.
$$
If we fix $q$, then by \eqref{eq:nw_ZL0}, $\tbeta_2(p_1,q) > \tbeta_2(p_2,q)$ for $p_1 < p_2$, and hence
\begin{gather*}
 \omega^\ZL_s(p_1) > \omega^\ZL_s(p_2),~s = 0,1,~~ \omega^\ZL_2(p_1) < \omega^\ZL_2(p_2).   
\end{gather*}
When $p$ is fixed, we have $\tbeta_2(p,q_1) < \tbeta_2(p,q_2)$ for $q_1 < q_2$.  
Using \eqref{eq:nw_ZL0}, we find
\begin{gather*}
 \omega^\ZL_s(q_1) < \omega^\ZL_s(q_2),~s = 0,1,~~ \omega^\ZL_2(q_1) > \omega^\ZL_2(q_2).  
\end{gather*}
This analysis corresponds to $x=-0.01$ in Table \ref{tab:weight_comparison_first_ZL}. 
\end{stencil}

\begin{stencil}
$\{ I_{-3}, I_{-2}, I_{-1}, I_{0}, I_{1} \}$ \\
The smoothness indicators $\beta_0 = 0$, $\beta_1 = \frac{4}{3} \delta^2$ and $\beta_2 = \frac{10}{3} \delta^2$ yield $\tau^\ZL_5 = \frac{1}{p} \ln (1+\beta_2)$.
Following the approach of Stencil \ref{case:3}, Section \ref{sec:Euler}, we obtain
$$
   d_0 < \omega^\ZL_0 \lessapprox 1,~~0 \lessapprox \omega^\ZL_s \ll d_s,~s = 1,2.
$$
Fixing $q$, we have $\tbeta_s(p_1,q) > \tbeta_s(p_2,q),~s = 1,2$ for $p_1 < p_2$, and 
\begin{gather*}
 \omega^\ZL_0(p_1) > \omega^\ZL_0(p_2),~~ \omega^\ZL_s(p_1) < \omega^\ZL_s(p_2),~s = 1,2.   
\end{gather*}
Assuming $p$ is fixed, then $\tbeta_s(p,q_1) < \tbeta_s(p,q_2),~s = 1,2$ for $q_1 < q_2$, and
\begin{gather*}
 \omega^\ZL_0(q_1) < \omega^\ZL_0(q_2),~~ \omega^\ZL_s(q_1) > \omega^\ZL_s(q_2),~s = 1,2.  
\end{gather*}
The analysis is associated with $x=0$ in Table \ref{tab:weight_comparison_first_ZL}.
\end{stencil}

\begin{stencil}
$\{ I_{-2}, I_{-1}, I_{0}, I_{1}, I_{2} \}$ \\
Using the analysis analogous to what is done for Stencil \ref{case:3} above, we have the following results, at $x=0.01$ in Table \ref{tab:weight_comparison_first_ZL},
\begin{gather*}
 0 \lessapprox \omega^\ZL_s \ll d_s,~s = 0,1,~~ d_2 < \omega^\ZL_2 \lessapprox 1, \\
 \omega^\ZL_s(p_1) < \omega^\ZL_s(p_2),~ s = 0,1,~ p_2 > p_1 > 0, \\
 \omega^\ZL_s(q_1) > \omega^\ZL_s(q_2),~ s = 0,1,~ q_2 > q_1 \geqslant 1.
\end{gather*}
\end{stencil}

\begin{stencil}
$\{ I_{-1}, I_{0}, I_{1}, I_{2}, I_{3} \}$ \\
The results corresponding to $x=0.02$ in Table \ref{tab:weight_comparison_first_ZL}, which can be obtained in the same manner as in Stencil \ref{case:2} above, are as follows:
\begin{gather*}
 0 \lessapprox \omega^\ZL_0 \ll d_0,~~\omega^\ZL_s > d_s,~s = 1,2, \\
 \omega^\ZL_0(p_1) < \omega^\ZL_0(p_2),~ p_2 > p_1 > 0, \\
 \omega^\ZL_0(q_1) > \omega^\ZL_0(q_2),~ q_2 > q_1 \geqslant 1.
\end{gather*}
\end{stencil}

It seems that the numerical solutions in Table \ref{tab:numerical_solution_first_ZL} after the first stage of the TVD Runge-Kutta time step look the same due to the round-off errors in double precision. 

\begin{sidewaystable}
\renewcommand{\arraystretch}{1.2}
\centering

\vspace{4.6in}
\begin{tabular}{ccccccccc} 
\hline
$x$ & -0.03 & -0.02 & -0.01 & 0 & 0.01 & 0.02 & 0.03 & 0.04 \\ 
\hline
$\omega^\ZL_0$ ($p=1,q=1$) & 0.1 & 0.1 & 0.142857 & 1 & 3.273e-41 & 2.145e-41 & 0.1 & 0.1 \\ 
$\omega^\ZL_0$ ($p=2,q=1$) & 0.1 & 0.1 & 0.142857 & 1 & 5.546e-41 & 3.456e-41 & 0.1 & 0.1 \\
$\omega^\ZL_0$ ($p=1,q=2$) & 0.1 & 0.1 & 0.142857 & 1 & 1.850e-81 & 2.173e-81 & 0.1 & 0.1 \\
$\omega^\ZL_0$ ($p=2,q=2$) & 0.1 & 0.1 & 0.142857 & 1 & 6.501e-81 & 6.816e-81 & 0.1 & 0.1 \\
\hline
$\omega^\ZL_1$ ($p=1,q=1$) & 0.6 & 0.6 & 0.857143 & 8.592e-40 & 2.864e-40 & 0.666667 & 0.6 & 0.6 \\
$\omega^\ZL_1$ ($p=2,q=1$) & 0.6 & 0.6 & 0.857143 & 1.268e-39 & 4.228e-40 & 0.666667 & 0.6 & 0.6 \\
$\omega^\ZL_1$ ($p=1,q=2$) & 0.6 & 0.6 & 0.857143 & 6.166e-80 & 2.055e-80 & 0.666667 & 0.6 & 0.6 \\
$\omega^\ZL_1$ ($p=2,q=2$) & 0.6 & 0.6 & 0.857143 & 1.454e-79 & 4.846e-80 & 0.666667 & 0.6 & 0.6 \\
\hline
$\omega^\ZL_2$ ($p=1,q=1$) & 0.3 & 0.3 & 8.272e-41 & 2.946e-40 & 1 & 0.333333 & 0.3 & 0.3 \\  
$\omega^\ZL_2$ ($p=2,q=1$) & 0.3 & 0.3 & 1.333e-40 & 4.992e-40 & 1 & 0.333333 & 0.3 & 0.3 \\
$\omega^\ZL_2$ ($p=1,q=2$) & 0.3 & 0.3 & 8.380e-81 & 1.665e-80 & 1 & 0.333333 & 0.3 & 0.3 \\
$\omega^\ZL_2$ ($p=2,q=2$) & 0.3 & 0.3 & 2.629e-80 & 5.851e-80 & 1 & 0.333333 & 0.3 & 0.3 \\    
\hline
\end{tabular}
\caption{The nonlinear weights $\omega^\ZL_s$ near the discontinuity at $x=0$, with the linear weights $(d_0, d_1, d_2) = (0.1, 0.6, 0.3)$, at the first stage \eqref{eq:RK1st} of the single TVD Runge-Kutta time step for the Riemann problem \eqref{eq:weight_comparison}.}
\label{tab:weight_comparison_first_ZL}

\vspace{0.2in}

\begin{tabular}{ccccccccc}
\hline
$x$ & -0.03 & -0.02 & -0.01 & 0 & 0.01 & 0.02 & 0.03 & 0.04 \\ 
\hline
$\hat{f}^\ZL_j$ ($p=1,q=1$) & 1 & 1 & 1 & 1-2.220e-16 & -7.501e-41 & 7.149e-42 & 0 & 0 \\  
$\hat{f}^\ZL_j$ ($p=2,q=1$) & 1 & 1 & 1 & 1-2.220e-16 & -1.167e-40 & 1.152e-41 & 0 & 0 \\
$\hat{f}^\ZL_j$ ($p=1,q=2$) & 1 & 1 & 1 & 1-2.220e-16 & -4.967e-81 & 7.242e-82 & 0 & 0 \\ 
$\hat{f}^\ZL_j$ ($p=2,q=2$) & 1 & 1 & 1 & 1-2.220e-16 & -1.349e-80 & 2.272e-81 & 0 & 0 \\
\hline
\end{tabular}
\caption{The numerical fluxes $\hat{f}^\ZL_j$ near the discontinuity at $x=0$ at the first stage \eqref{eq:RK1st} of the single TVD Runge-Kutta time step for the Riemann problem \eqref{eq:weight_comparison}.}
\label{tab:numerical_flux_first_ZL}

\vspace{0.2in}

\begin{tabular}{cccccccc} 
\hline
$x$ & -0.025 & -0.015 & -0.005 & 0.005 & 0.015 & 0.025 & 0.035  \\ 
\hline
$\ubar^\ZL_j$ ($p=1,q=1$) & 1 & 1 & 1 & 0.5-1.110e-16 & -4.108e-41 & 3.574e-42 & 0 \\  
$\ubar^\ZL_j$ ($p=2,q=1$) & 1 & 1 & 1 & 0.5-1.110e-16 & -6.410e-41 & 5.760e-42 & 0 \\
$\ubar^\ZL_j$ ($p=1,q=2$) & 1 & 1 & 1 & 0.5-1.110e-16 & -2.846e-81 & 3.621e-82 & 0 \\
$\ubar^\ZL_j$ ($p=2,q=2$) & 1 & 1 & 1 & 0.5-1.110e-16 & -7.883e-81 & 1.136e-81 & 0 \\
$\ubar_j$                 & 1 & 1 & 1 & 0.5           & 0          & 0         & 0 \\    
\hline
\end{tabular}
\caption{The numerical solutions $\ubar^\ZL_j$, compared with the exact solution $\ubar_j$, near the discontinuity at $x=0$ at the first stage \eqref{eq:RK1st} of the single TVD Runge-Kutta time step for the Riemann problem \eqref{eq:weight_comparison}.}
\label{tab:numerical_solution_first_ZL}

\end{sidewaystable}

\setcounter{stencil}{9}
\begin{stencil}
$\{ I_{-2}, I_{-1}, I_{0}, I_{1}, I_{2} \}$ \\
We can estimate the smoothness indicators,
\begin{align*}
 \beta^\ZL_0 &\approx \frac{13}{12} (1-\nu)^2 \delta^2  + \frac{9}{4} (1-\nu)^2 \delta^2, \\
 \beta^\ZL_1 &\approx \frac{13}{12} (1-2\nu)^2 \delta^2 + \frac{1}{4} \delta^2, \\
 \beta^\ZL_2 &\approx \frac{13}{12} \nu^2 \delta^2      + \frac{1}{4} \nu^2 \delta^2.
\end{align*} 
Since $1-2\nu \leqslant 1-\nu$ and $\frac{1}{4} < \frac{9}{16} \leqslant \frac{9}{4} (1-\nu)^2 < \frac{9}{4}$, we obtain $\beta^\ZL_0 > \beta^\ZL_1, \beta^\ZL_2$ for $0 < \nu \leqslant 0.5$.
With $q$ fixed, we have $\tbeta_s(p_1,q) < \tbeta_s(p_2,q),~s = 1,2$ for $p_1 < p_2$, and hence
$$
   \omega^\ZL_0(p_1) < \omega^\ZL_0(p_2).
$$
When fixing $p$, we find $\tbeta_s(p,q_1) > \tbeta_s(p,q_2),~s = 1,2$ for $q_1 < q_2$, and hence
$$
   \omega^\ZL_0(q_1) > \omega^\ZL_0(q_2).  
$$
The graphs of $\omega^\ZL_s,~s = 0,1,2$ versus $\nu$ are shown in Figure \ref{fig:weightZL12}.
If we fix $\nu$ for $0 < \nu \leqslant 0.5$, the value of $\omega^\ZL_0$ with $p=2,\, q=1$ is larger than other combinations, which is consistent with the analysis.
\begin{figure}[htbp]
\centering
\includegraphics[width=0.325\textwidth]{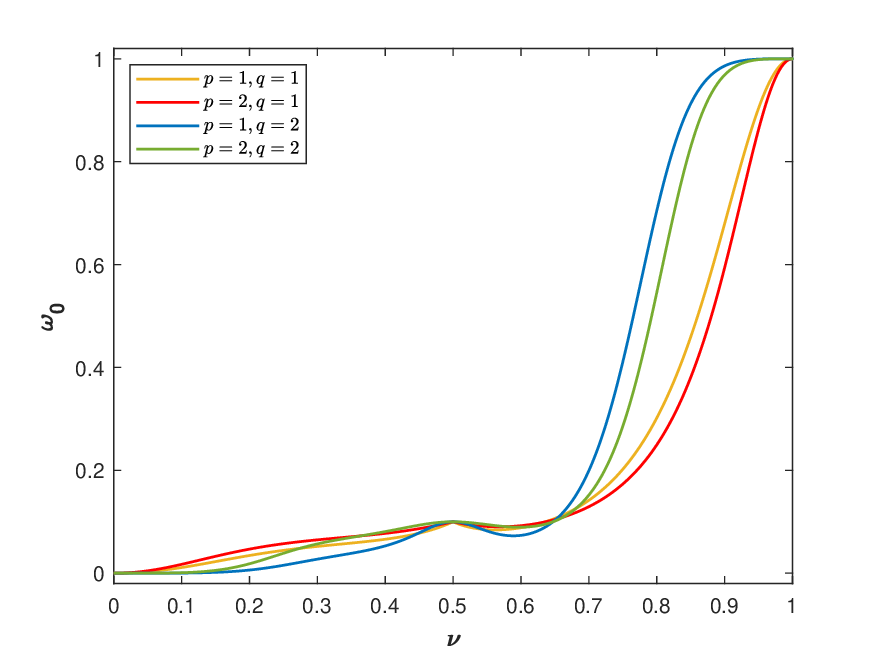}
\includegraphics[width=0.325\textwidth]{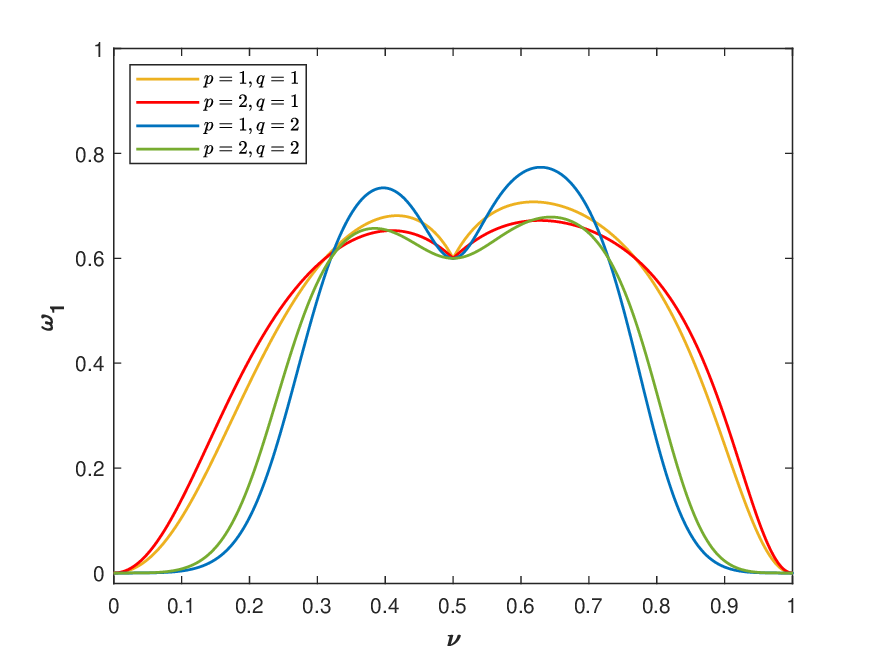}
\includegraphics[width=0.325\textwidth]{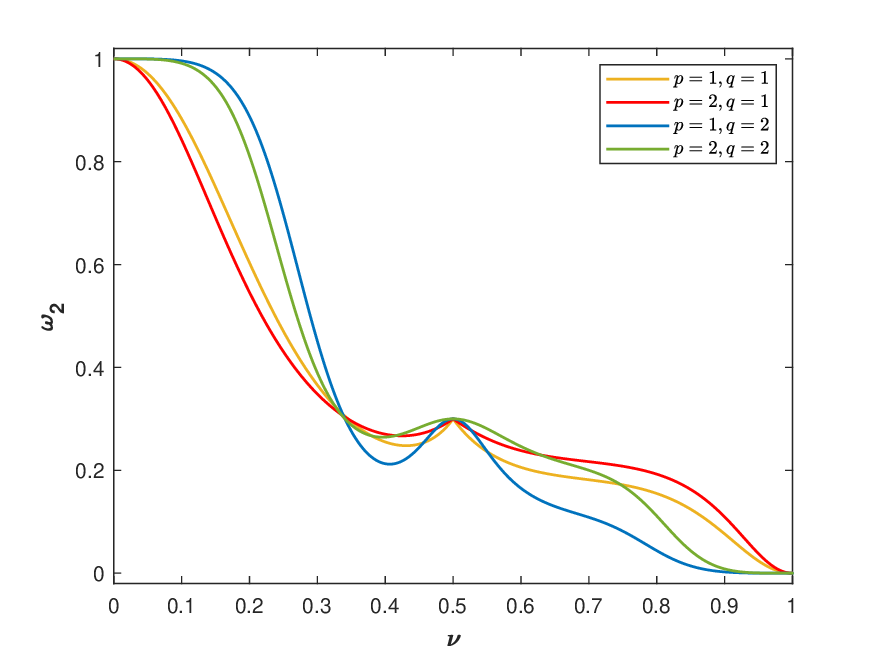}
\vspace{-0.3cm}
\caption{The nonlinear weights $\omega^\ZL_0$ (left), $\omega^\ZL_1$ (middle) and $\omega^\ZL_2$ (right) in Stencil \ref{case:10} with different $p$ and $q$.}
\label{fig:weightZL12}
\end{figure}

Recall from Section \ref{sec:RK2nd}, our goal is to increase the value of $\omega^\ZL_0$ in this stencil in order to reduce the errors $e^{\ZL,\,2}_0$ \eqref{eq:e2b0} and $e^{\ZL,\,2}_1$ \eqref{eq:e2b1}.
From the analysis above, increasing $p$ and decreasing $q$ can achieve our goal.
However, it is not that obvious in Table \ref{tab:numerical_solution_second_ZL} as it is at the early stage of approximations and we only show 6 digits after the integer part. 
\end{stencil}

\begin{sidewaystable}
\renewcommand{\arraystretch}{1.2}
\centering

\vspace{4.6in}
\begin{tabular}{cccccccccccc} 
\hline 
$x$ & -0.03 & -0.02 & -0.01 & 0 & 0.01 & 0.02 & 0.03 & 0.04 & 0.05 & 0.06 & 0.07 \\ 
\hline 
$\omega^\ZL_0$ ($p=1,q=1$) & 0.1 & 0.1 & 0.142857 & 1 & 0.1-9.714e-17 & 2.827e-40 & 7.196e-41 & 0.1 & 0.1 & 0.1 & 0.1 \\  
$\omega^\ZL_0$ ($p=2,q=1$) & 0.1 & 0.1 & 0.142857 & 1 & 0.1-4.163e-17 & 4.321e-40 & 1.106e-40 & 0.1 & 0.1 & 0.1 & 0.1 \\
$\omega^\ZL_0$ ($p=1,q=2$) & 0.1 & 0.1 & 0.142857 & 1 & 0.1           & 1.203e-79 & 2.343e-80 & 0.1 & 0.1 & 0.1 & 0.1 \\
$\omega^\ZL_0$ ($p=2,q=2$) & 0.1 & 0.1 & 0.142857 & 1 & 0.1           & 3.211e-79 & 6.370e-80 & 0.1 & 0.1 & 0.1 & 0.1 \\
\hline
$\omega^\ZL_1$ ($p=1,q=1$) & 0.6 & 0.6 & 0.857143 & 4.489e-39 & 0.6+3.331e-16 & 1.496e-39 & 0.666667 & 0.6 & 0.6 & 0.6 & 0.6 \\   
$\omega^\ZL_1$ ($p=2,q=1$) & 0.6 & 0.6 & 0.857143 & 7.178e-39 & 0.6+1.110e-16 & 2.393e-39 & 0.666667 & 0.6 & 0.6 & 0.6 & 0.6 \\
$\omega^\ZL_1$ ($p=1,q=2$) & 0.6 & 0.6 & 0.857143 & 1.745e-78 & 0.6           & 5.817e-79 & 0.666667 & 0.6 & 0.6 & 0.6 & 0.6 \\
$\omega^\ZL_1$ ($p=2,q=2$) & 0.6 & 0.6 & 0.857143 & 5.360e-78 & 0.6           & 1.787e-78 & 0.666667 & 0.6 & 0.6 & 0.6 & 0.6 \\
\hline
$\omega^\ZL_2$ ($p=1,q=1$) & 0.3 & 0.3 & 2.775e-40 & 2.544e-39 & 0.3-2.776e-16 & 1 & 0.333333 & 0.3 & 0.3 & 0.3 & 0.3 \\    
$\omega^\ZL_2$ ($p=2,q=1$) & 0.3 & 0.3 & 4.265e-40 & 3.889e-39 & 0.3-1.665e-16 & 1 & 0.333333 & 0.3 & 0.3 & 0.3 & 0.3 \\
$\omega^\ZL_2$ ($p=1,q=2$) & 0.3 & 0.3 & 9.036e-80 & 1.082e-78 & 0.3           & 1 & 0.333333 & 0.3 & 0.3 & 0.3 & 0.3 \\
$\omega^\ZL_2$ ($p=2,q=2$) & 0.3 & 0.3 & 2.457e-79 & 2.890e-78 & 0.3           & 1 & 0.333333 & 0.3 & 0.3 & 0.3 & 0.3 \\  
\hline
\end{tabular}
\caption{The nonlinear weights $\omega^\ZL_s$ near the discontinuity at $x = 0.005$, with the linear weights $(d_0, d_1, d_2) = (0.1, 0.6, 0.3)$, at the second stage \eqref{eq:RK2nd} of the single TVD Runge-Kutta time step for the Riemann problem \eqref{eq:weight_comparison}.}
\label{tab:weight_comparison_second_ZL}

\vspace{0.2in}

\begin{tabular}{cccccccccccc} 
\hline
$x$ & -0.03 & -0.02 & -0.01 & 0 & 0.01 & 0.02 & 0.03 & 0.04 & 0.05 & 0.06 & 0.07 \\ 
\hline
$\hat{f}^\ZL_j$ ($p=1,q=1$) & 1 & 1 & 1 & 1-2.220e-16 & 0.208333 & -2.061e-40 & 1.894e-41 & -2.144e-42 & 1.191e-43 & 0 & 0 \\    
$\hat{f}^\ZL_j$ ($p=2,q=1$) & 1 & 1 & 1 & 1-2.220e-16 & 0.208333 & -3.240e-40 & 2.939e-41 & -3.385e-42 & 1.920e-43 & 0 & 0 \\
$\hat{f}^\ZL_j$ ($p=1,q=2$) & 1 & 1 & 1 & 1-2.220e-16 & 0.208333 & -7.919e-80 & 4.462e-81 & -1.733e-82 & 1.207e-83 & 0 & 0 \\ 
$\hat{f}^\ZL_j$ ($p=2,q=2$) & 1 & 1 & 1 & 1-2.220e-16 & 0.208333 & -2.308e-79 & 1.225e-80 & -5.089e-82 & 3.787e-83 & 0 & 0 \\    
\hline
\end{tabular}
\caption{The numerical fluxes $\hat{f}^\N_j$ near the discontinuity at $x = 0.005$ at the second stage \eqref{eq:RK2nd} of the single TVD Runge-Kutta time step for the Riemann problem \eqref{eq:weight_comparison}.}
\label{tab:numerical_flux_second_ZL}

\vspace{0.2in}

\begin{tabular}{ccccccccccc} 
\hline
$x$ & -0.025 & -0.015 & -0.005 & 0.005 & 0.015 & 0.025 & 0.035 & 0.045 & 0.055 & 0.065 \\ 
\hline 
$\ubar^\ZL_j$ ($p=1,q=1$) & 1 & 1 & 1 & 0.223958 & 0.026042 & -2.723e-41 & 2.635e-42 & -2.829e-43 & 1.489e-44 & 0 \\   
$\ubar^\ZL_j$ ($p=2,q=1$) & 1 & 1 & 1 & 0.223958 & 0.026042 & -4.273e-41 & 4.097e-42 & -4.471e-43 & 2.400e-44 & 0 \\
$\ubar^\ZL_j$ ($p=1,q=2$) & 1 & 1 & 1 & 0.223958 & 0.026042 & -1.037e-80 & 5.794e-82 & -2.317e-83 & 1.509e-84 & 0 \\
$\ubar^\ZL_j$ ($p=2,q=2$) & 1 & 1 & 1 & 0.223958 & 0.026042 & -3.010e-80 & 1.595e-81 & -6.834e-83 & 4.733e-84 & 0 \\ 
$\ubar_j$                 & 1 & 1 & 1 & 0.5      & 0        & 0          & 0         & 0          & 0          & 0 \\     
\hline
\end{tabular}
\caption{The numerical solutions $\ubar^\ZL_j$, compared with the exact solution $\ubar_j$, near the discontinuity at $x = 0.005$ at the second stage \eqref{eq:RK2nd} of the single TVD Runge-Kutta time step for the Riemann problem \eqref{eq:weight_comparison}.}
\label{tab:numerical_solution_second_ZL}

\end{sidewaystable}

\setcounter{stencil}{18}
\begin{stencil}
$\{ I_{-2}, I_{-1}, I_0, I_1, I_2 \}$ \\
Given the estimates of the smoothness indicators,
\begin{align*}
 \beta^\ZL_0 &= \frac{13}{12} \left[ (1-\nu) \delta - e^\ZL_0 \right]^2 +
                \frac{9}{4}   \left( \delta - \nu \delta - e^\ZL_0 \right)^2, \\
 \beta^\ZL_1 &= \frac{13}{12} \left[ (1-\nu) \delta - e^\ZL_0 - (\nu \delta + e^\ZL_0 - e^\ZL_1) \right]^2 + 
                \frac{1}{4}   \left( \delta - e^\ZL_{2} \right)^2, \\
 \beta^\ZL_2 &\approx \frac{13}{12} \left[ (\nu \delta + e^\ZL_0 - e^\ZL_1) - e^\ZL_1 \right]^2 + 
                      \frac{1}{4}   \left[ 3(\nu \delta + e^\ZL_0 - e^\ZL_1) - e^\ZL_1 \right]^2,
\end{align*}
we have, by \eqref{eq:e2b0}, \eqref{eq:e2b1} and \eqref{eq:e2b01}, $\beta^\ZL_0 > \beta^\ZL_1 > \beta^\ZL_2$ for $0 < \nu \leqslant 0.5$.
When fixing $q$, we have $\tbeta_s(p_1,q) < \tbeta_s(p_2,q),~s = 1,2$ for $p_1 < p_2$, and hence
$$
   \omega^\ZL_0(p_1) < \omega^\ZL_0(p_2),~~ \omega^\ZL_2(p_1) > \omega^\ZL_2(p_2).
$$
If we fix $p$, then $\tbeta_s(p,q_1) > \tbeta_s(p,q_2),~s = 1,2$ for $q_1 < q_2$, and hence
$$
   \omega^\ZL_0(q_1) > \omega^\ZL_0(q_2),~~ \omega^\ZL_2(q_1) < \omega^\ZL_2(q_2).  
$$
Figure \ref{fig:weightZL25} shows the graphs of $\omega^\ZL_s,~s = 0,1,2$ versus $\nu$.
For any fixed $\nu$, the value of $\omega^\ZL_0$ with $p=2,\, q=1$ is larger than other combinations, matching the analysis above.
We also note from Figure \ref{fig:weightZL25} that $\omega^\ZL_1(p_1) < \omega^\ZL_1(p_2)$ for $p_1 < p_2$ and $\omega^\ZL_1(q_1) > \omega^\ZL_1(q_2)$ for $q_1 < q_2$.
\begin{figure}[htbp]
\centering
\includegraphics[width=0.325\textwidth]{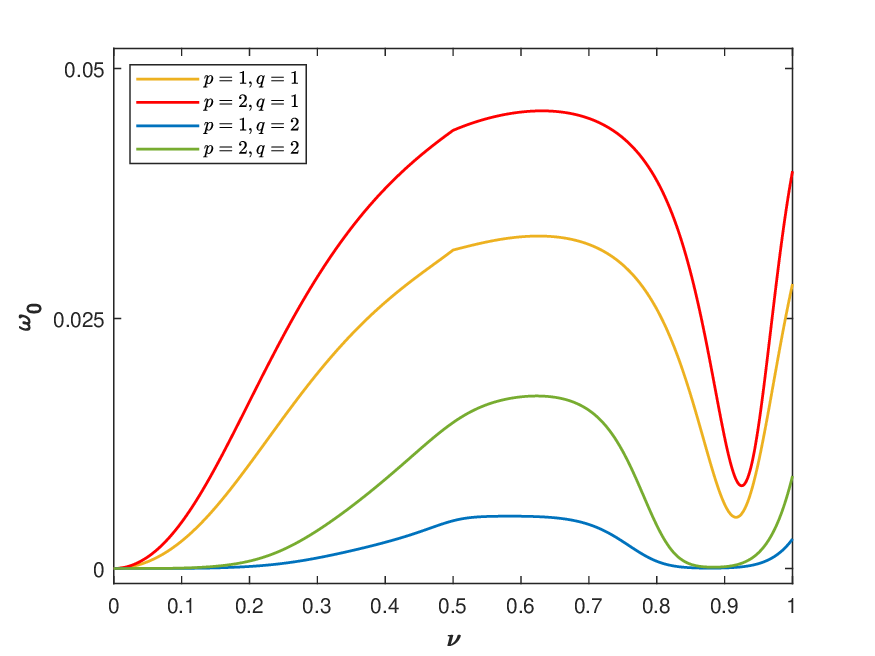}
\includegraphics[width=0.325\textwidth]{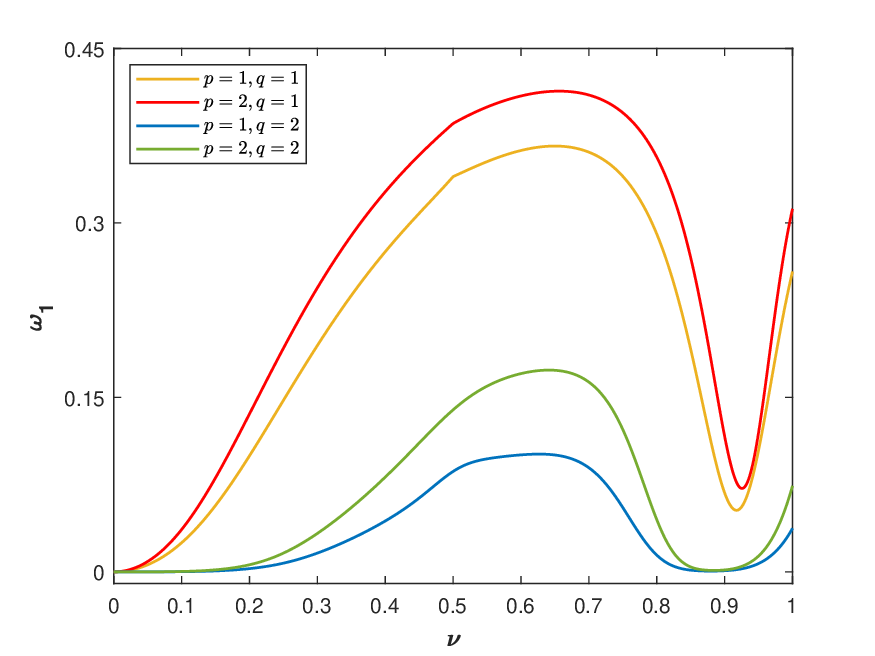}
\includegraphics[width=0.325\textwidth]{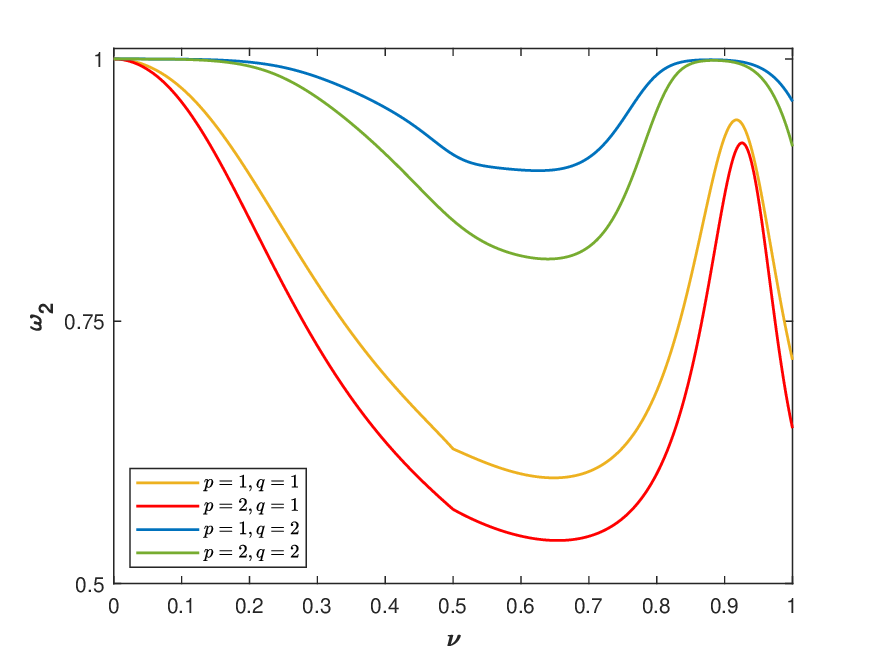}
\vspace{-0.3cm}
\caption{The nonlinear weights $\omega^\ZL_0$ (left), $\omega^\ZL_1$ (middle) and $\omega^\ZL_2$ (right) in Stencil \ref{case:19} with different $p$ and $q$.}
\label{fig:weightZL25}
\end{figure}

In Section \ref{sec:RK3rd}, increasing the value of $\omega^\ZL_0$ in this stencil would reduce the errors $\es^{\ZL,\,1}_0$ \eqref{eq:e1b0} and $\es^{\ZL,\,1}_1$ \eqref{eq:e1b1}.
So the increase of $p$ and the decrease of $q$ leads to the better numerical solutions after the single TVD Runge-Kutta time step and thus the reduction of numerical dissipation around the jump discontinuity, which is verified in Table \ref{tab:numerical_solution_third_ZL}.
\end{stencil}

\begin{sidewaystable}
\renewcommand{\arraystretch}{1.2}
\centering
\vspace{4.6in}
\begin{tabular}{cccccccc} 
\hline
$x$ & -0.03 & -0.02 & -0.01 & 0 & 0.01 & 0.02 & 0.03 \\ 
\hline
$\omega^\ZL_0$ ($p=1,q=1$) & 0.1 & 0.1 & 0.142857 & 1 & 0.031848 & 3.940e-3 & 1.890e-39 \\  
$\omega^\ZL_0$ ($p=2,q=1$) & 0.1 & 0.1 & 0.142857 & 1 & 0.043815 & 5.949e-3 & 2.843e-39 \\
$\omega^\ZL_0$ ($p=1,q=2$) & 0.1 & 0.1 & 0.142857 & 1 & 0.004774 & 2.966e-5 & 5.357e-78 \\
$\omega^\ZL_0$ ($p=2,q=2$) & 0.1 & 0.1 & 0.142857 & 1 & 0.014596 & 8.156e-5 & 1.353e-77 \\   
\hline
$\omega^\ZL_1$ ($p=1,q=1$) & 0.6 & 0.6 & 0.857143 & 1.403e-39 & 0.339855 & 0.102831 & 2.269e-37 \\  
$\omega^\ZL_1$ ($p=2,q=1$) & 0.6 & 0.6 & 0.857143 & 2.059e-39 & 0.385533 & 0.113228 & 2.326e-37 \\
$\omega^\ZL_1$ ($p=1,q=2$) & 0.6 & 0.6 & 0.857143 & 1.647e-79 & 0.086346 & 0.005228 & 2.448e-74 \\
$\omega^\ZL_1$ ($p=2,q=2$) & 0.6 & 0.6 & 0.857143 & 3.798e-79 & 0.139327 & 0.005537 & 2.453e-74 \\  
\hline
$\omega^\ZL_2$ ($p=1,q=1$) & 0.3 & 0.3 & 1.261e-40 & 5.284e-40 & 0.628296 & 0.893229 & 1 \\  
$\omega^\ZL_2$ ($p=2,q=1$) & 0.3 & 0.3 & 1.988e-40 & 8.563e-40 & 0.570652 & 0.880822 & 1 \\ 
$\omega^\ZL_2$ ($p=1,q=2$) & 0.3 & 0.3 & 1.898e-80 & 4.924e-80 & 0.908879 & 0.994742 & 1 \\
$\omega^\ZL_2$ ($p=2,q=2$) & 0.3 & 0.3 & 5.599e-80 & 1.568e-79 & 0.846077 & 0.994382 & 1 \\
\hline
$x$ & 0.04 & 0.05 & 0.06 & 0.07 & 0.08 & 0.09 & 0.1 \\ 
\hline
$\omega^\ZL_0$ ($p=1,q=1$) & 2.458e-38 & 0.1 & 0.1 & 0.1 & 0.1 & 0.1 & 0.1 \\  
$\omega^\ZL_0$ ($p=2,q=1$) & 3.688e-38 & 0.1 & 0.1 & 0.1 & 0.1 & 0.1 & 0.1 \\
$\omega^\ZL_0$ ($p=1,q=2$) & 2.719e-75 & 0.1 & 0.1 & 0.1 & 0.1 & 0.1 & 0.1 \\
$\omega^\ZL_0$ ($p=2,q=2$) & 6.800e-75 & 0.1 & 0.1 & 0.1 & 0.1 & 0.1 & 0.1 \\ 
\hline
$\omega^\ZL_1$ ($p=1,q=1$) & 0.666667 & 0.6 & 0.6 & 0.6 & 0.6 & 0.6 & 0.6 \\  
$\omega^\ZL_1$ ($p=2,q=1$) & 0.666667 & 0.6 & 0.6 & 0.6 & 0.6 & 0.6 & 0.6 \\
$\omega^\ZL_1$ ($p=1,q=2$) & 0.666667 & 0.6 & 0.6 & 0.6 & 0.6 & 0.6 & 0.6 \\
$\omega^\ZL_1$ ($p=2,q=2$) & 0.666667 & 0.6 & 0.6 & 0.6 & 0.6 & 0.6 & 0.6 \\  
\hline
$\omega^\ZL_2$ ($p=1,q=1$) & 0.333333 & 0.3 & 0.3 & 0.3 & 0.3 & 0.3 & 0.3 \\  
$\omega^\ZL_2$ ($p=2,q=1$) & 0.333333 & 0.3 & 0.3 & 0.3 & 0.3 & 0.3 & 0.3 \\ 
$\omega^\ZL_2$ ($p=1,q=2$) & 0.333333 & 0.3 & 0.3 & 0.3 & 0.3 & 0.3 & 0.3 \\
$\omega^\ZL_2$ ($p=2,q=2$) & 0.333333 & 0.3 & 0.3 & 0.3 & 0.3 & 0.3 & 0.3 \\
\hline
\end{tabular}
\caption{The nonlinear weights $\omega^\ZL_s$ near the discontinuity at $x=0.005$, with the linear weights $(d_0, d_1, d_2) = (0.1, 0.6, 0.3)$, at the third stage \eqref{eq:RK3rd} of the single TVD Runge-Kutta time step for the Riemann problem \eqref{eq:weight_comparison}.}
\label{tab:weight_comparison_third_ZL}
\end{sidewaystable}

\begin{sidewaystable}
\renewcommand{\arraystretch}{1.2}
\centering

\vspace{5.0in}
\begin{tabular}{cccccccc} 
\hline
$x$ & -0.03 & -0.02 & -0.01 & 0 & 0.01 & 0.02 & 0.03 \\ 
\hline
$\hat{f}^\ZL_j$ ($p=1,q=1$) & 1 & 1 & 1 & 1-2.220e-16 & 0.056811 & 0.006619 & -9.080e-40 \\  
$\hat{f}^\ZL_j$ ($p=2,q=1$) & 1 & 1 & 1 & 1-2.220e-16 & 0.047506 & 0.006589 & -8.945e-40 \\
$\hat{f}^\ZL_j$ ($p=1,q=2$) & 1 & 1 & 1 & 1-2.220e-16 & 0.088029 & 0.008557 & -1.060e-76 \\ 
$\hat{f}^\ZL_j$ ($p=2,q=2$) & 1 & 1 & 1 & 1-2.220e-16 & 0.079344 & 0.008555 & -1.059e-76 \\ 
\hline
$x$ & 0.04 & 0.05 & 0.06 & 0.07 & 0.08 & 0.09 & 0.1 \\ 
\hline
$\hat{f}^\ZL_j$ ($p=1,q=1$) & 2.180e-40 & -1.694e-42 & 1.608e-43 & -1.266e-44 & 4.965e-46 & 0 & 0 \\  
$\hat{f}^\ZL_j$ ($p=2,q=1$) & 3.274e-40 & -2.651e-42 & 2.522e-43 & -2.010e-44 & 8.000e-46 & 0 & 0 \\
$\hat{f}^\ZL_j$ ($p=1,q=2$) & 2.361e-77 & -4.885e-82 & 2.552e-83 & -1.099e-84 & 5.029e-86 & 0 & 0 \\
$\hat{f}^\ZL_j$ ($p=2,q=2$) & 5.903e-77 & -1.400e-81 & 7.168e-83 & -3.304e-84 & 1.578e-85 & 0 & 0 \\  
\hline
\end{tabular}
\caption{The numerical fluxes $\hat{f}^\ZL_j$ near the discontinuity at $x=0.005$ at the third stage \eqref{eq:RK3rd} of the single TVD Runge-Kutta time step for the Riemann problem \eqref{eq:weight_comparison}.}
\label{tab:numerical_flux_third_ZL}

\vspace{0.2in}

\begin{tabular}{cccccccc} 
\hline
$x$ & -0.025 & -0.015 & -0.005 & 0.005 & 0.015 & 0.025 & 0.035 \\ 
\hline
$\ubar^\ZL_j$ ($p=1,q=1$) & 1 & 1 & 1 & 0.463702 & 0.034092 & 0.002206 & -3.736e-40 \\  
$\ubar^\ZL_j$ ($p=2,q=1$) & 1 & 1 & 1 & 0.466803 & 0.031000 & 0.002196 & -4.046e-40 \\
$\ubar^\ZL_j$ ($p=1,q=2$) & 1 & 1 & 1 & 0.453296 & 0.043852 & 0.002852 & -4.320e-77 \\
$\ubar^\ZL_j$ ($p=2,q=2$) & 1 & 1 & 1 & 0.456191 & 0.040957 & 0.002852 & -5.496e-77 \\
$\ubar_j$                 & 1 & 1 & 1 & 0.5      & 0        & 0        & 0 \\ 
\hline
$x$ & 0.045 & 0.055 & 0.065 & 0.075 & 0.085 & 0.095 \\
\hline
$\ubar^\ZL_j$ ($p=1,q=1$) & 7.305e-41 & -6.082e-43 & 5.782e-44 & -4.384e-45 & 1.655e-46 & 0 \\  
$\ubar^\ZL_j$ ($p=2,q=1$) & 1.097e-40 & -9.519e-43 & 9.078e-44 & -6.968e-45 & 2.667e-46 & 0 \\
$\ubar^\ZL_j$ ($p=1,q=2$) & 7.869e-78 & -1.703e-82 & 8.872e-84 & -3.832e-85 & 1.676e-86 & 0 \\
$\ubar^\ZL_j$ ($p=2,q=2$) & 1.968e-77 & -4.875e-82 & 2.499e-83 & -1.154e-84 & 5.259e-86 & 0 \\
$\ubar_j$                 & 0         & 0          & 0         & 0          & 0         & 0 \\ 
\hline
\end{tabular}
\caption{The numerical solutions $\ubar^\ZL_j$, compared with the exact solution $\ubar_j$, near the discontinuity at $x=0.005$ at the third stage \eqref{eq:RK3rd} of the single TVD Runge-Kutta time step for the Riemann problem \eqref{eq:weight_comparison}.}
\label{tab:numerical_solution_third_ZL}

\vspace{0.2in}

\begin{tabular}{ccccccccc} 
\hline
$x$ & 0.9655 & 0.975 & 0.985 & 0.995 & 1.005 & 1.015 & 1.025 & 1.035 \\ 
\hline
$\ubar^\ZL_j$ ($p=1,q=1$) & 0.988142 & 0.947999 & 0.830315 & 0.625163 & 0.380986 & 0.170123 & 0.046944 & 0.009850 \\  
$\ubar^\ZL_j$ ($p=2,q=1$) & 0.990070 & 0.951865 & 0.834157 & 0.627000 & 0.380018 & 0.166594 & 0.042425 & 0.007464 \\  
$\ubar^\ZL_j$ ($p=1,q=2$) & 0.970122 & 0.926589 & 0.816626 & 0.622006 & 0.385993 & 0.183243 & 0.068118 & 0.028137 \\
$\ubar^\ZL_j$ ($p=2,q=2$) & 0.974881 & 0.934869 & 0.824114 & 0.625452 & 0.384047 & 0.176093 & 0.058424 & 0.022929 \\
$\ubar_j$                 & 1        & 1        & 1        & 1        & 0        & 0        & 0        & 0 \\ 
\hline
\end{tabular}
\caption{The numerical solutions $\ubar^\ZL_j$ at $T=1$, compared with the exact solution $\ubar_j$, near the discontinuity at $x = 1$ for the Riemann problem \eqref{eq:weight_comparison}.}
\label{tab:numerical_solution_final_ZL}

\end{sidewaystable}

Tables \ref{tab:numerical_solution_third_ZL} and \ref{tab:numerical_solution_final_ZL} show that the WENO scheme that yields more accurate numerical solutions at the first TVD Runge-Kutta time step, will perform better at the final time, as is observed in Section \ref{sec:FV_WENO_discontinuity}.

\section{Finite volume WENO schemes for 2D scalar conservation laws} \label{sec:FV_WENO_2d}
Now we extend the Cartesian grid to two dimensions.
Consider the domain $[a,b] \times [c,d]$ covered by cells
$$
   I_{ij} = [x_{i-1/2}, x_{i+1/2}] \times [y_{j-1/2}, y_{j+1/2}],~ 1 \leqslant i \leqslant N_x,~ 1 \leqslant j \leqslant N_y,
$$
where
$$ 
   a = x_{1/2} < x_{3/2} < \cdots < x_{N_x-1/2} < x_{N_x+1/2} = b,
$$
and
$$ 
   c = y_{1/2} < y_{3/2} < \cdots < y_{N_y-1/2} < y_{N_y+1/2} = d.
$$
For the cell $I_{ij}$, the cell center $(x_i, y_j)$ is
$$
   x_i = \frac{1}{2} (x_{i-1/2}+ x_{i+1/2}), \quad y_i = \frac{1}{2} (y_{j-1/2}+ y_{j+1/2}), 
$$ 
and we denote the grid sizes by 
$$
   \Delta x_i = x_{i+1/2} - x_{i-1/2}, \quad i = 1, 2, \cdots, N_x
$$
and
$$
   \Delta y_j = y_{j+1/2} - y_{j-1/2}, \quad j = 1, 2, \cdots, N_y.
$$
We use
$$ 
   \Delta x = \max_{1 \leqslant i \leqslant N_x} \Delta x_i, \quad \Delta y = \max_{1 \leqslant j \leqslant N_y} \Delta y_j
$$
to denote the maximum grid sizes.

We continue to describe the finite volume WENO schemes for two-dimensional scalar conservation laws
\begin{equation} \label{eq:2d_scalar_hyperbolic}
 u_t(x,y,t) + f_x(u(x,y,t)) + g_y(u(x,y,t)) = 0.
\end{equation}
Integrating \eqref{eq:2d_scalar_hyperbolic} over the cell $I_{ij}$, we obtain
\begin{equation} \label{eq:2d_scalar_hyperbolic_integral}
\begin{aligned}
 \frac{\dx \ubar(x_i,y_j,t)}{\dx t} = {} & - \frac{1}{\Delta x_i \Delta y_j} \Bigg[ \int_{y_{j-1/2}}^{y_{j+1/2}} f(u(x_{i+1/2},y,t)) - f(u(x_{i-1/2},y,t)) \dx y \\
                                      {} & + \int_{x_{i-1/2}}^{x_{i+1/2}} f(u(x,y_{j+1/2},t)) - f(u(x,y_{j-1/2},t)) \dx x \Bigg],                                              
\end{aligned}
\end{equation}
with the cell average
$$
   \ubar(x_i,y_j,t) = \frac{1}{\Delta x_i \Delta y_j} \int_{y_{j-1/2}}^{y_{j+1/2}} \int^{x_{i+1/2}}_{x_{i-1/2}} u(\xi,\eta,t) \, \dx \xi \, \dx \eta.
$$
Here we have four integrals to approximate.
By setting $\eta = \frac{2(y - y_j)}{\Delta y_j} = \frac{(y - y_j)}{\delta_j}$,  
\begin{align*}
 \frac{1}{\Delta y_j} \int_{y_{j-1/2}}^{y_{j+1/2}} f(u(x_{i+1/2},y,t)) \dx y &= \frac{1}{\Delta y_j} \int_{-1}^1 f \left( u \left( x_{i+1/2},y_j + \delta_j \eta,t \right) \right) \frac{\Delta y_j}{2} \dx \eta \\
 &= \frac{1}{2} \int_{-1}^1 f \left( u \left( x_{i+1/2},y_j + \delta_j \eta,t \right) \right) \dx \eta \\
 & \approx \frac{1}{2} \sum_{\gamma=1}^n c_{\gamma} f \left( u \left( x_{i+1/2},y_j + \delta_j \eta_{\gamma},t \right) \right) \\
 & \approx \frac{1}{2} \sum_{\gamma=1}^n c_{\gamma} h \left( u^-_{i+1/2,y_j + \delta_j \eta_{\gamma}}, u^+_{i+1/2,y_j + \delta_j \eta_{\gamma}} \right) ,
\end{align*}
where $c_{\gamma}$ and $\eta_{\gamma}$ are Gaussian quadrature nodes and weights, and $u^{\pm}_{i+1/2,y}$ are the fifth-order accurate approximation obtained by WENO reconstruction described in Appendix.
We use $\hat{f}_{i+1/2,j}$ to denote the last term, which is an approximation to the first integral in \eqref{eq:2d_scalar_hyperbolic_integral}.
The flux $\hat{g}_{i+1/2,j}$ can be obtained in a similar way.
So we approximate \eqref{eq:2d_scalar_hyperbolic_integral} by the following conservative scheme
$$
   \frac{\dx \ubar_{ij}(t)}{\dx t}  = - \frac{\hat{f}_{i+1/2,j} - \hat{f}_{i-1/2,j}}{\Delta x_i} - \frac{\hat{g}_{i,j+1/2} - \hat{g}_{i,j-1/2}}{\Delta y_j},
$$
where $\ubar_{ij}(t)$ is the numerical approximation to the cell average $\ubar(x_i,y_j,t)$.

\section{Numerical results} \label{sec:nr}
In this section, we present some numerical experiments to compare the proposed WENO scheme, referred as WENO-ZL, with the WENO-JS, WENO-M, WENO-Z and WENO-ZR $(p=2)$ schemes.
We use the one-dimensional linear advection equation to verify the order of accuracy of the WENO schemes in terms of $L_1, L_2$ and $L_{\infty}$ error norms:
\begin{align*}
 & L_1 = \frac{1}{N} \sum_{i=1}^N \left| u(x_i, T) - u_i(T) \right|, \\
 & L_2 = \sqrt{\frac{1}{N} \sum_{i=1}^N \left( u(x_i, T) - u_i(T) \right)^2}, \\
 & L_{\infty} = \max_{0 \leqslant i \leqslant N} \left| u(x_i, T)- u_i(T) \right|,
\end{align*}
with $u(x_i, T)$ the exact solution and $u_i(T)$ the numerical approximation at the final time $t = T$.
The rest of experiments exhibit the numerical solutions from WENO-ZL, in comparison with WENO-JS, WENO-M, WENO-Z and WENO-ZR.
We choose $\epsilon = 10^{-40}$ for the WENO-M, WENO-Z, WENO-ZR and WENO-ZL schemes whereas $\epsilon = 10^{-6}$ for WENO-JS as in \cite{Jiang}.
The CFL number is set to be 0.4, unless otherwise stated. 
We employ the explicit third-order TVD Runge-Kutta method for time integration.

\subsection{1D scalar case} 
\begin{example} \label{ex:advection_1d_order}
We first examine the order of accuracy for the one-dimensional linear advection equation
$$
   u_t + u_x = 0, \quad -1 \leqslant x \leqslant 1,
$$
with the initial condition $u(x,0) = \sin(\pi x)$ and the periodic boundary condition. 
The exact solution is given by
$$
   u(x,t)= \sin \left( \pi (x-t) \right).
$$
The numerical solution is computed up to the final time $T=8$ with the time step $\Delta t = 0.1 \cdot \Delta x$. 

The $L_1, L_2$ and $L_{\infty}$ errors versus $N$, as well as the order of accuracy, for the WENO-JS, WENO-M, WENO-Z, WENO-ZR and WENO-ZL ($p=q=2$) schemes are displayed in Tables \ref{tab:advection_1d_order_L1}, \ref{tab:advection_1d_order_L2} and \ref{tab:advection_1d_order_Linf}, respectively. 
The fourth-order accuracy is achieved for all schemes.
But we also notice that the CFL number would affect the order of accuracy.
\end{example}

\begin{table}[h!]
\renewcommand{\arraystretch}{1.1}
\scriptsize
\centering
\caption{$L_1$ error and order of accuracy for Example \ref{ex:advection_1d_order}.}      
\begin{tabular}{clcrlcrlcrlcrlc} 
\hline  
N & \multicolumn{2}{l}{WENO-JS} & & \multicolumn{2}{l}{WENO-M} & & \multicolumn{2}{l}{WENO-Z} & & \multicolumn{2}{l}{WENO-ZR} & & \multicolumn{2}{l}{WENO-ZL} \\ 
    \cline{2-3}                     \cline{5-6}                    \cline{8-9}                    \cline{11-12}                   \cline{14-15}
  & Error & Order               & & Error & Order              & & Error & Order              & & Error & Order               & & Error & Order \\
\hline
10  & 1.00e-1 & --     & & 3.18e-2 & --     & & 2.88e-2 & --     & & 2.36e-2 & --     & & 2.36e-2 & --     \\  
20  & 5.47e-3 & 4.1945 & & 8.47e-4 & 5.2282 & & 8.79e-4 & 5.0338 & & 8.17e-4 & 4.8520 & & 8.18e-4 & 4.8509 \\  
40  & 1.81e-4 & 4.9164 & & 2.80e-5 & 4.9201 & & 2.80e-5 & 4.9699 & & 2.79e-5 & 4.8705 & & 2.79e-5 & 4.8715 \\
80  & 5.90e-6 & 4.9379 & & 1.12e-6 & 4.6447 & & 1.12e-6 & 4.6470 & & 1.12e-6 & 4.6418 & & 1.12e-6 & 4.6418 \\ 
160 & 2.15e-7 & 4.7818 & & 6.53e-8 & 4.0994 & & 6.53e-8 & 4.0996 & & 6.53e-8 & 4.0993 & & 6.53e-8 & 4.0993 \\
\hline
\end{tabular}
\label{tab:advection_1d_order_L1}
\end{table}

\begin{table}[h!]
\renewcommand{\arraystretch}{1.1}
\scriptsize
\centering
\caption{$L_2$ error and order of accuracy for Example \ref{ex:advection_1d_order}.}      
\begin{tabular}{clcrlcrlcrlcrlc} 
\hline  
N & \multicolumn{2}{l}{WENO-JS} & & \multicolumn{2}{l}{WENO-M} & & \multicolumn{2}{l}{WENO-Z} & & \multicolumn{2}{l}{WENO-ZR} & & \multicolumn{2}{l}{WENO-ZL} \\ 
    \cline{2-3}                     \cline{5-6}                    \cline{8-9}                    \cline{11-12}                   \cline{14-15}
  & Error & Order               & & Error & Order              & & Error & Order              & & Error & Order               & & Error & Order \\
\hline
10  & 1.12e-1 & --     & & 3.63e-2 & --     & & 3.24e-2 & --     & & 2.67e-2 & --     & & 2.67e-2 & --     \\  
20  & 5.90e-3 & 4.2535 & & 9.42e-4 & 5.2682 & & 9.89e-4 & 5.0353 & & 9.12e-4 & 4.8706 & & 9.13e-4 & 4.8688 \\  
40  & 2.05e-4 & 4.8493 & & 3.11e-5 & 4.9201 & & 3.16e-5 & 4.9691 & & 3.10e-5 & 4.8762 & & 3.10e-5 & 4.8774 \\
80  & 6.64e-6 & 4.9457 & & 1.24e-6 & 4.6458 & & 1.24e-6 & 4.6648 & & 1.24e-6 & 4.6430 & & 1.24e-6 & 4.6431 \\ 
160 & 2.40e-7 & 4.7927 & & 7.25e-8 & 4.0995 & & 7.25e-8 & 4.1020 & & 7.25e-8 & 4.0993 & & 7.25e-8 & 4.0994 \\
\hline
\end{tabular}
\label{tab:advection_1d_order_L2}
\end{table}

\begin{table}[h!]
\renewcommand{\arraystretch}{1.1}
\scriptsize
\centering
\caption{$L_{\infty}$ error and order of accuracy for Example \ref{ex:advection_1d_order}.}      
\begin{tabular}{clcrlcrlcrlcrlc} 
\hline  
N & \multicolumn{2}{l}{WENO-JS} & & \multicolumn{2}{l}{WENO-M} & & \multicolumn{2}{l}{WENO-Z} & & \multicolumn{2}{l}{WENO-ZR} & & \multicolumn{2}{l}{WENO-ZL} \\ 
    \cline{2-3}                     \cline{5-6}                    \cline{8-9}                    \cline{11-12}                   \cline{14-15}
  & Error & Order               & & Error & Order              & & Error & Order              & & Error & Order               & & Error & Order \\
\hline
10  & 1.58e-1 & --     & & 4.96e-2 & --     & & 4.35e-2 & --     & & 3.63e-2 & --     & & 3.62e-2 & --     \\  
20  & 8.25e-3 & 4.2618 & & 1.30e-3 & 5.2509 & & 1.39e-3 & 4.9637 & & 1.29e-3 & 4.8165 & & 1.29e-3 & 4.8119 \\ 
40  & 3.18e-4 & 4.6988 & & 4.39e-5 & 4.8887 & & 4.48e-5 & 4.9601 & & 4.39e-5 & 4.8760 & & 4.39e-5 & 4.8771 \\
80  & 1.08e-5 & 4.8769 & & 1.76e-6 & 4.6441 & & 1.77e-6 & 4.6648 & & 1.76e-6 & 4.6431 & & 1.76e-6 & 4.6433 \\  
160 & 3.85e-7 & 4.8125 & & 1.03e-7 & 4.0994 & & 1.03e-7 & 4.1054 & & 1.03e-7 & 4.0993 & & 1.03e-7 & 4.0993 \\
\hline
\end{tabular}
\label{tab:advection_1d_order_Linf}
\end{table}

\begin{figure}[htbp]
\centering
\includegraphics[width=0.325\textwidth]{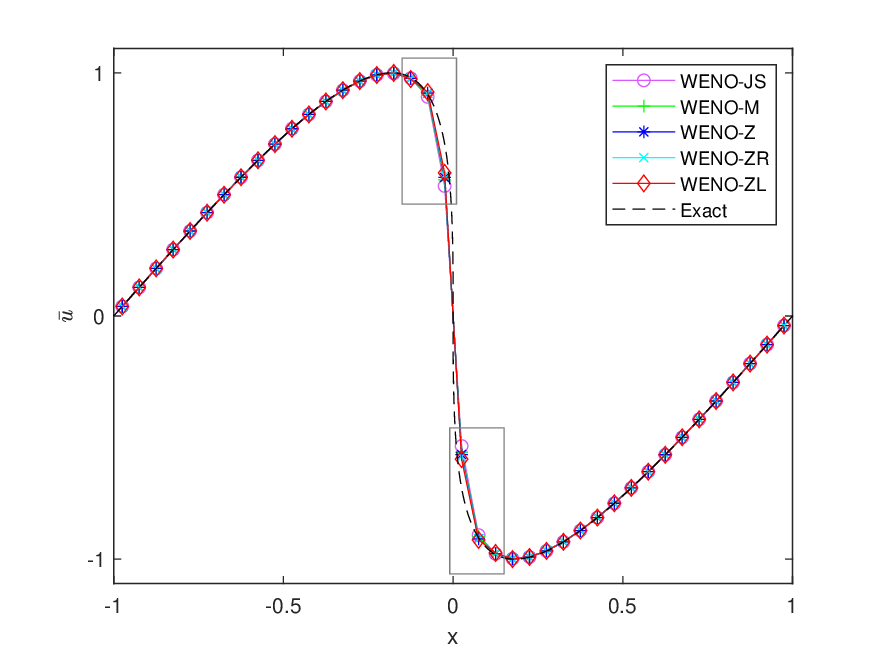}
\includegraphics[width=0.325\textwidth]{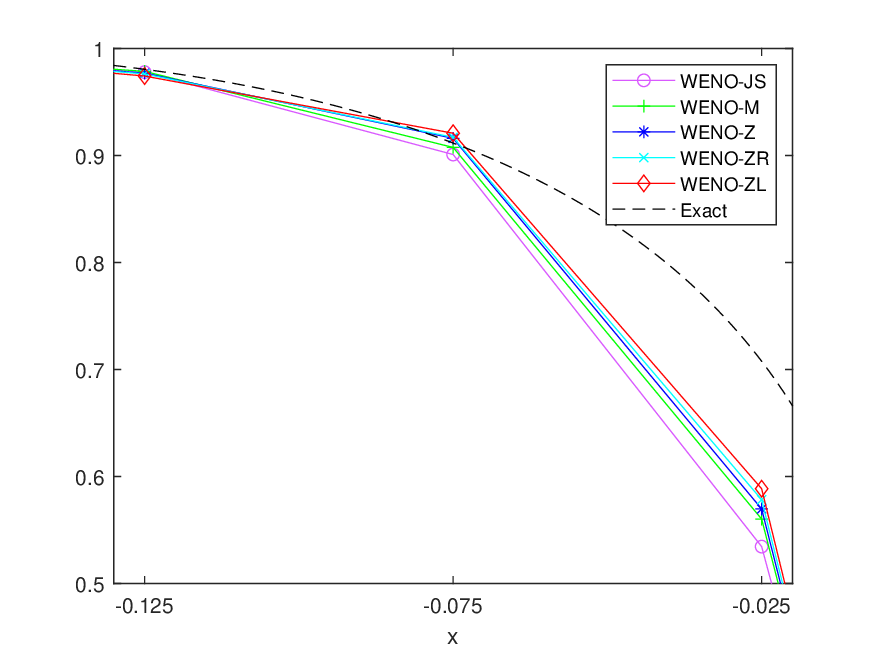}
\includegraphics[width=0.325\textwidth]{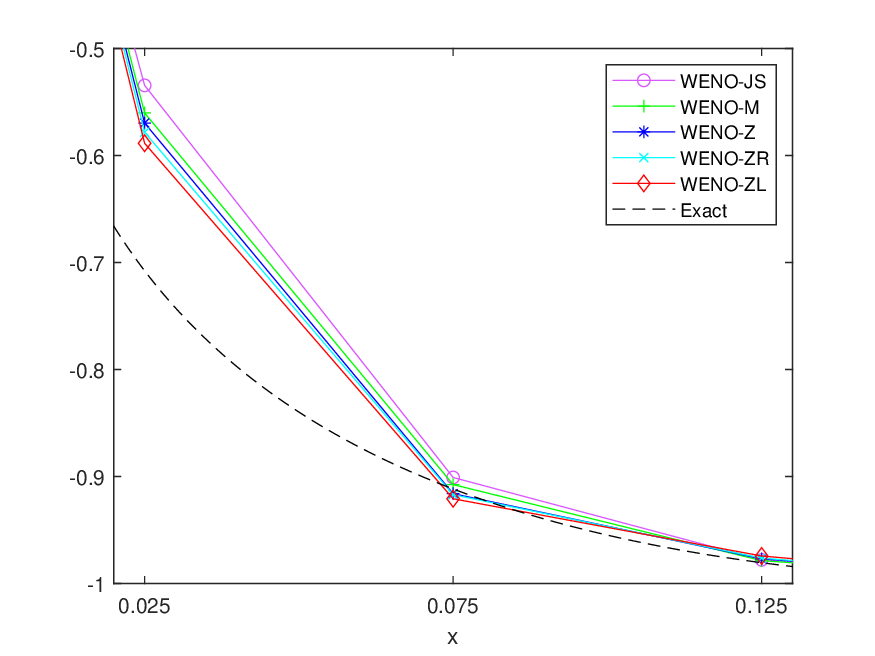}
\caption{Solution profiles for Example \ref{ex:burgers_1d} at $T = 1/\pi$ (left), close-up view of the solutions in the boxes (middle, right) computed by WENO-JS (purple), WENO-M (green), WENO-Z (blue), WENO-ZR (cyan) and WENO-ZL (red) with $N = 40$.
The dashed black lines are the exact solution.}
\label{fig:burgers_1d}
\end{figure}

\begin{example} \label{ex:burgers_1d}
We solve the Burgers' equation with the initial condition:
\begin{align*}
 u_t + \left(\frac{1}{2} u^2 \right)_x &= 0,\\
                                u(x,0) &= - \sin (\pi x).
\end{align*}  
The exact solution is smooth up to the final time $T = 1/\pi$.
We divide the computational domain $[-1, \, 1]$ into $N=40$ uniform cells.
Figure \ref{fig:burgers_1d} shows the solution profile at the final time by each scheme (WENO-ZL $p=5, \, q=1$).
We can see that WENO-JS, WENO-M, WENO-Z, WENO-ZR and WENO-ZL with less dissipation yield increasingly sharper approximations around the shock.
\end{example}

\begin{figure}[htbp]
\centering
\includegraphics[width=0.495\textwidth]{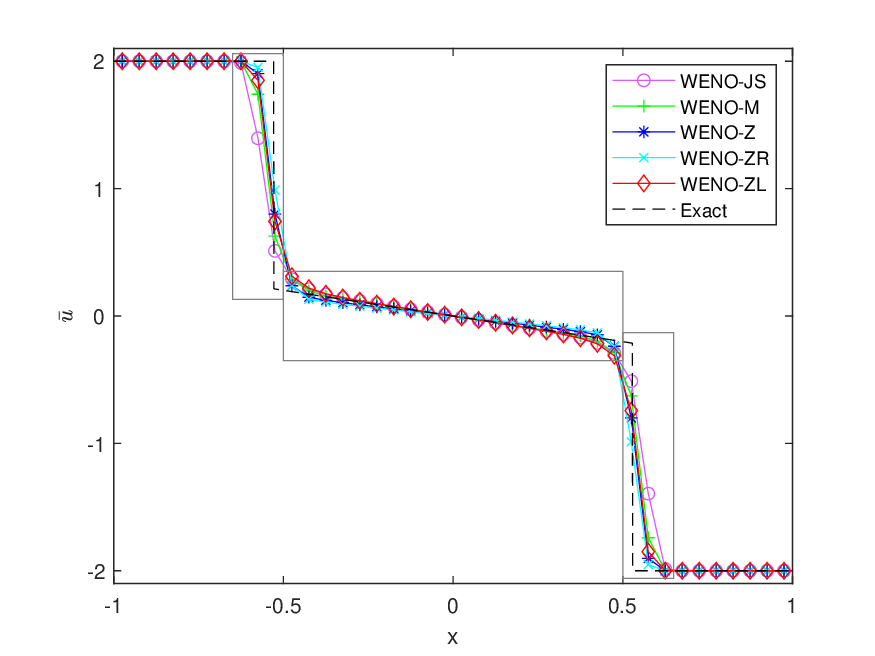}
\includegraphics[width=0.495\textwidth]{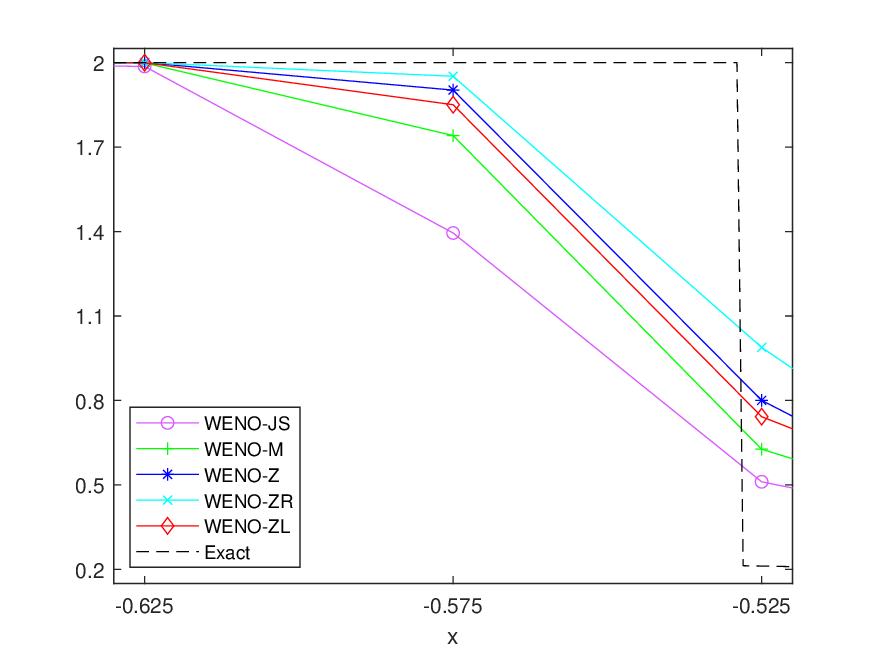}
\includegraphics[width=0.495\textwidth]{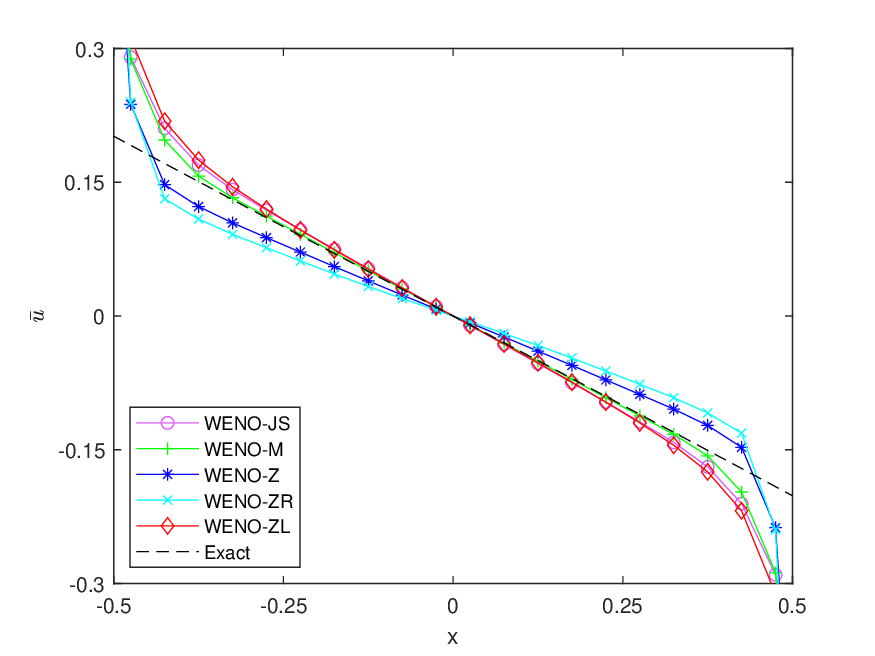}
\includegraphics[width=0.495\textwidth]{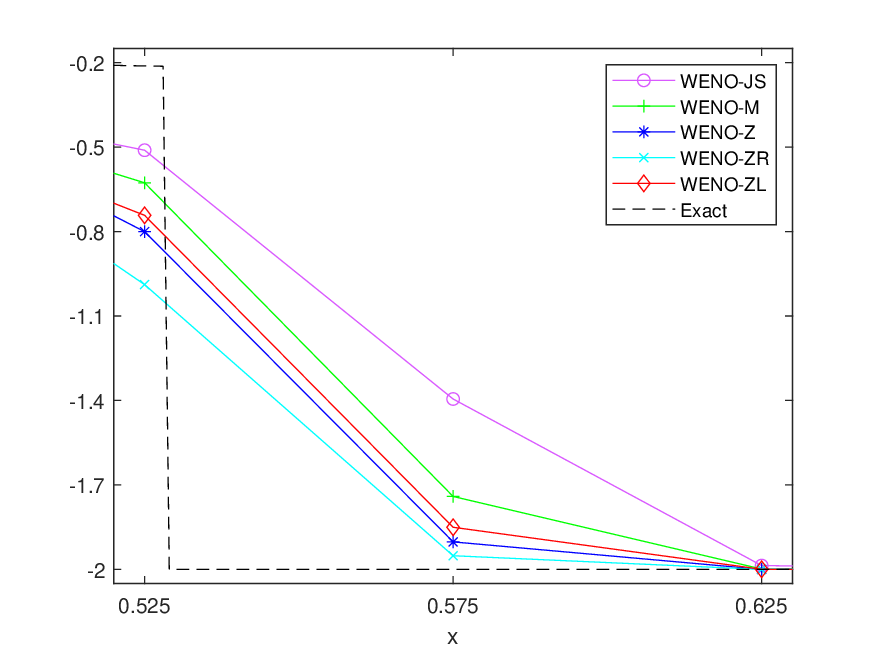}
\caption{Solution profiles for Example \ref{ex:nonconvex} with $u_L = 2$ and $u_R = -2$ at $T = 1$ (left), close-up view of the solutions in the boxes (top right, bottom left, bottom right) solved by WENO-JS (purple), WENO-M (green), WENO-Z (blue), WENO-ZR (cyan) and WENO-ZL (red) with $N = 40$.
The dashed black lines are the exact solution.}
\label{fig:nonconvex_IC0}
\end{figure}

\begin{figure}[htbp]
\centering
\includegraphics[width=0.495\textwidth]{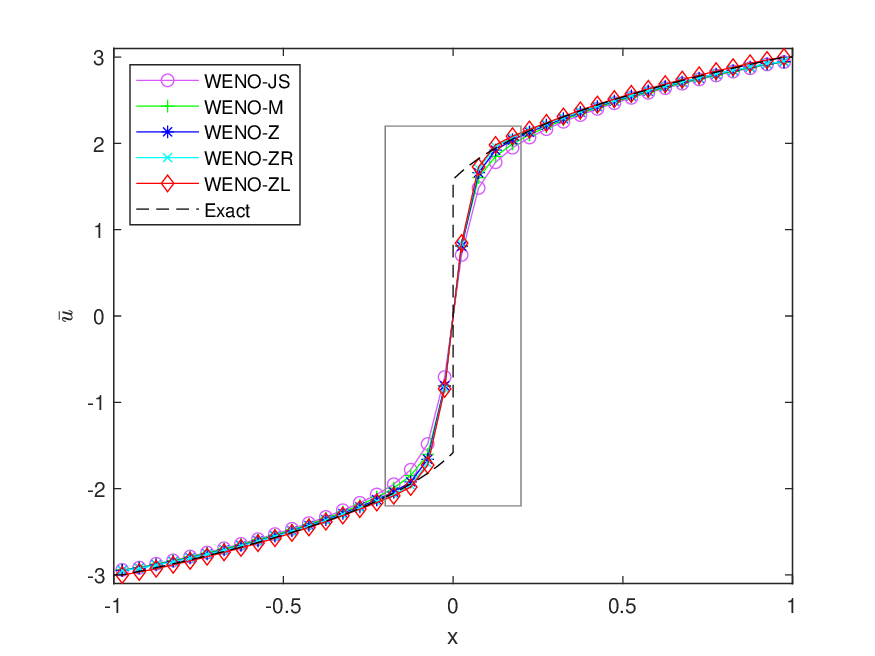}
\includegraphics[width=0.495\textwidth]{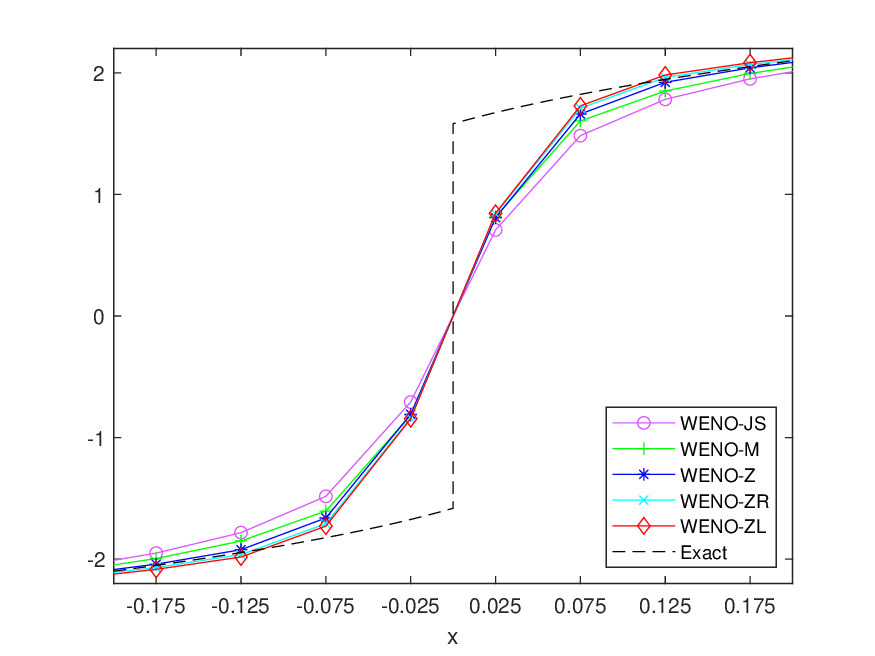}
\caption{Solution profiles for Example \ref{ex:nonconvex} with $u_L = -3$ and $u_R = 3$ at $T = 0.05$ (left), close-up view of the solutions in the box (right) solved by WENO-JS (purple), WENO-M (green), WENO-Z (blue), WENO-ZR (cyan) and WENO-ZL (red) with $N = 40$.
The dashed black lines are the exact solution.}
\label{fig:nonconvex_IC1}
\end{figure}

\begin{example} \label{ex:nonconvex}
We could test the performance of WENO schemes for the Riemann problems of the nonconvex flux  
$$
   f(u) = \frac{1}{4} (u^2-1)(u^2-4)
$$
with the initial condition 
$$
   u(x,0) = \left\{ 
             \begin{array}{ll} 
              u_L, & x \leqslant 0, \\ 
              u_R, & x > 0.
             \end{array} 
            \right.
$$
The computational domain $[-1, \, 1]$ is divided into $N=40$ uniform cells.
If $u_L = 2$ and $u_R = -2$, then the exact solution consists of two shocks and one rarefaction wave in between. 
The final time is $T=1$.
We present the numerical solutions at the final time with WENO-ZL ($p=q=2$), as shown in Figure \ref{fig:nonconvex_IC0}.
If $u_L = -3$ and $u_R = 3$, then the exact solution is a stationary shock at $x=0$.
Figure \ref{fig:nonconvex_IC1} shows the numerical results by five WENO schemes (WENO-ZL $p=2, \, q=1$) up to the final time $T=0.05$.
\end{example}

\begin{figure}[htbp]
\centering
\includegraphics[width=0.495\textwidth]{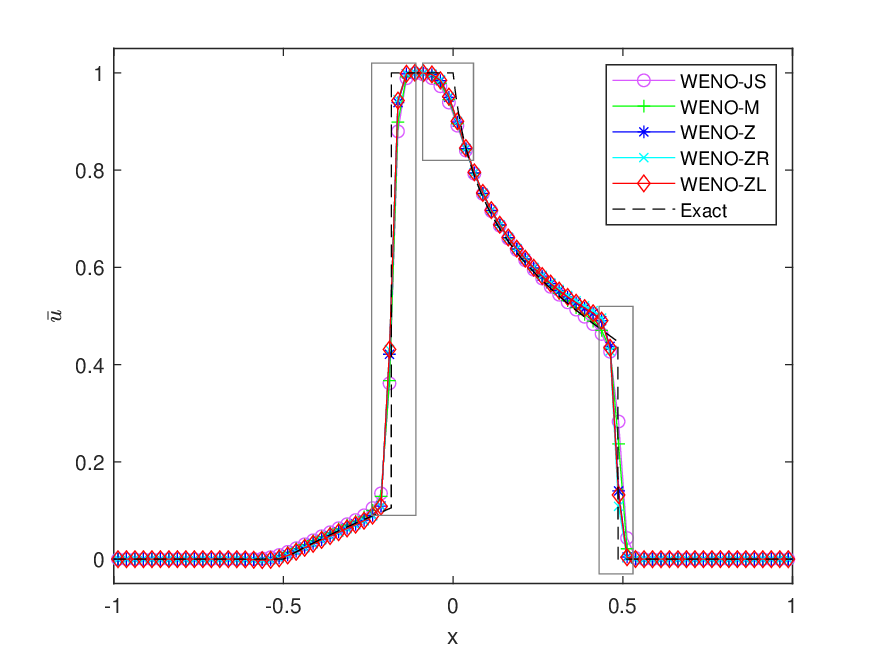}
\includegraphics[width=0.495\textwidth]{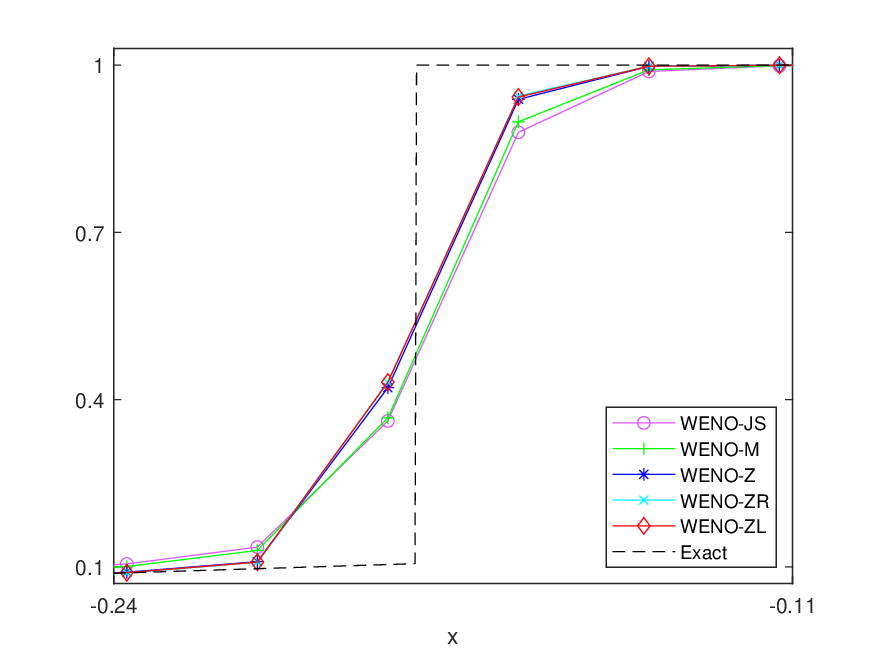}
\includegraphics[width=0.495\textwidth]{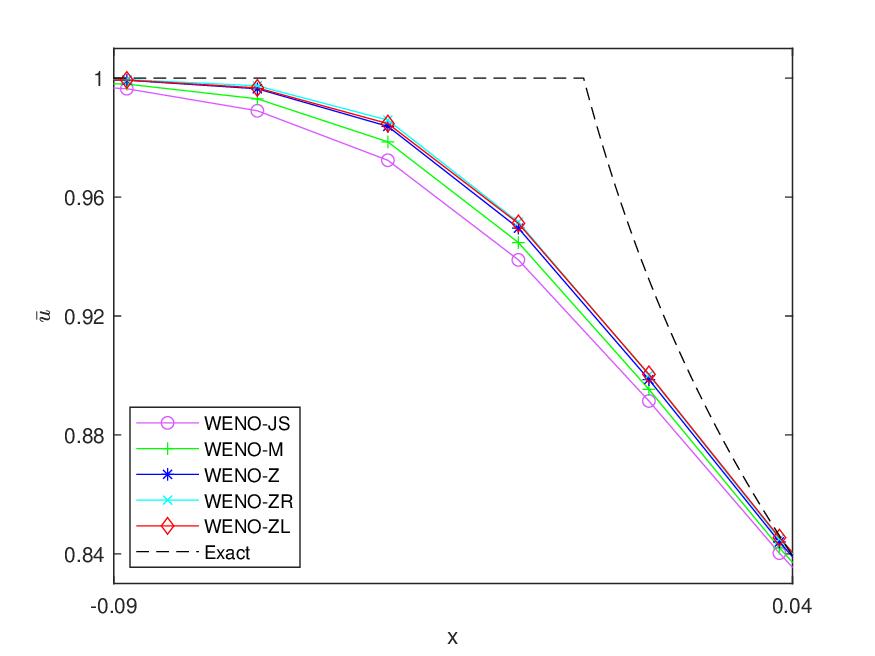}
\includegraphics[width=0.495\textwidth]{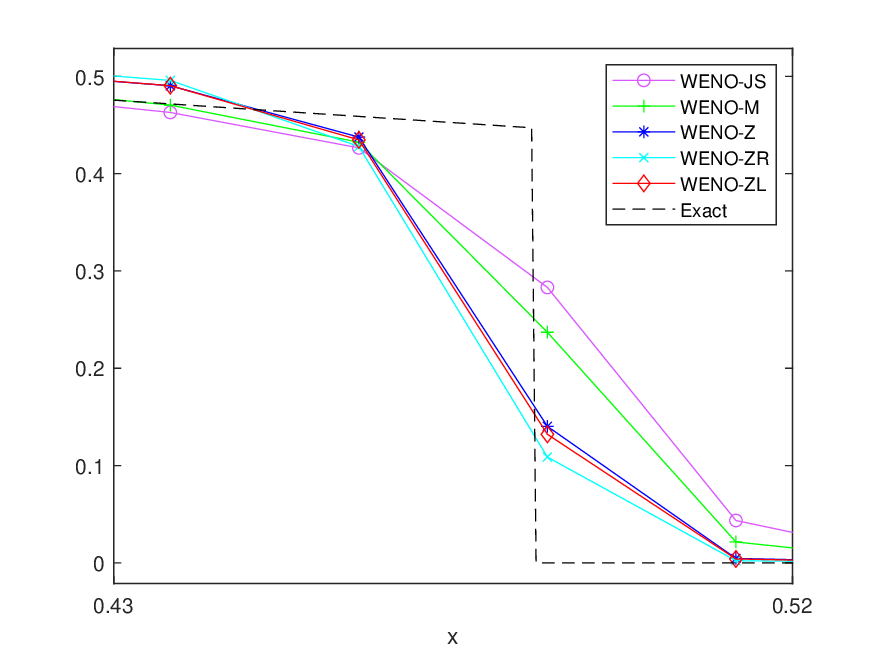}
\caption{Solution profiles for Example \ref{ex:buckley} at $T = 0.3$ (top left), close-up view of the solutions in the boxes from left to right (top right, bottom left, bottom right) calculated by WENO-JS (purple), WENO-M (green), WENO-Z (blue), WENO-ZR (cyan) and WENO-ZL (red) with $N = 80$.
The dashed black lines are the exact solution.}
\label{fig:buckley}
\end{figure}

\begin{example} \label{ex:buckley}
The Buckley-Leverett equation, a model for two-phase flow in porous media, is another typical problem of the nonconvex flux.
The flux has the form
$$
   f(u)= \frac{4u^2}{4u^2+(1-u)^2}.
$$
The initial condition is set as 
$$
   u(x,0) = \left\{ 
             \begin{array}{ll} 
              1, & -\frac{1}{2} \leqslant x \leqslant 0, \\ 
              0, & \text{otherwise}.
             \end{array} 
            \right.
$$
We divide the computational domain $[-1, \, 1]$ into $N=80$ uniform cells.
The final time is $T=0.3$.
We plot the numerical solutions at the final time, as shown in Figure \ref{fig:buckley}.
Around the shock front, WENO-Z, WENO-ZR and WENO-ZL ($p=q=1$) yield sharper solution than WENO-JS and WENO-M.
The scheme WENO-ZR gives more accurate solution profile in the regions of rarefaction and contact discontinuity.
\end{example}

\subsection{1D system case}
For one-dimensional system problem, we consider the Euler equations of gas dynamics
\begin{equation} \label{eq:euler_1d}
 \bfu_t + \bff(\bfu)_x = 0, 
\end{equation}
where the column vector $\bfu$ of the conserved variables and the flux vector $\bff$ in the $x$ direction are
$$
   \bfu = \left[ \rho,~ \rho u,~ E \right]^T, \quad \bff(\bfu) = \left[ \rho u,~ \rho u^2 + P,~ u(E+P) \right]^T,
$$
with $\rho, u$ and $P$ the primitive variables representing density, velocity and pressure, respectively.
The specific kinetic energy $E$ is 
$$
   E = \frac{P}{\gamma - 1} + \frac{1}{2} \rho u^2 
$$
with $\gamma = 1.4$ for the ideal gas. 

\begin{figure}[htbp]
\centering
\includegraphics[width=0.495\textwidth]{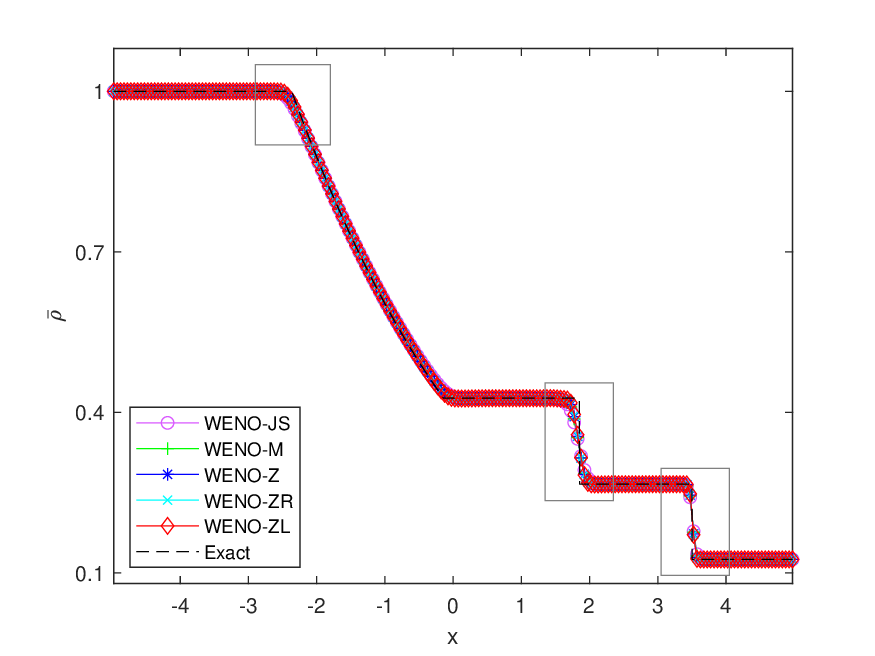}
\includegraphics[width=0.495\textwidth]{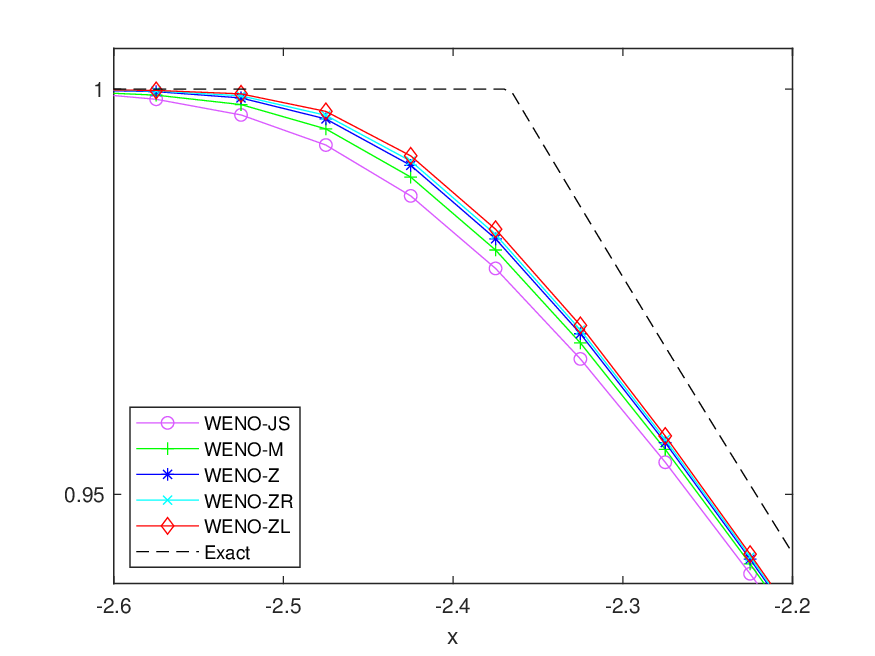}
\includegraphics[width=0.495\textwidth]{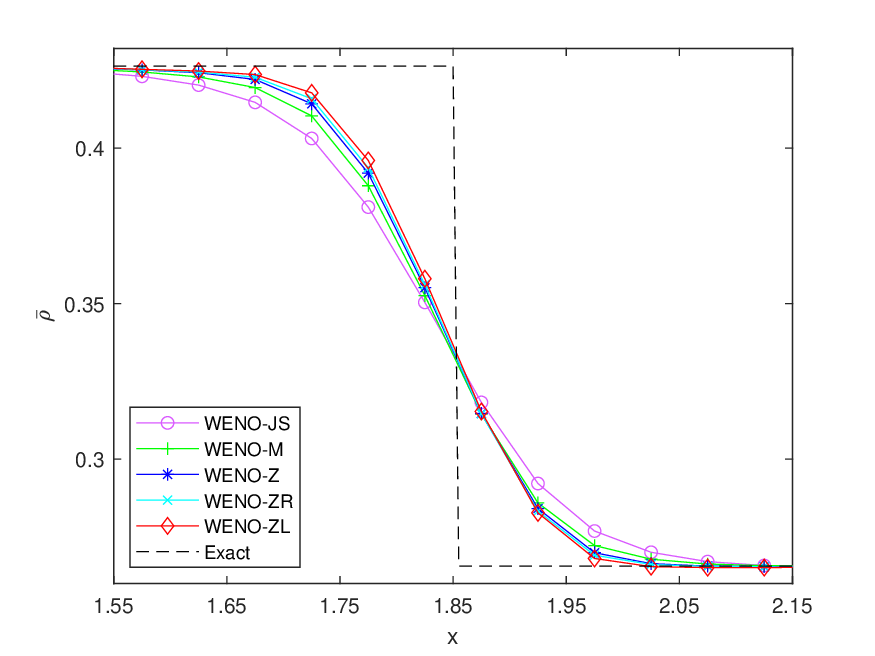}
\includegraphics[width=0.495\textwidth]{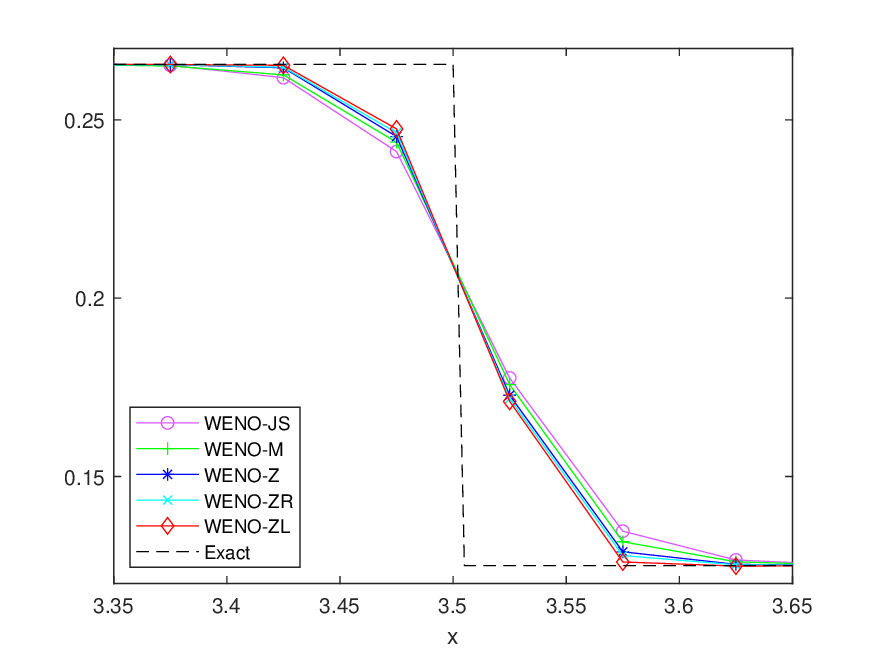}
\caption{Density profiles for the Sod problem \eqref{eq:euler_1d} and \eqref{eq:sod} at $T=2$ (top left), close-up view of the solutions in the boxes from left to right (top right, bottom left, bottom right) approximated by WENO-JS (purple), WENO-M (green), WENO-Z (blue), WENO-ZR (cyan) and WENO-ZL (red) with $N = 200$. 
The dashed black lines are the exact solution.}
\label{fig:sod}
\end{figure}

\begin{figure}[!t]
\centering
\includegraphics[width=0.495\textwidth]{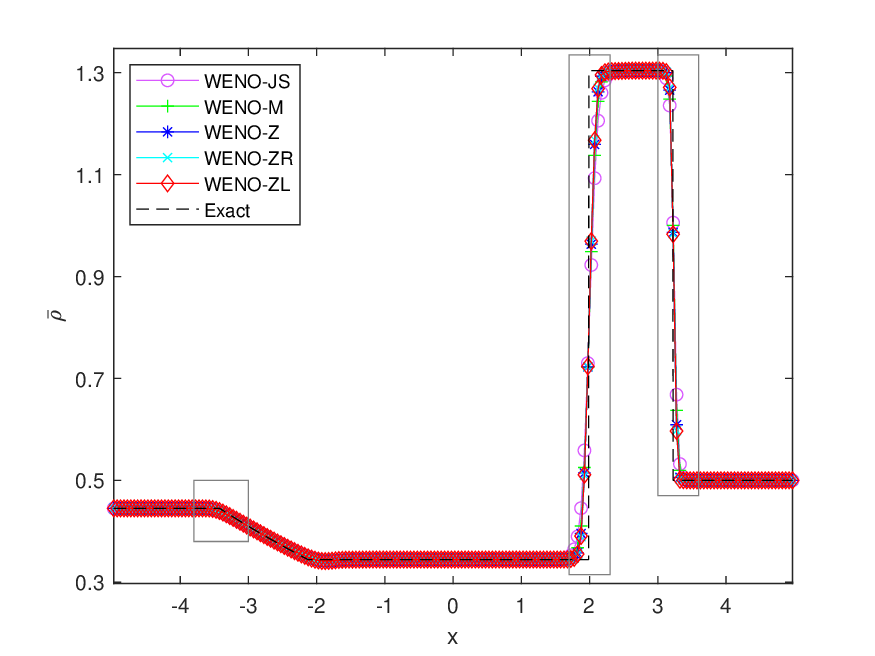}
\includegraphics[width=0.495\textwidth]{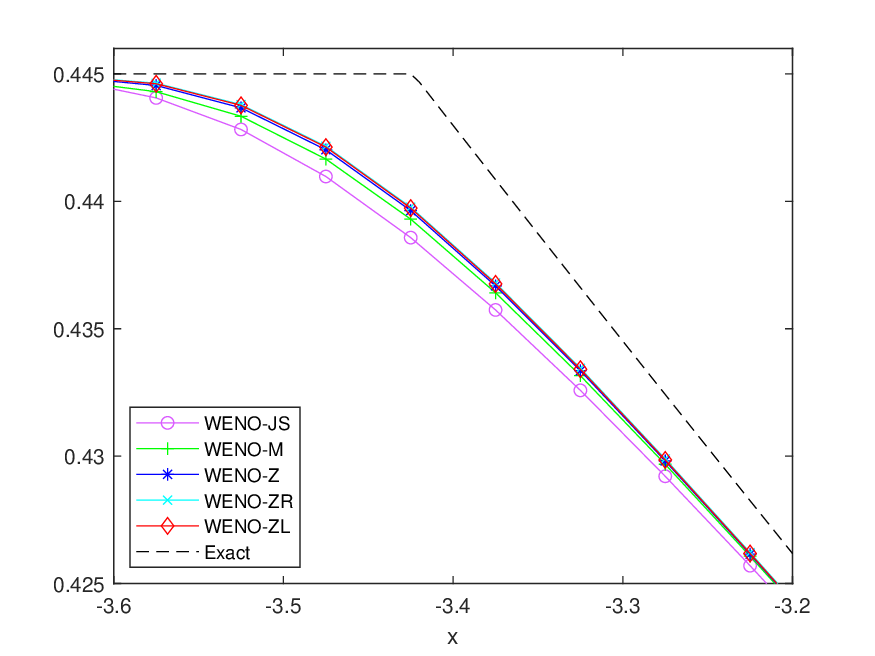}
\includegraphics[width=0.495\textwidth]{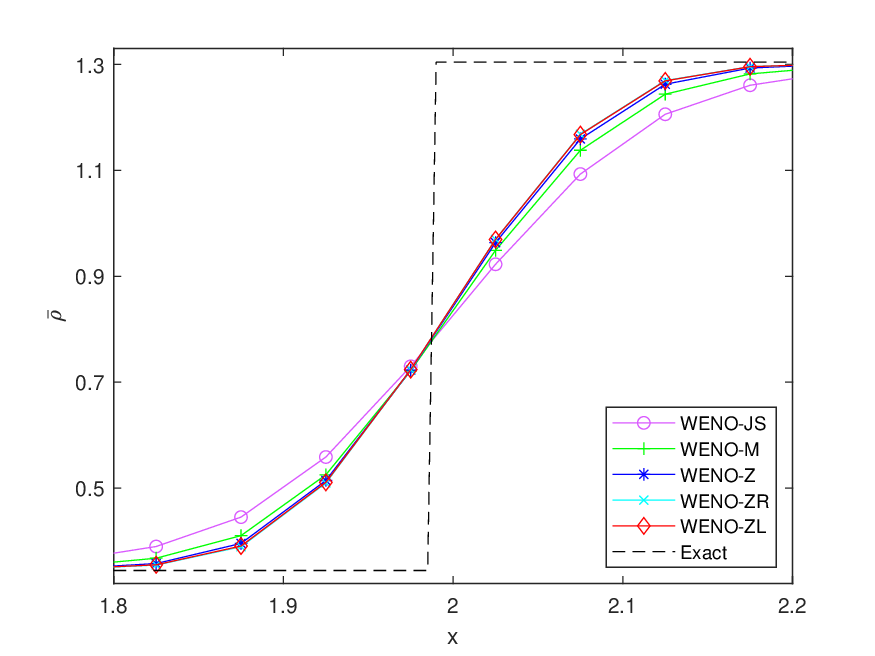}
\includegraphics[width=0.495\textwidth]{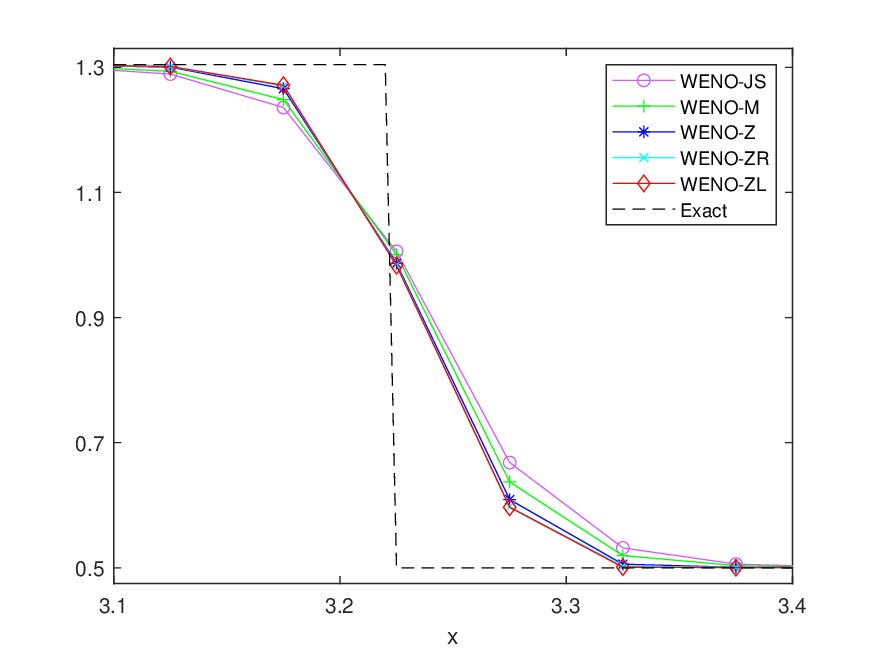}
\caption{Density profiles for the Lax problem \eqref{eq:euler_1d} and \eqref{eq:lax} at $T=1.3$ (top left), close-up view of the solutions in the boxes from left to right (top right, bottom left, bottom right) approximated by WENO-JS (purple), WENO-M (green), WENO-Z (blue), WENO-ZR (cyan) and WENO-ZL (red) with $N = 200$.
The dashed black lines are the exact solution.}
\label{fig:lax}
\end{figure}

\begin{example} \label{ex:euler_1d}
We begin with the standard shock tube problems.
The computational domain is $[-5, \, 5]$ with $N=200$ uniform cells.

The Sod problem consists of the Euler equations \eqref{eq:euler_1d} and the initial condition
\begin{equation} \label{eq:sod}
 (\rho, u, P ) = \left\{ 
                  \begin{array}{ll} 
                   (1,~0,~1),       & x \leqslant 0, \\ 
                   (0.125,~0,~0.1), & x > 0.
                  \end{array} 
                 \right.
\end{equation}
Figure \ref{fig:sod} presents the numerical density by the WENO schemes (WENO-ZL $p=5, \, q=1$) up to the final time $T=2$.
The density profile produced by WENO-ZL shows the less dissipation around the regions of rarefaction, contact discontinuity and shock wave than the other WENO schemes.

The Lax problem is the the Euler equations \eqref{eq:euler_1d} with the initial condition
\begin{equation} \label{eq:lax}
 (\rho, u, P ) = \left\{ 
                  \begin{array}{ll} 
                   (0.445,~0.698,~3.528), & x \leqslant 0, \\ 
                   (0.5,~0,~0.571),       & x > 0, 
                  \end{array} 
                 \right. 
\end{equation}
The final time is $T=1.3$.
We plot the numerical solutions of the density $\rho$ at the final time in Figure \ref{fig:lax}.
The density profile with WENO-ZL ($p=2, \, q=1$) is the most accurate near the rarefaction, contact discontinuity and shock regions.
\end{example}

\begin{figure}[h!]
\centering
\includegraphics[width=0.495\textwidth]{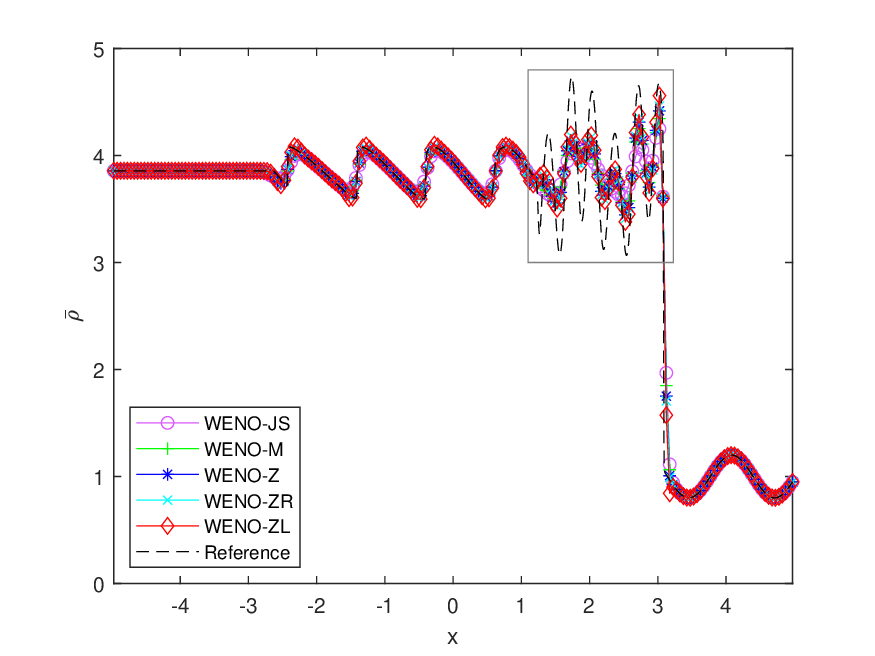}
\includegraphics[width=0.495\textwidth]{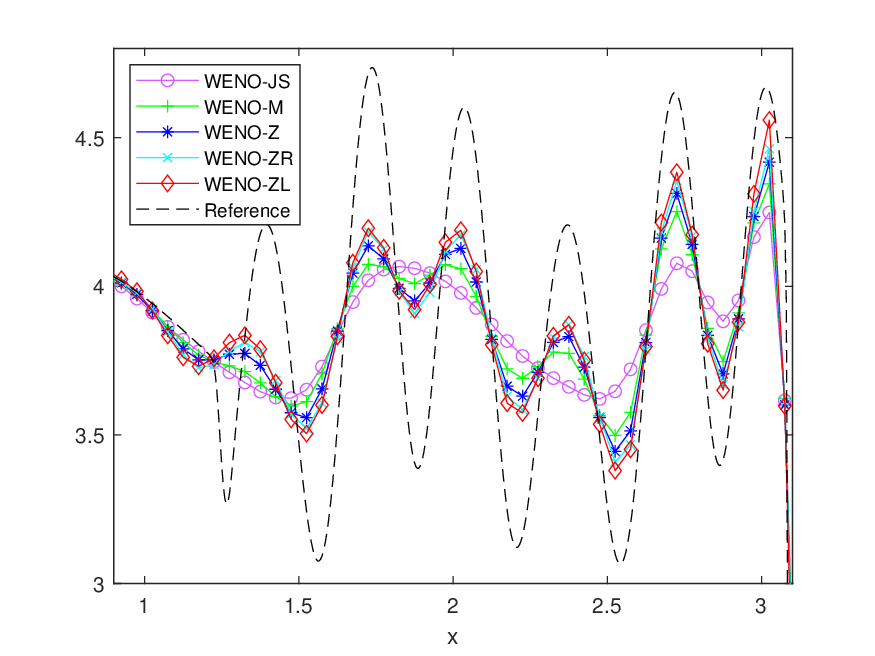}
\caption{Density profiles for Example \ref{ex:shock_entropy_wave} with $k=5$ at $T=2$ (left), close-up view of the solutions in the box (right) computed by WENO-JS (purple), WENO-M (green), WENO-Z (blue), WENO-ZR (cyan) and WENO-ZL (red) with $N = 200$.
The dashed black lines are generated by finite difference WENO-M with 2001 grid points.}
\label{fig:shock_entropy_wave_k5}
\end{figure}

\begin{figure}[h!]
\centering
\includegraphics[width=0.495\textwidth]{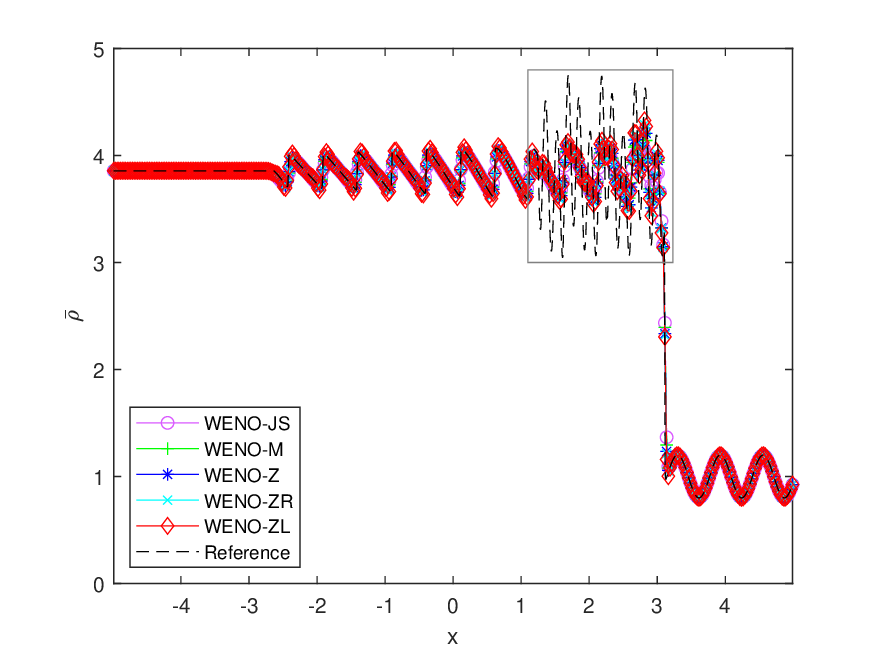}
\includegraphics[width=0.495\textwidth]{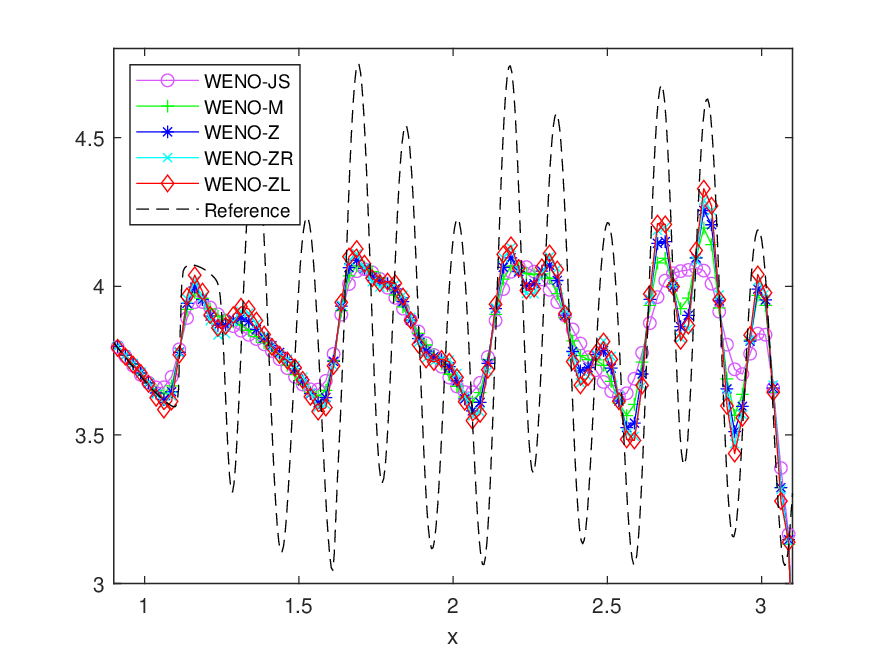}
\caption{Density profiles for Example \ref{ex:shock_entropy_wave} with $k=10$ at $T=2$ (left), close-up view of the solutions in the box (right) computed by WENO-JS (purple), WENO-M (green), WENO-Z (blue), WENO-ZR (cyan) and WENO-ZL (red) with $N = 400$.
The dashed black lines are generated by finite difference WENO-M with 2001 grid points.}
\label{fig:shock_entropy_wave_k10}
\end{figure}

\begin{example} \label{ex:shock_entropy_wave}
The shock entropy wave interaction problem \cite{ShuOsherI} is the Euler equations \eqref{eq:euler_1d} with a right moving Mach 3 shock and an entropy wave in density, of which the initial condition is given by
$$
   (\rho, u, P ) = \left\{ 
                    \begin{array}{ll} 
                     (3.857143,~2.629369,~10.333333), & x < -4, \\ 
                     (1 + 0.2 \sin(kx),~0,~1),        & x \geqslant -4,
                    \end{array} 
                   \right. 
$$
with $k$ the wave number of the entropy wave.

For $k=5$, we take a uniform grid with $N=200$ cells on the computational domain $[-5, \, 5]$. 
The numerical solution, computed by finite difference WENO-M with a high resolution of 2001 grid points, is used as the reference solution. 
The numerical solutions at $T=2$ are displayed in Figure \ref{fig:shock_entropy_wave_k5}.

For $k=10$, the computational domain $[-5, \, 5]$ is divided into $N=400$ uniform cells. 
We still take the numerical solution by finite difference WENO-M with 2001 points as the reference solution. 
Figure \ref{fig:shock_entropy_wave_k10} shows the approximate density profile at $T=2$.

We observe the improved performance of WENO-ZL ($p=5, \, q=1$) in capturing the fine structure of the density profile over the rest of WENO schemes.
\end{example}

\begin{figure}[h!]
\centering
\includegraphics[width=0.325\textwidth]{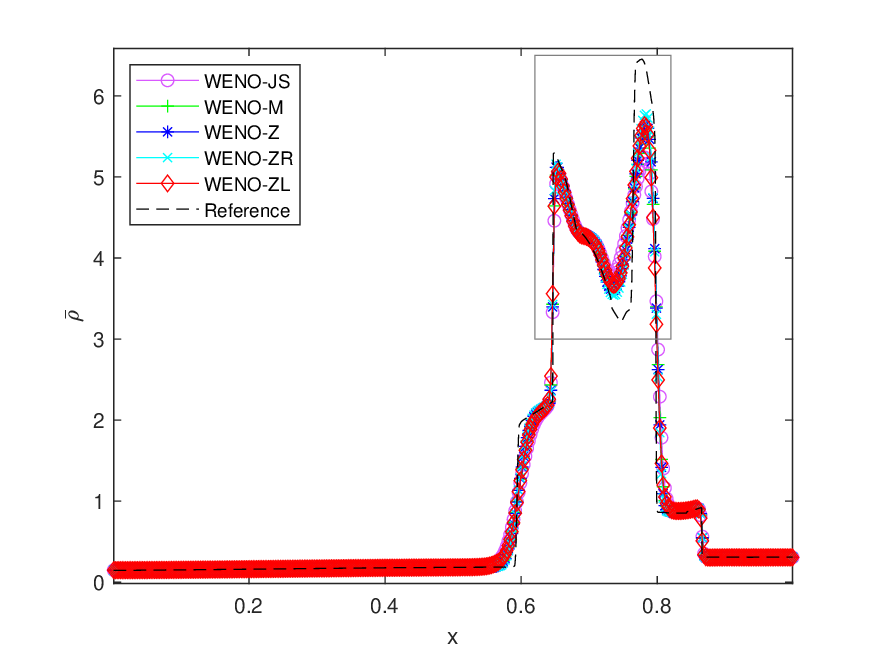}
\includegraphics[width=0.325\textwidth]{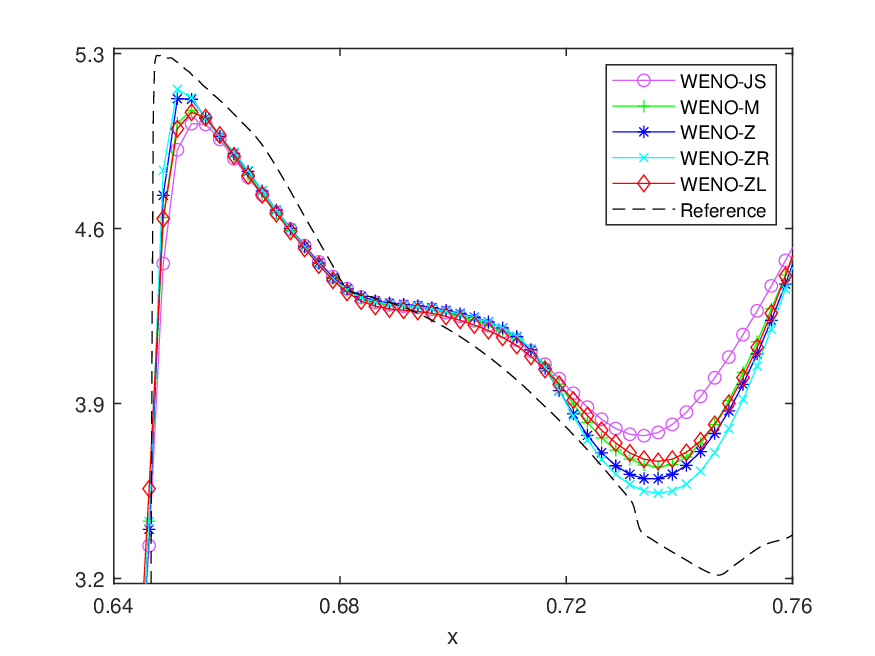}
\includegraphics[width=0.325\textwidth]{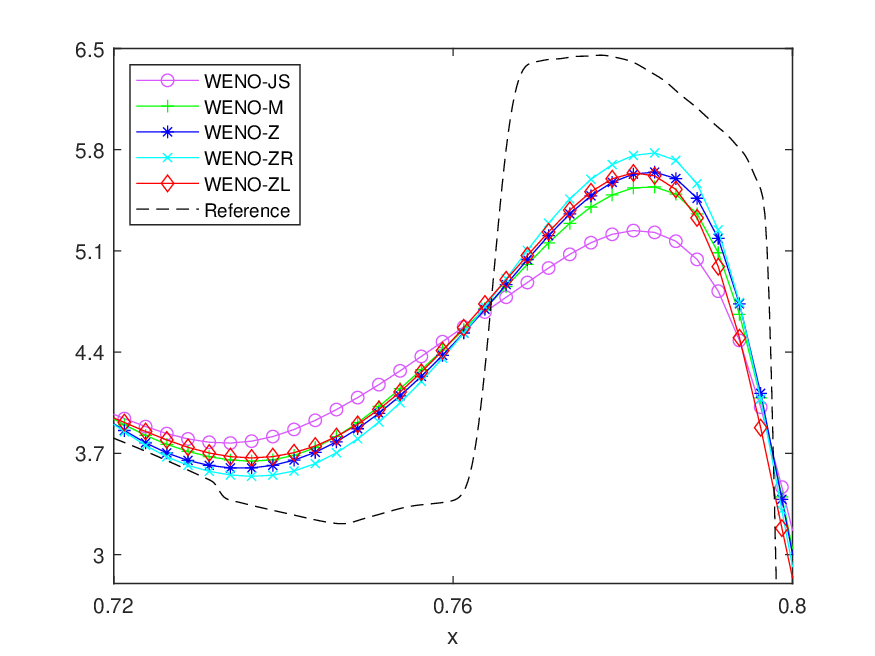}
\caption{Density profiles for Example \ref{ex:interacting_blastwave} at $T=0.038$ (left), close-up views of the solutions in the box (middle, right) computed by WENO-JS (purple), WENO-M (green), WENO-Z (blue), WENO-ZR (cyan) and WENO-ZL (red) with $N = 400$.
The dashed black lines are generated by FD WENO-M with $4001$ grid points.}
\label{fig:interacting_blastwave}
\end{figure}

\begin{example} \label{ex:interacting_blastwave}
Consider the blastwaves interaction problem, which is composed of the Euler equations \eqref{eq:euler_1d} and the initial condition
$$
   (\rho, u, P ) = \left\{ 
                    \begin{array}{ll} 
                     (1,~0,~1000), & 0 \leqslant x < 0.1, \\ 
                     (1,~0,~0.01), & 0.1 \leqslant x < 0.9, \\
                     (1,~0,~100),  & 0.9 \leqslant x \geqslant 1.
                    \end{array} 
                   \right. 
$$
We divide the computational domain $[0, \, 1]$ into $N=400$ uniform cells.
Notice that it is easy to produce the negative number of pressure during the simulation and thus the complex number of sound speed if we use larger $p$ and smaller $q$.
Then we apply $p=1/7, \, q=2$ for WENO-ZL to ensure the positivity of the pressure.
The numerical solutions of the density $rho$ at the final time $T=0.038$ are shown in Figure \ref{fig:interacting_blastwave}, where the dashed lines are the approximations using finite difference WENO-M with $4001$ grid points for the reference solution.
It is worth mentioning that WENO-ZR has better resolution for density than WENO-ZL as we sacrifice the accuracy for the positivity of the pressure.
Relaxing the fixed value of $p$ in each cell without using the positivity-preserving method for better accuracy is still under investigation. 
\end{example}

\subsection{2D scalar case}
\begin{example} \label{ex:advection_2d_order}
We consider the two-dimensional linear advection equation
$$
   u_t + u_x + u_y = 0,\quad -1 \leqslant x, y \leqslant 1,
$$
with the initial condition
$$
   u(x,y,0) = \left\{ 
                   \begin{array}{ll} 
                    1, & \text{if } (x,y) \in S, \\ 
                    0, & \text{otherwise},
                   \end{array}
                   \right.
$$
with $S = \{ (x,y) | |x \pm y| < 1/\sqrt{2} \}$ a unit square centered at the origin, and the periodic boundary condition. 
We use $\Delta x = \Delta y$ and run the WENO schemes up to the final time $T=4$ with the time step $\Delta t = \cfl \cdot \Delta x$.

The $L_1, L_2$ and $L_{\infty}$ errors versus $N$, as well as the order of accuracy, for the WENO-JS, WENO-M, WENO-Z, WENO-ZR and WENO-ZL ($p=5, \, q=1$) schemes are provided in Tables \ref{tab:advection_2d_order_L1}, \ref{tab:advection_2d_order_L2} and \ref{tab:advection_2d_order_Linf}, respectively. 
We see that WENO-ZL performs slightly better than the other WENO schemes in terms of $L_1$ and $L_2$ errors, while WENO-ZR behaves better in terms of $L_{\infty}$ accuracy.
\end{example}

\begin{table}[h!]
\renewcommand{\arraystretch}{1.1}
\scriptsize
\centering
\caption{$L_1$ error and order of accuracy for Example \ref{ex:advection_2d_order}.}      
\begin{tabular}{clcrlcrlcrlcrlc} 
\hline
N & \multicolumn{2}{l}{WENO-JS} & & \multicolumn{2}{l}{WENO-M} & & \multicolumn{2}{l}{WENO-Z} & & \multicolumn{2}{l}{WENO-ZR} & & \multicolumn{2}{l}{WENO-ZL} \\ 
    \cline{2-3}                     \cline{5-6}                    \cline{8-9}                    \cline{11-12}                   \cline{14-15}
  & Error & Order               & & Error & Order              & & Error & Order              & & Error & Order               & & Error & Order \\
\hline
10  & 1.23e-1 & --     & & 1.09e-1 & --     & & 1.07e-1 & --     & & 1.04e-1 & --     & & 1.04e-1 & --     \\  
20  & 8.16e-2 & 0.5925 & & 6.70e-2 & 0.6951 & & 6.26e-2 & 0.7701 & & 6.03e-2 & 0.7907 & & 5.98e-2 & 0.7980 \\  
40  & 4.98e-2 & 0.7123 & & 4.01e-2 & 0.7422 & & 3.72e-2 & 0.7512 & & 3.60e-2 & 0.7421 & & 3.60e-2 & 0.7319 \\
80  & 2.99e-2 & 0.7356 & & 2.43e-2 & 0.7219 & & 2.31e-2 & 0.6866 & & 2.26e-2 & 0.6739 & & 2.27e-2 & 0.6662 \\ 
160 & 1.75e-2 & 0.7734 & & 1.43e-2 & 0.7627 & & 1.40e-2 & 0.7246 & & 1.38e-2 & 0.7140 & & 1.39e-2 & 0.7078 \\
\hline
\end{tabular}
\label{tab:advection_2d_order_L1}
\end{table}

\begin{table}[h!]
\renewcommand{\arraystretch}{1.1}
\scriptsize
\centering
\caption{$L_2$ error and order of accuracy for Example \ref{ex:advection_2d_order}.}      
\begin{tabular}{clcrlcrlcrlcrlc} 
\hline  
N & \multicolumn{2}{l}{WENO-JS} & & \multicolumn{2}{l}{WENO-M} & & \multicolumn{2}{l}{WENO-Z} & & \multicolumn{2}{l}{WENO-ZR} & & \multicolumn{2}{l}{WENO-ZL} \\ 
    \cline{2-3}                     \cline{5-6}                    \cline{8-9}                    \cline{11-12}                   \cline{14-15}
  & Error & Order               & & Error & Order              & & Error & Order              & & Error & Order               & & Error & Order \\
\hline 
10  & 1.73e-1 & --     & & 1.56e-1 & --     & & 1.53e-1 & --     & & 1.49e-1 & --     & & 1.48e-1 & --     \\  
20  & 1.35e-1 & 0.3596 & & 1.18e-1 & 0.4127 & & 1.13e-1 & 0.4420 & & 1.10e-1 & 0.4429 & & 1.09e-1 & 0.4434 \\  
40  & 1.06e-1 & 0.3541 & & 9.21e-2 & 0.3523 & & 8.88e-2 & 0.3443 & & 8.70e-2 & 0.3324 & & 8.66e-2 & 0.3255 \\
80  & 8.29e-2 & 0.3486 & & 7.33e-2 & 0.3293 & & 7.15e-2 & 0.3128 & & 7.05e-2 & 0.3032 & & 7.01e-2 & 0.3043 \\ 
160 & 6.44e-2 & 0.3645 & & 5.75e-2 & 0.3490 & & 5.68e-2 & 0.3335 & & 5.62e-2 & 0.3279 & & 5.59e-2 & 0.3268 \\   
\hline
\end{tabular}
\label{tab:advection_2d_order_L2}
\end{table}

\begin{table}[h!]
\renewcommand{\arraystretch}{1.1}
\scriptsize
\centering
\caption{$L_{\infty}$ error and order of accuracy for Example \ref{ex:advection_2d_order}.}      
\begin{tabular}{clcrlcrlcrlcrlc} 
\hline  
N & \multicolumn{2}{l}{WENO-JS} & & \multicolumn{2}{l}{WENO-M} & & \multicolumn{2}{l}{WENO-Z} & & \multicolumn{2}{l}{WENO-ZR} & & \multicolumn{2}{l}{WENO-ZL} \\ 
    \cline{2-3}                     \cline{5-6}                    \cline{8-9}                    \cline{11-12}                   \cline{14-15}
  & Error & Order               & & Error & Order              & & Error & Order              & & Error & Order               & & Error & Order \\
\hline
10  & 4.26e-1 & --      & & 3.77e-1 & --      & & 3.65e-1 & --      & & 3.52e-1 & --      & & 3.50e-1 & --     \\  
20  & 4.21e-1 &  0.0154 & & 3.63e-1 &  0.0560 & & 3.52e-1 &  0.0540 & & 3.40e-1 &  0.0505 & & 3.43e-1 &  0.0285 \\ 
40  & 4.57e-1 & -0.1160 & & 3.95e-1 & -0.1225 & & 3.86e-1 & -0.1315 & & 3.77e-1 & -0.1468 & & 3.81e-1 & -0.1515 \\
80  & 4.71e-1 & -0.0441 & & 4.05e-1 & -0.0369 & & 4.03e-1 & -0.0623 & & 3.94e-1 & -0.0643 & & 4.02e-1 & -0.0767 \\  
160 & 5.43e-1 & -0.2049 & & 5.02e-8 & -0.3070 & & 5.02e-1 & -0.3181 & & 4.98e-1 & -0.3390 & & 5.09e-1 & -0.3402 \\ 
\hline
\end{tabular}
\label{tab:advection_2d_order_Linf}
\end{table}

\begin{figure}[htbp]
\centering
\includegraphics[width=0.495\textwidth]{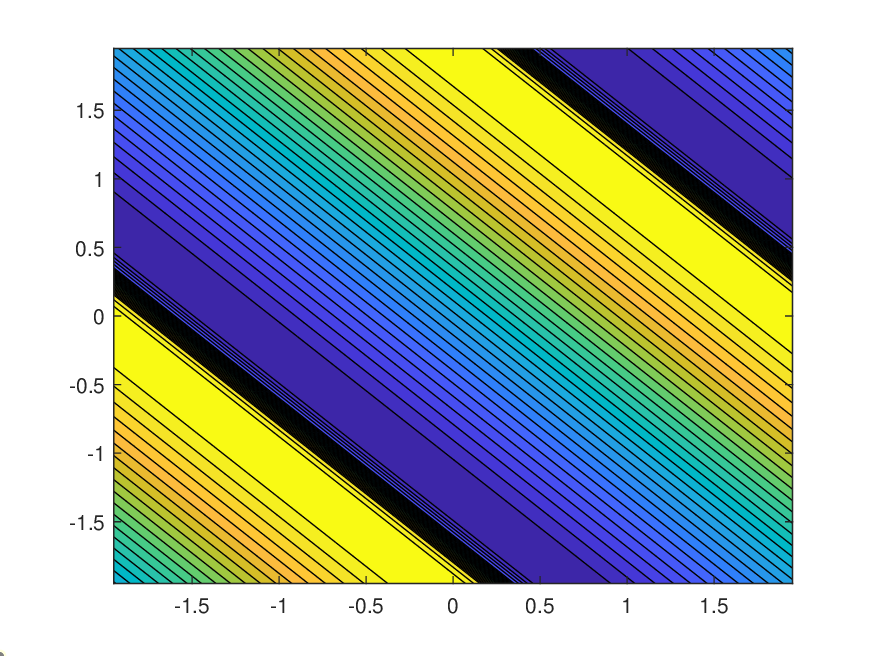}
\includegraphics[width=0.495\textwidth]{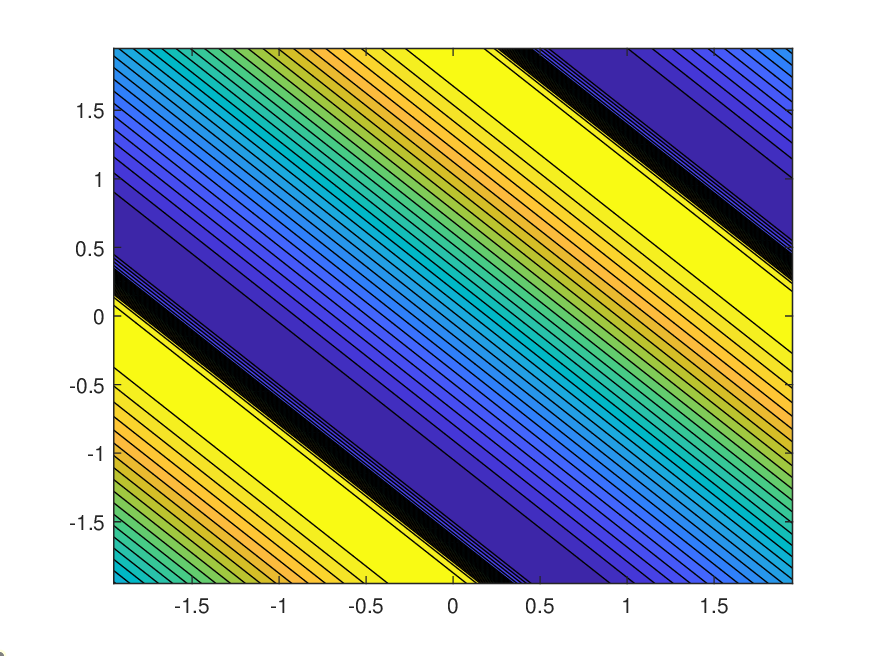}
\includegraphics[width=0.495\textwidth]{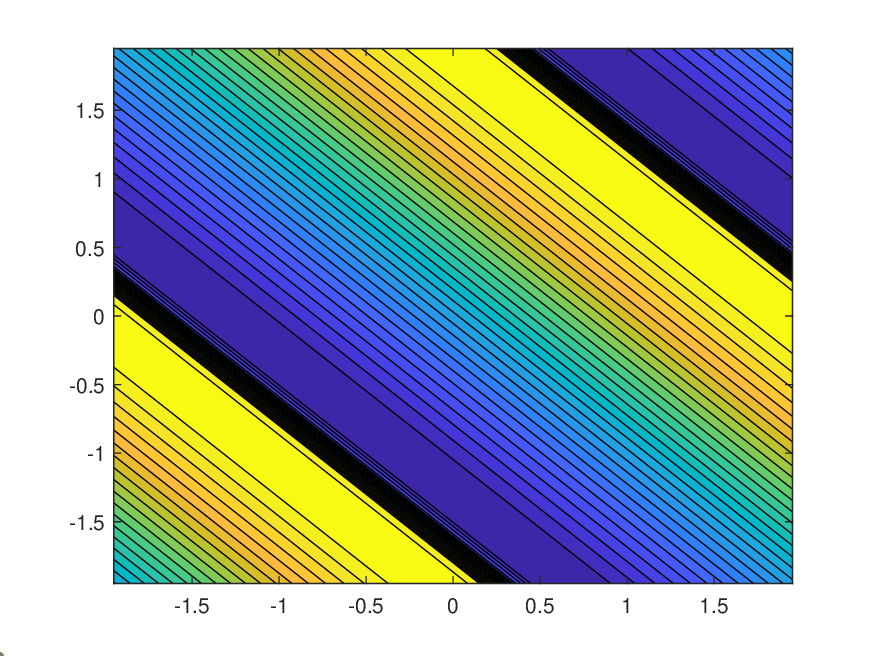}
\includegraphics[width=0.495\textwidth]{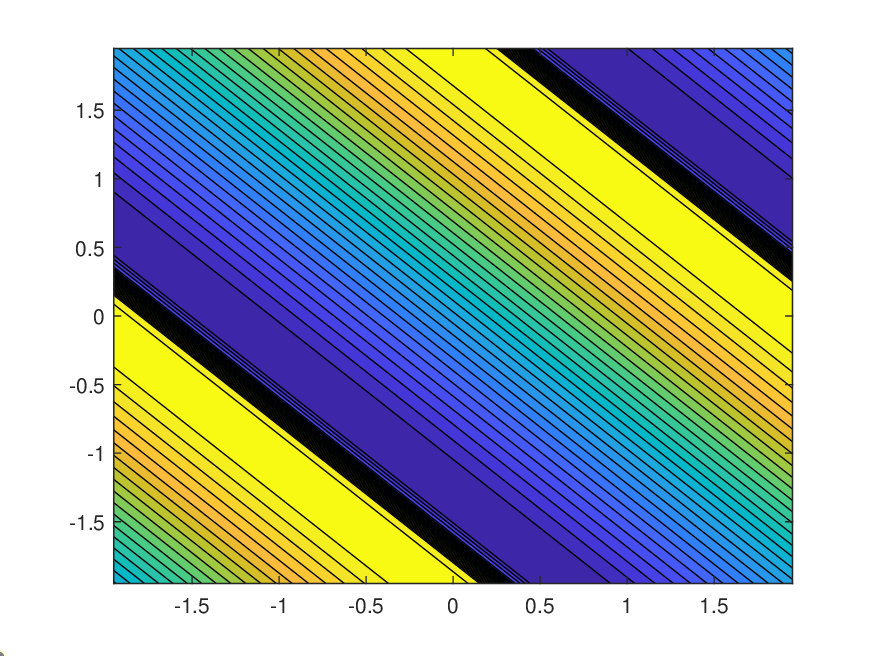}
\includegraphics[width=0.495\textwidth]{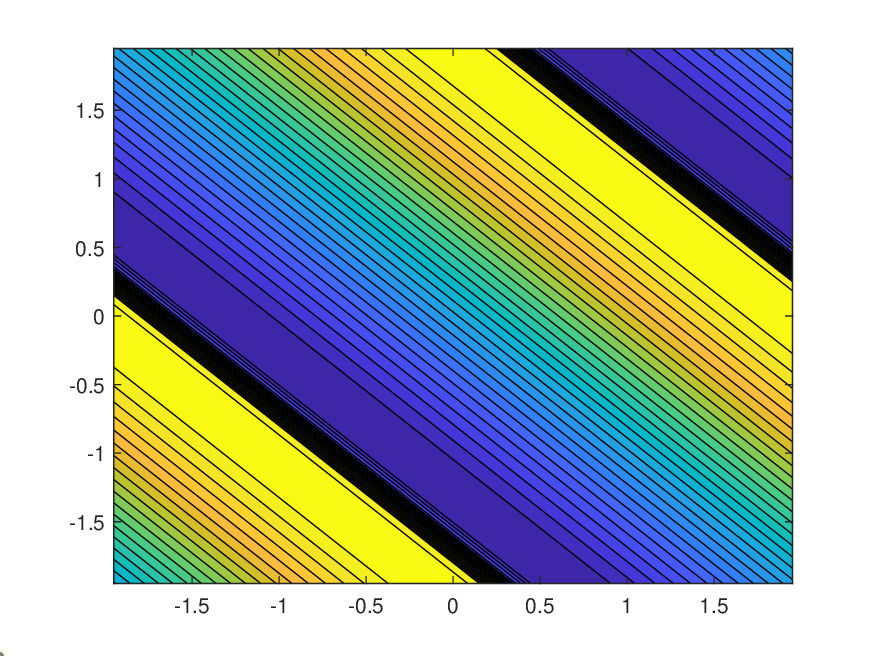}
\includegraphics[width=0.495\textwidth]{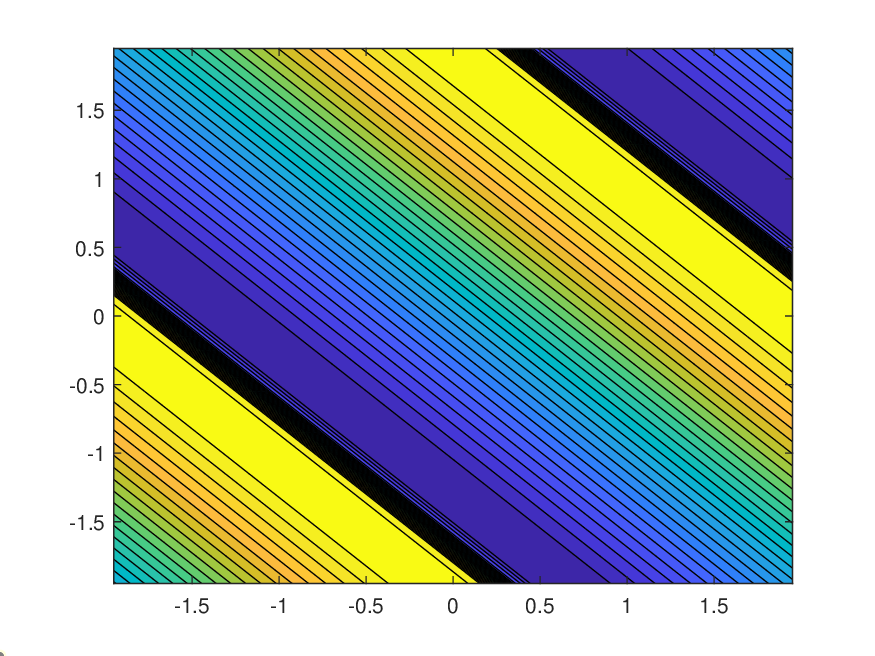}
\caption{Solution in the filled contour plot for Example \ref{ex:burgers_2d} at $T=2/\pi$ by WENO-JS (top left), WENO-M (top right), WENO-Z (middle left), WENO-ZR (middle right), WENO-ZL ($p=q=1$) (bottom left) and WENO-ZL ($p=3,\, q=1$) (bottom right) with $N_x = N_y = 40$.
Each contour plot displays contours at $30$ levels of the density.}
\label{fig:burgers_2d}
\end{figure}

\begin{example} \label{ex:burgers_2d}
Consider the two-dimensional Burgers' equation
$$
   u_t + \left(\frac{1}{2} u^2 \right)_x + \left(\frac{1}{2} u^2 \right)_y = 0
$$
with the initial condition
$$
   u(x,y,0) = \frac{1}{4} + \frac{1}{2} \sin \left( \pi \frac{x+y}{2} \right),
$$
and the periodic boundary condition in each direction. 
The exact solution is smooth up to the final time $T = 2/\pi$.
We divide the computational domain $[-2, \, 2] \times [-2, \, 2]$ into $N_x \times N_y = 40 \times 40$ uniform cells.
The combinations $(p,\, q)=(1,\, 1)$ and $(p,\, q)=(3\, ,1)$ are chosen for WENO-ZL.
We observe that all WENO schemes generate comparable numerical results at the final time in Figure \ref{fig:burgers_2d}, which indicates that our WENO-ZL schemes work well in multi-dimensional problems. 
\end{example}

\begin{figure}[htbp]
\centering
\includegraphics[width=0.495\textwidth]{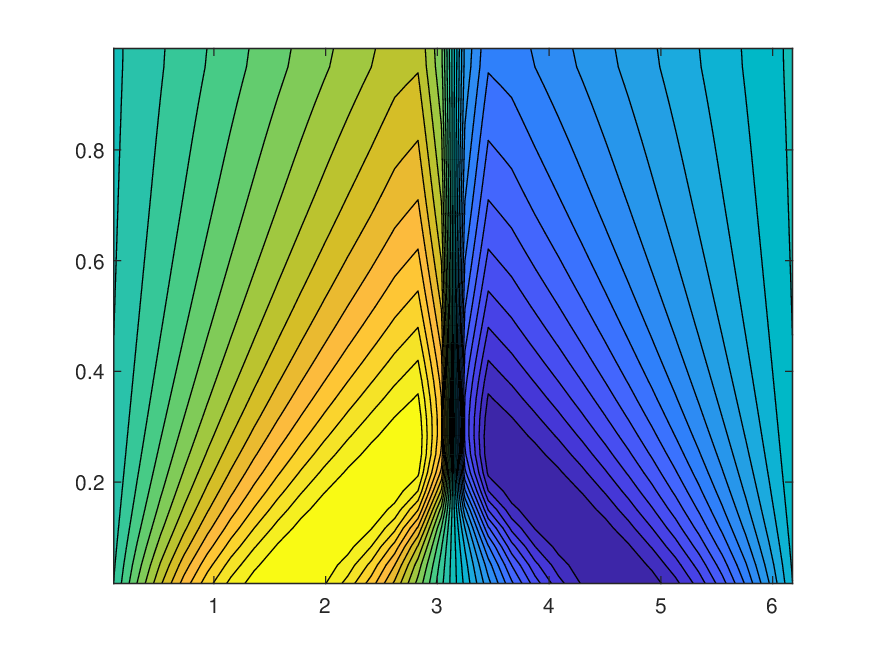}
\includegraphics[width=0.495\textwidth]{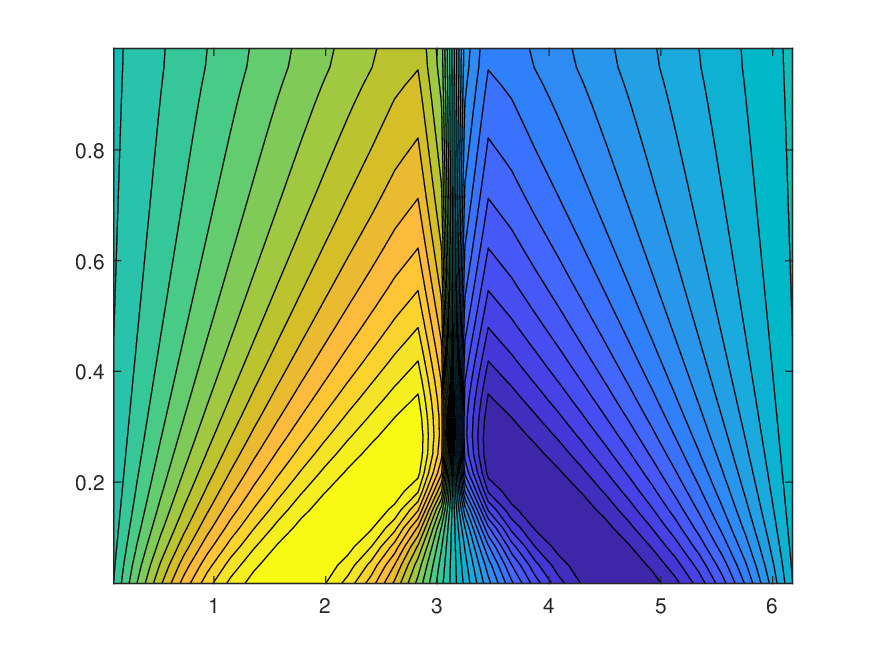}
\includegraphics[width=0.495\textwidth]{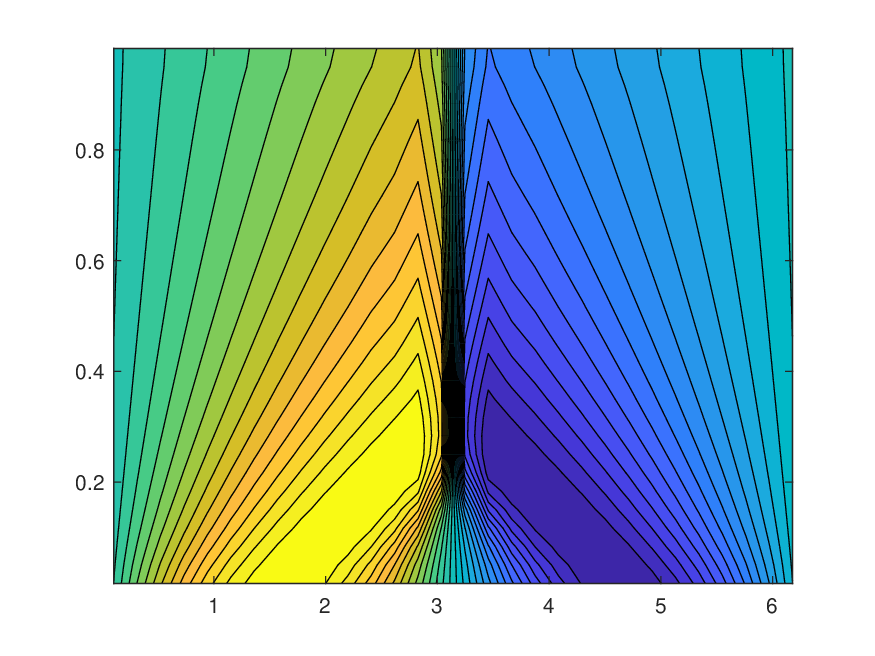}
\includegraphics[width=0.495\textwidth]{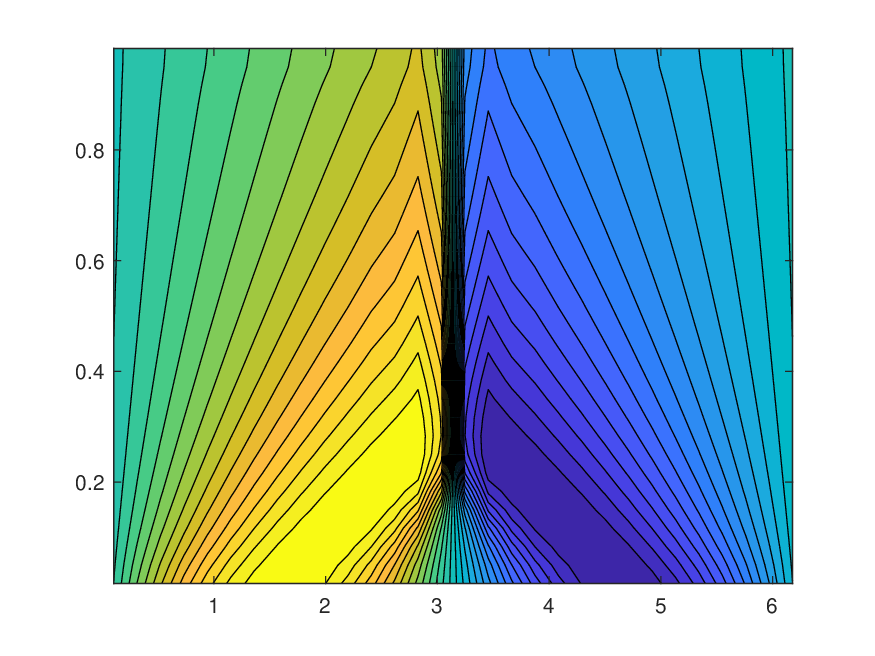}
\includegraphics[width=0.495\textwidth]{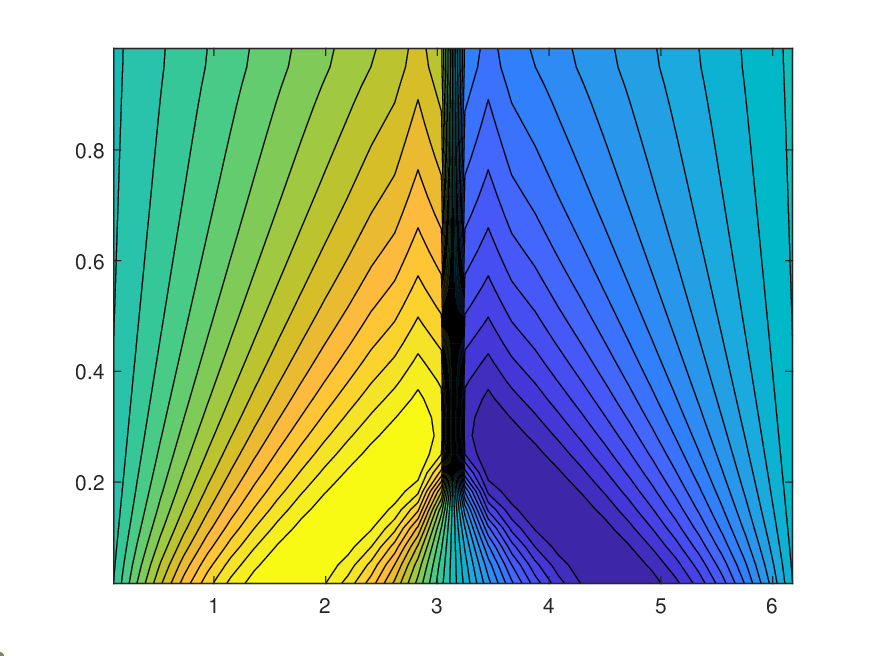}
\includegraphics[width=0.495\textwidth]{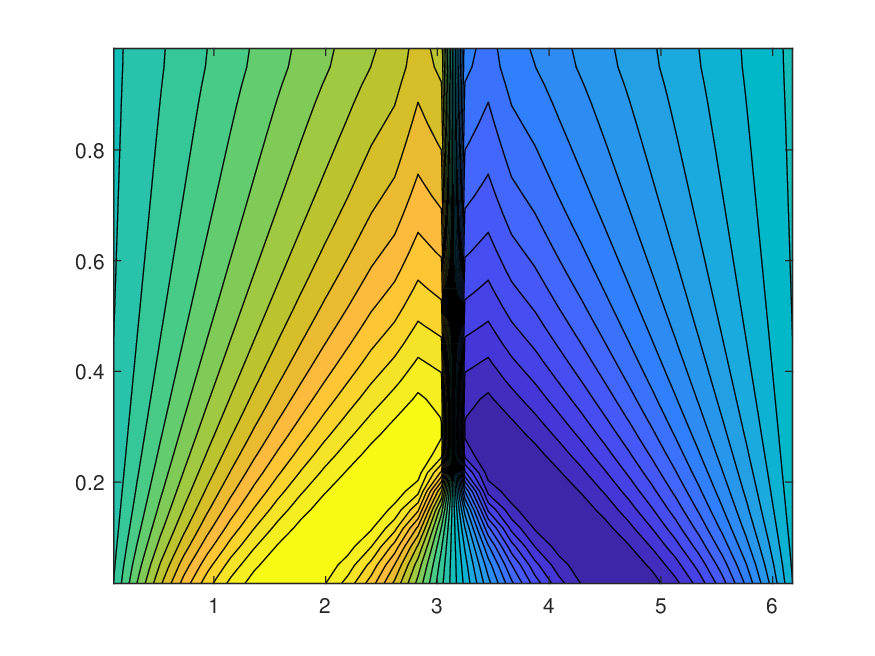}
\caption{Solution in the filled contour plot for Example \ref{ex:burgers_advection} with $\alpha=0,\, \beta=5$ at $T=1$ by WENO-JS (top left), WENO-M (top right), WENO-Z (middle left), WENO-ZR (middle right), WENO-ZL ($p=q=1$) (bottom left) and WENO-ZL ($p=3,\, q=1$) (bottom right) with $N_x = N_y = 30$.
Each contour plot displays contours at $30$ levels of the density.}
\label{fig:burgers_advection_IC1}
\end{figure}

\begin{figure}[htbp]
\centering
\includegraphics[width=0.495\textwidth]{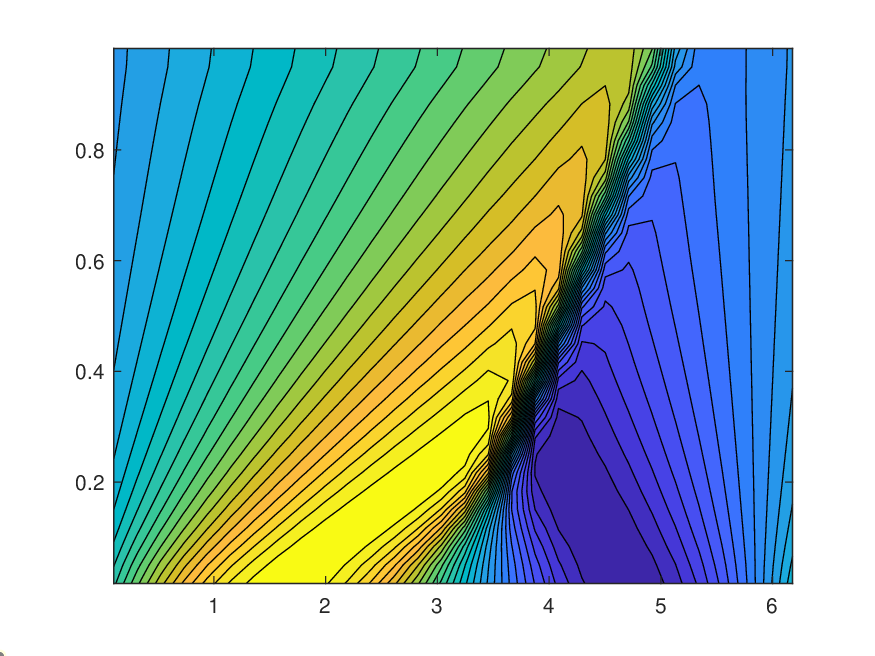}
\includegraphics[width=0.495\textwidth]{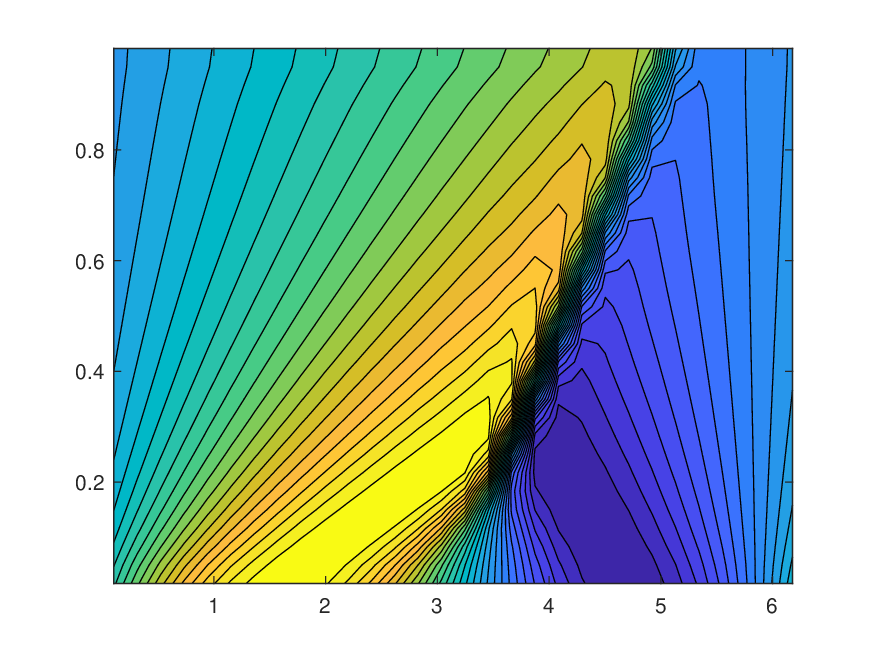}
\includegraphics[width=0.495\textwidth]{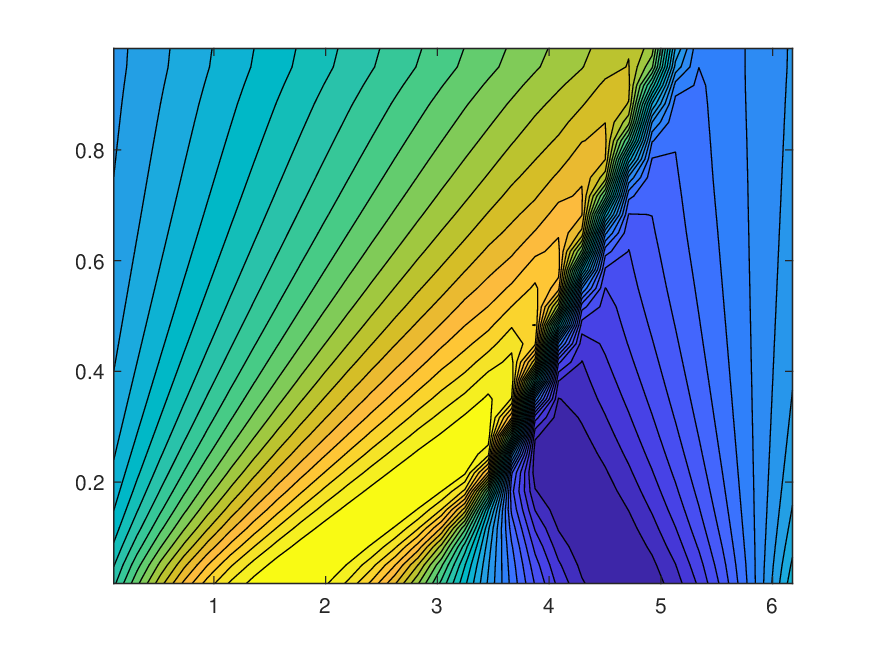}
\includegraphics[width=0.495\textwidth]{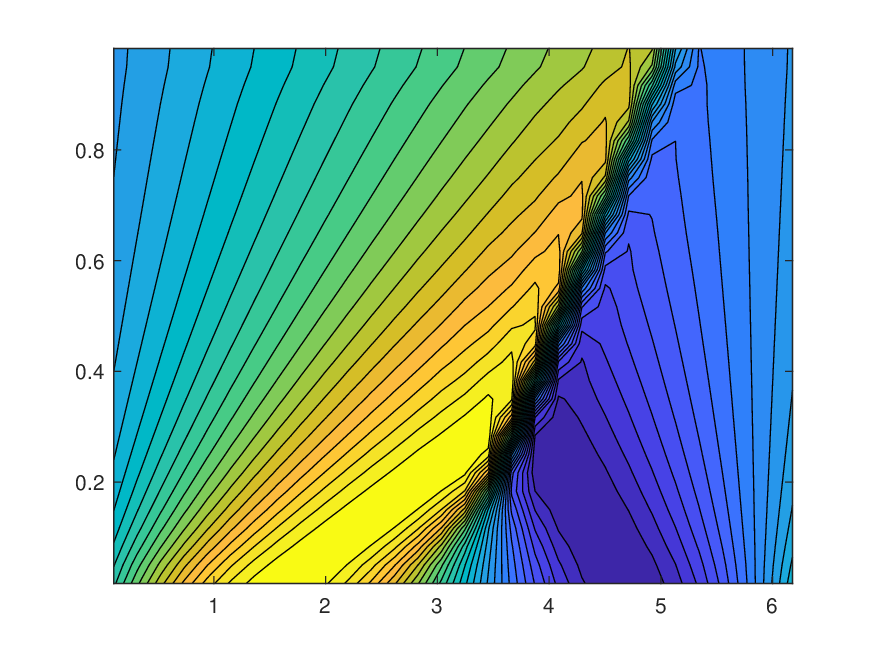}
\includegraphics[width=0.495\textwidth]{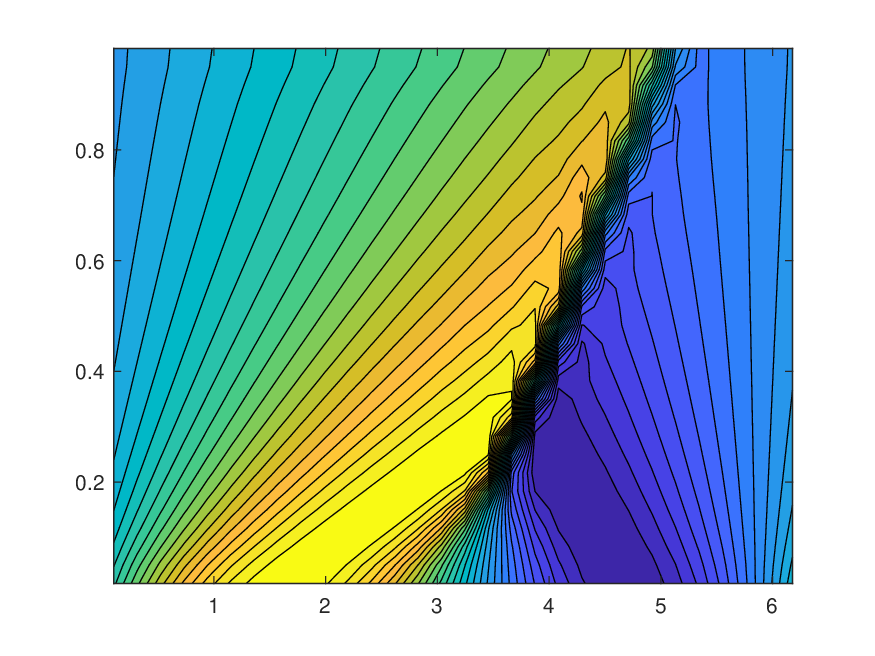}
\includegraphics[width=0.495\textwidth]{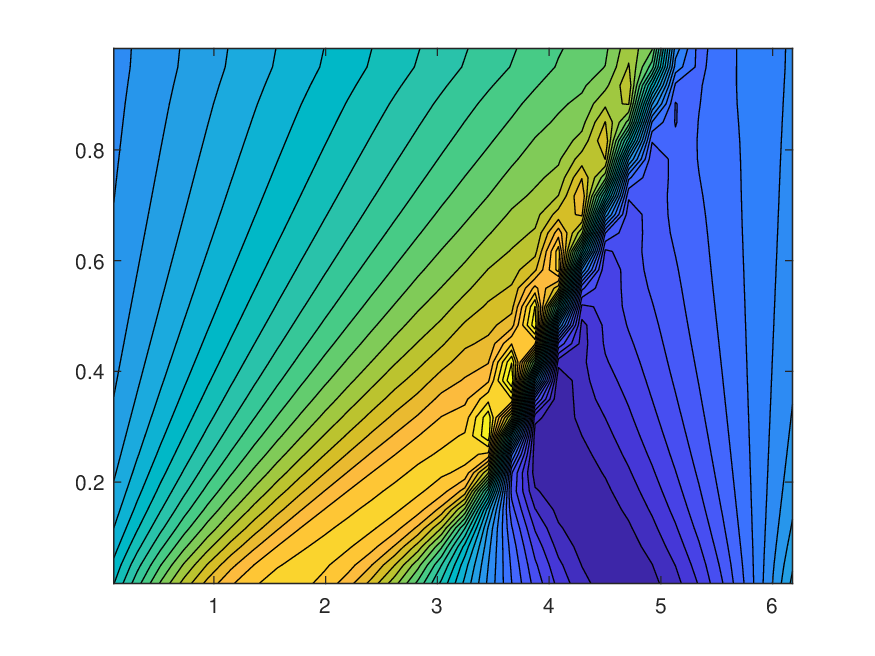}
\caption{Solution in the filled contour plot for Example \ref{ex:burgers_advection} with $\alpha=2,\, \beta=5$ at $T=1$ by WENO-JS (top left), WENO-M (top right), WENO-Z (middle left), WENO-ZR (middle right), WENO-ZL ($p=q=1$) (bottom left) and WENO-ZL ($p=3,\, q=1$) (bottom right) with $N_x = N_y = 30$.
Each contour plot displays contours at $30$ levels of the density.}
\label{fig:burgers_advection_IC2}
\end{figure}

\begin{example} \label{ex:burgers_advection}
The two-dimensional boundary layer problem \cite{ShuOsherII} is of the form
\begin{align*}
 u_t + \left(\frac{1}{2} u^2 \right)_x + u_y &= 0, \\
                                    u(x,y,0) &= \alpha + \beta \sin x, \\
                                    u(x,0,t) &= \alpha + \beta \sin x,
\end{align*}
with $u$ periodic in $x$ with the period $2 \pi$, to steady state. 
In steady state, it resembles the one-dimensional Burgers' equation with the initial condition, if $y$ is identified as the time $t$.
The computational domain $[0,\, 2 \pi] \times [0,\, 1]$ is divided into $N_x \times N_y = 30 \times 30$ uniform cells.
We employ the WENO schemes (WENO-ZL $(p,\, q)=(1,\, 1)$ and $(p,\, q)=(3\, ,1)$) to compute the boundary layer problem to the steady state at the final time $T=1$. 
Figures \ref{fig:burgers_advection_IC1} and \ref{fig:burgers_advection_IC2} contain the level curves for $\alpha = 0,\, \beta = 5$ and $\alpha = 2,\, \beta = 5$. 
\end{example}

\section{Conclusions} \label{sec:conclusions} 
In this paper, we have analyzed the finite volume WENO scheme by examining the first step of the explicit third-order TVD Runge-Kutta method and locating specific terms accountable for the numerical dissipation based on the Riemann problem of the advection equation.
A new set of Z-type nonlinear weights has been designed subsequently in order to reduce the numerical dissipation and thus sharpen the approximation near the discontinuity.
Besides, the proposed WENO-ZL scheme does not increase substantial extra computational expenses.
Numerical examples are presented by comparing with other WENO schemes, to demonstrate some advantages of the WENO-ZL scheme.
Relaxing the fixed values of $p$ in \eqref{eq:global_smooth_indicator_ZL} and $q$ in \eqref{eq:weights_ZL} needs further investigation for better resolution.
The theoretical justification for the improved behaviors of WENO-M, WENO-Z, WENO-ZR and WENO-ZL over WENO-JS at the final time, not just the first third-order TVD Runge-Kutta time step, is being explored in our ongoing work.

\section*{Acknowledgments}
This work is funded by National Research Foundation of Korea under the grant number 2021R1A2C3009648 and POSTECH Basic Science Research Institute under the NRF grant number NRF2021R1A6A1A1004294412.

\section*{Appendix}
\subsection*{Errors for the single third-order TVD Runge-Kutta time step}
\begin{equation}
\begin{aligned} 
 e^{\Nc 1}_{-2} &= - \frac{\nu}{6} \omega^{\Nc 1}_{2,-3/2} \delta, \\
 e^{\Nc 1}_{-1} &= \frac{\nu}{6} \left( \omega^{\Nc 1}_{2,-3/2} + 2 \A^{\Nc 1}_{-1/2} \right) \delta, \\
 e^{\Nc 1}_{0}  &= \frac{\nu}{6} \left( - 2 \A^{\Nc 1}_{-1/2} + \B^{\Nc 1}_{1/2} \right) \delta, \\
 e^{\Nc 1}_{1}  &= - \frac{\nu}{6} \left( \B^{\Nc 1}_{1/2} + 2 \omega^{\Nc 1}_{0,3/2} \right) \delta, \\
 e^{\Nc 1}_{2}  &= \frac{1}{3} \nu \omega^{\Nc 1}_{0,3/2} \delta,
\end{aligned}
\end{equation}

\begin{equation} 
\begin{aligned}
 e^{\Nc 2}_{-2} &= - \frac{\nu}{24} \omega^{\Nc 2}_{2,-3/2} (1-\nu) \delta + \frac{\nu}{24} \omega^{\Nc 2}_{2,-3/2} e^{\Nc 1}_0,\\                 
 e^{\Nc 2}_{-1} &= \frac{1}{24} \left[ \left( \omega^{\Nc 2}_{2,-3/2} + 2 \A^{\Nc 2}_{-1/2} \right) - \left( \omega^{\Nc 2}_{2,-3/2} + \C^{\Nc 2}_{-1/2} \right) \nu \right] \nu \delta + \frac{\nu}{24} \left[ - \left( \omega^{\Nc 2}_{2,-3/2} + \C^{\Nc 2}_{-1/2} \right) e^{\Nc 1}_0 + \omega^{\Nc 2}_{2,-1/2} e^{\Nc 1}_1 \right], \\
 e^{\Nc 2}_{0}  &= - \frac{3}{4} \nu \delta + \frac{1}{24} \left[ \left( -2 \A^{\Nc 2}_{-1/2} + \D^{\Nc 2}_{1/2} + 2 \A^{\Nc 2}_{1/2} \right) + \left( \C^{\Nc 2}_{-1/2} - \D^{\Nc 2}_{1/2} \right) \nu \right] \nu \delta \\
                & + \frac{\nu}{24} \left[ \left( \C^{\Nc 2}_{-1/2} - \D^{\Nc 2}_{1/2} + 6/\nu \right) e^{\Nc 1}_0 - \left( \omega^{\Nc 2}_{2,-1/2} + \C^{\Nc 2}_{1/2} \right) e^{\Nc 1}_1 + \omega^{\Nc 2}_{2,1/2} e^{\Nc 1}_2 \right], \\                 
 e^{\Nc 2}_{1}  &= \frac{1}{4} \nu \delta + \frac{1}{24} \left[ - \left( \D^{\Nc 2}_{1/2} + 2 \A^{\Nc 2}_{1/2} + 2 \omega^{\Nc 2}_{0,3/2} \right) + \left( \D^{\Nc 2}_{1/2} + \E^{\Nc 2}_{3/2} \right) \nu \right] \nu \delta \\
                & + \frac{\nu}{24} \left[ \left( \D^{\Nc 2}_{1/2} + \E^{\Nc 2}_{3/2} \right) e^{\Nc 1}_0 + \left( \C^{\Nc 2}_{1/2} - \D^{\Nc 2}_{3/2} + 6/\nu \right) e^{\Nc 1}_1 - \left( \omega^{\Nc 2}_{2,1/2} + \C^{\Nc 2}_{3/2} \right) e^{\Nc 1}_2 \right], \\                 
 e^{\Nc 2}_{2}  &= \frac{1}{24} \left[ 2 \omega^{\Nc 2}_{0,3/2} - \left( \B^{\Nc 2}_{3/2} + 2 \omega^{\Nc 2}_{0,5/2} \right) \nu \right] \nu \delta + \frac{\nu}{24} \left[ - \left( \E^{\Nc 2}_{3/2} + 2 \omega^{\Nc 2}_{0,5/2} \right) e^{\Nc 1}_0 \right. \\
                & \left. + \left( \D^{\Nc 2}_{3/2} + \E^{\Nc 2}_{5/2} \right) e^{\Nc 1}_1 + \left( \C^{\Nc 2}_{3/2} - \D^{\Nc 2}_{5/2} + 6/\nu \right) e^{\Nc 1}_2 \right], \\
 e^{\Nc 2}_{3}  &= \frac{1}{12} \omega^{\Nc 2}_{0,5/2} \nu^2 \delta + \frac{\nu}{24} \left[ 2 \omega^{\Nc 2}_{0,5/2} e^{\Nc 1}_0 - \left( \E^{\Nc 2}_{5/2} + 2 \omega^{\Nc 2}_{0,7/2} \right) e^{\Nc 1}_1 + \left( \D^{\Nc 2}_{5/2} + \E^{\Nc 2}_{7/2} \right) e^{\Nc 1}_2 \right], \\
 e^{\Nc 2}_{4}  &= \frac{\nu}{24} \left[ 2 \omega^{\Nc 2}_{0,7/2} e^{\Nc 1}_1 - \left( \E^{\Nc 2}_{7/2} + 2 \omega^{\Nc 2}_{0,9/2} \right) e^{\Nc 1}_2 \right], \\
 e^{\Nc 2}_{5}  &= \frac{\nu}{12} \omega^{\Nc 2}_{0,9/2} e^{\Nc 1}_2,
\end{aligned}
\end{equation}

\begin{equation}
\begin{aligned}
 \es^{\Nc 1}_{-2} &= - \frac{\nu}{9} \omega^{\Nc 3}_{2,-3/2} (1-\nu) \delta + \frac{\nu}{9} \omega^{\Nc 3}_{2,-3/2} e^{\Nc 2}_0, \\                
 \es^{\Nc 1}_{-1} &= \frac{1}{9} \left[ \left( \omega^{\Nc 3}_{2,-3/2} + 2 \A^{\Nc 3}_{-1/2} \right) - \left( \omega^{\Nc 3}_{2,-3/2} + \C^{\Nc 3}_{-1/2} \right) \nu \right] \nu \delta + \frac{\nu}{9} \left[ - \left( \omega^{\Nc 3}_{2,-3/2} + \C^{\Nc 3}_{-1/2} \right) e^{\Nc 2}_0 + \omega^{\Nc 3}_{2,-1/2} e^{\Nc 2}_1 \right], \\            
 \es^{\Nc 1}_0    &= \frac{1}{3} \nu \delta + \frac{1}{9} \left[ \left( -2 \A^{\Nc 3}_{-1/2} + \B^{\Nc 3}_{1/2} \right) + \left( \C^{\Nc 3}_{-1/2} - \D^{\Nc 3}_{1/2} \right) \nu \right] \nu \delta \\
                  & + \frac{\nu}{9} \left[ \left( \C^{\Nc 3}_{-1/2} - \D^{\Nc 3}_{1/2} + 6/\nu \right) e^{\Nc 2}_0 - \left( \omega^{\Nc 3}_{2,-1/2} + \C^{\Nc 3}_{1/2} \right) e^{\Nc 2}_1 + \omega^{\Nc 3}_{2,1/2} e^{\Nc 2}_2 \right], \\                
 \es^{\Nc 1}_1    &= \frac{1}{9} \left[ - \left( \B^{\Nc 3}_{1/2} + 2 \omega^{\Nc 3}_{0,3/2} \right) + \left( \D^{\Nc 3}_{1/2} + \E^{\Nc 3}_{3/2} \right) \nu \right] \nu \delta \\
                  & + \frac{\nu}{9} \left[ \left( \D^{\Nc 3}_{1/2} + \E^{\Nc 3}_{3/2} \right) e^{\Nc 2}_0 + \left( \C^{\Nc 3}_{1/2} - \D^{\Nc 3}_{3/2} + 6/\nu \right) e^{\Nc 2}_1 - \left( \omega^{\Nc 3}_{2,1/2} + \C^{\Nc 3}_{3/2} \right) e^{\Nc 2}_2 + \omega^{\Nc 3}_{2,3/2} e^{\Nc 2}_3 \right], \\                
 \es^{\Nc 1}_2    &= \frac{1}{9} \left[ 2 \omega^{\Nc 3}_{0,3/2} - \left( \E^{\Nc 3}_{3/2} + 2 \omega^{\Nc 3}_{0,5/2} \right) \nu \right] \nu \delta + \frac{\nu}{9} \left[ - \left( \E^{\Nc 3}_{3/2} + 2 \omega^{\Nc 3}_{0,5/2} \right) e^{\Nc 2}_0 + \left( \D^{\Nc 3}_{3/2} + \E^{\Nc 3}_{5/2} \right) e^{\Nc 2}_1 \right. \\
                  & \left. + \left( \C^{\Nc 3}_{3/2} - \D^{\Nc 3}_{5/2} + 6/\nu \right) e^{\Nc 2}_2 - \left( \omega^{\Nc 3}_{2,3/2} + \C^{\Nc 3}_{5/2} \right) e^{\Nc 2}_3 + \omega^{\Nc 3}_{2,5/2} e^{\Nc 2}_4 \right], \\
 \es^{\Nc 1}_3    &= \frac{2}{9} \omega^{\Nc 3}_{0,5/2} \nu^2 \delta + \frac{\nu}{9} \left[ 2 \omega^{\Nc 3}_{0,5/2} e^{\Nc 2}_0 - \left( \E^{\Nc 3}_{5/2} + 2 \omega^{\Nc 3}_{0,7/2} \right) e^{\Nc 2}_1 \right. \\
                  & \left. + \left( \D^{\Nc 3}_{5/2} + \E^{\Nc 3}_{7/2} \right) e^{\Nc 2}_2 + \left( \C^{\Nc 3}_{5/2} - \D^{\Nc 3}_{7/2} + 6/\nu \right) e^{\Nc 2}_3 - \left( \omega^{\Nc 3}_{2,5/2} + \C^{\Nc 3}_{7/2} \right) e^{\Nc 2}_4 + \omega^{\Nc 3}_{2,7/2} e^{\Nc 2}_5 \right], \\
 \es^{\Nc 1}_4    &= \frac{\nu}{9} \left[ 2 \omega^{\Nc 3}_{0,7/2} e^{\Nc 2}_1 - \left( \E^{\Nc 3}_{7/2} + 2 \omega^{\Nc 3}_{0,9/2} \right) e^{\Nc 2}_2 + \left( \D^{\Nc 3}_{7/2} + \E^{\Nc 3}_{9/2} \right) e^{\Nc 2}_3 \right. \\
                   & \left. + \left( \C^{\Nc 3}_{7/2} - \D^{\Nc 3}_{9/2} + 6/\nu \right) e^{\Nc 2}_4 - \left( \omega^{\Nc 3}_{2,7/2} + \C^{\Nc 3}_{9/2} \right) e^{\Nc 2}_5 \right], \\
 \es^{\Nc 1}_5    &= \frac{\nu}{9} \left[ 2 \omega^{\Nc 3}_{0,9/2} e^{\Nc 2}_2 - \left( \E^{\Nc 3}_{9/2} + 2 \omega^{\Nc 3}_{0,11/2} \right) e^{\Nc 2}_3 \right. \\
                   & \left. + \left( \D^{\Nc 3}_{9/2} + \E^{\Nc 3}_{11/2} \right) e^{\Nc 2}_4 + \left( \C^{\Nc 3}_{9/2} - \D^{\Nc 3}_{11/2} + 6/\nu \right) e^{\Nc 2}_5 \right], \\
 \es^{\Nc 1}_6    &= \frac{\nu}{9} \left[ 2 \omega^{\Nc 3}_{0,11/2} e^{\Nc 2}_3 - \left( \E^{\Nc 3}_{11/2} + 2 \omega^{\Nc 3}_{0,13/2} \right) e^{\Nc 2}_4 + \left( \D^{\Nc 3}_{11/2} + \E^{\Nc 3}_{13/2} \right) e^{\Nc 2}_5 \right], \\
 \es^{\Nc 1}_7    &= \frac{\nu}{9} \left[ 2 \omega^{\Nc 3}_{0,13/2} e^{\Nc 2}_4 - \left( \E^{\Nc 3}_{13/2} + 2 \omega^{\Nc 3}_{0,15/2} \right) e^{\Nc 2}_5 \right], \\
 \es^{\Nc 1}_8    &= \frac{2 \nu}{9} \omega^{\Nc 3}_{0,15/2} e^{\Nc 2}_5,
\end{aligned}
\end{equation}
where
\begin{align*}
\A^{\Nc k}_{j} &=    \omega^{\Nc k}_{1,j} + 2 \omega^{\Nc k}_{2,j}, \\ 
\B^{\Nc k}_{j} &=  5 \omega^{\Nc k}_{0,j} +   \omega^{\Nc k}_{1,j}, \\ 
\C^{\Nc k}_{j} &=  2 \omega^{\Nc k}_{1,j} + 5 \omega^{\Nc k}_{2,j}, \\ 
\D^{\Nc k}_{j} &= 11 \omega^{\Nc k}_{0,j} + 5 \omega^{\Nc k}_{1,j} + 2 \omega^{\Nc k}_{2,j}, \\
\E^{\Nc k}_{j} &=  7 \omega^{\Nc k}_{0,j} +   \omega^{\Nc k}_{1,j},
\end{align*}
with $\N \in \{\JS, \M, \Z, \ZR, \ZL \}$ and $j = -3/2, \cdots, 15/2$. 
The superscript $k=1,2,3$, represents the respective first, second and third stages of the single third-order TVD Runge-Kutta step.

\subsection*{2D Finite volume WENO schemes}
Let $\xi_{\gamma} = \sqrt{\frac{3}{5}}$.

Evaluating $p(x)$ and $p^{(s)}(x)$ at $x = x_{i-\xi_{\gamma}}$ gives
\begin{align*}
 v^{(S^5)}_{i-\xi_{\gamma}} &= - \frac{9+22\sqrt{15}}{2400} \vbar_{i-2} + \frac{29+41\sqrt{15}}{600} \vbar_{i-1} + \frac{1093}{1200} \vbar_i + \frac{29-41\sqrt{15}}{2400} \vbar_{i+1} - \frac{9-22\sqrt{15}}{2400} \vbar_{i+2}, \\
 v^{(0)}_{i-\xi_{\gamma}} &= \frac{2-3\sqrt{15}}{60} \vbar_{i-2} + \frac{-1+3\sqrt{15}}{15} \vbar_{i-1} + \frac{62-9\sqrt{15}}{60} \vbar_i, \\
 v^{(1)}_{i-\xi_{\gamma}} &= \frac{2+3\sqrt{15}}{60} \vbar_{i-1} + \frac{14}{15} \vbar_i + \frac{2-3\sqrt{15}}{60} \vbar_{i+1}, \\
 v^{(2)}_{i-\xi_{\gamma}} &= \frac{62+9\sqrt{15}}{60} \vbar_i - \frac{1+3\sqrt{15}}{15} \vbar_{i+1} + \frac{2+3\sqrt{15}}{60} \vbar_{i+2},
\end{align*} 
There are positive linear weights $d_0 = \frac{1008+71\sqrt{15}}{5240}, d_1 = \frac{403}{655}$ and $d_2 = \frac{1008-71\sqrt{15}}{5240}$ with their sum equal one for consistency, such that
$$
   v^{(S^5)}_{i-\xi_{\gamma}} = \sum_{s=0}^2 d_s v^{(s)}_{i-\xi_{\gamma}}.
$$
The WENO reconstruction is 
$$
   v_{i-\xi_{\gamma}} = \sum_{s=0}^2 \omega^\N_s v^{(s)}_{i-\xi_{\gamma}},
$$
where the nonlinear weights $\omega^\JS_k$ \eqref{eq:weights_JS}, $\omega^\M_k$ \eqref{eq:weights_M}, $\omega^\Z_k$ \eqref{eq:weights_Z}, $\omega^\ZR_k$ \eqref{eq:weights_ZR}, and $\omega^\ZL_k$ \eqref{eq:weights_ZL} are applied.

Evaluating $p(x)$ and $p^{(s)}(x)$ at $x = x_i$ gives
\begin{align*}
 v^{(S^5)}_i &= \frac{3}{640} \vbar_{i-2} - \frac{29}{480} \vbar_{j-1} + \frac{1067}{960} \vbar_j - \frac{29}{480} \vbar_{j+1} + \frac{3}{640} \bar{v}_{j+2}, \\
 v^{(0)}_i &= - \frac{1}{24} \vbar_{i-2} + \frac{1}{12} \vbar_{i-1} + \frac{23}{24} \vbar_i, \\
 v^{(1)}_i &= - \frac{1}{24} \vbar_{i-1} + \frac{13}{12} \vbar_i - \frac{1}{24} \vbar_{i+1}, \\
 v^{(2)}_i &=   \frac{23}{24} \vbar_i + \frac{1}{12} \vbar_{i+1} - \frac{1}{24} \vbar_{i+2},
\end{align*}
There are linear weights $d_0 = d_2 = - \frac{9}{80}$ and $d_1 = \frac{49}{40}$ such that
$$
   v^{(S^5)}_i = \sum_{s=0}^2 d_s v^{(s)}_i.
$$
Note that $d_0$ and $d_2$ are negative. 
Then the WENO procedure cannot be applied directly to obtain a stable scheme and the splitting technique in \cite{Shi} could be utilized to treat the negative weights $d_0$ and $d_2$.
The linear weights are split into positive and negative parts: 
$$
   \tilde{\gamma}^+_s = \frac{1}{2} \left( d_s + 3 | d_s | \right),~\tilde{\gamma}^-_s = \tilde{\gamma}^+_s - d_s,~s=0,1,2.
$$
Then $d_s = \tilde{\gamma}^+_s - \tilde{\gamma}^-_s$ and 
\begin{align*}
 \tilde{\gamma}^+_0 &= \frac{9}{80}, & \tilde{\gamma}^+_1 &= \frac{49}{20}, & \tilde{\gamma}^+_2 &= \frac{9}{80}; \\
 \tilde{\gamma}^-_0 &= \frac{9}{40}, & \tilde{\gamma}^-_1 &= \frac{49}{40}, & \tilde{\gamma}^+_2 &= \frac{9}{40}.
\end{align*}
We scale them by
$$
   \sigma^{\pm} = \sum_{s=0}^2 \tilde{\gamma}_s^{\pm},~\gamma_s^{\pm} = \tilde{\gamma}_s^{\pm} / \sigma^{\pm},~s=0,1,2.
$$
Then $\sigma^+ = \frac{107}{40}$, $\sigma^- = \frac{67}{40}$ and the linear positive and negative weights $\gamma_s^{\pm}$ are given by
\begin{align*}
 \gamma^+_0 &= \frac{9}{214}, & \gamma^+_1 &= \frac{98}{107}, & \gamma^+_2 &= \frac{9}{214}; \\
 \gamma^-_0 &= \frac{9}{67}, & \gamma^-_1 &= \frac{49}{67}, & \gamma^+_2 &= \frac{9}{67},
\end{align*}
which satisfy 
\begin{equation} \label{eq:linear_weight_relation}
 d_s = \sigma^+ \gamma_s^+ - \sigma^- \gamma_s^-.
\end{equation}
The WENO reconstruction is 
$$
   v_i = \sum_{s=0}^2 \omega^\N_s v^i,
$$
where the nonlinear weights, based on the relation \eqref{eq:linear_weight_relation} for the linear weights, are defined by
$$
   \omega^\N_s = \sigma^+ \omega^{\Nc +}_k - \sigma^- \omega^{\Nc -}_k.
$$
The nonlinear positive and negative weights $\omega^{\JS,\pm}_k$ are defined as
\begin{equation} \label{eq:weights_JS_pm}
 \omega^{\JS,\pm}_s = \frac{\alpha^{\JS,\pm}_s}{\sum^2_{r=0} \alpha^{\JS,\pm}_r},~\alpha^{\JS,\pm}_s = \frac{\gamma^{\pm}_s}{(\beta_s + \epsilon)^2},~s=0,1,2.
\end{equation}
The mapped nonlinear positive and negative weights $\omega^{\M,\pm}_s$ are defined as
$$
   \omega^{\M,\pm}_s = \frac{\alpha^{*,\pm}_s}{\sum_{r=0}^2 \alpha^{*,\pm}_r},~\alpha^{*,\pm}_s = g_s(\omega^{\JS,\pm}_s),~s=0,1,2,
$$
with $\omega^{\JS,\pm}_s$ in \eqref{eq:weights_JS_pm}.
The Z-type nonlinear positive and negative weights $\omega^{\Z,\pm}_s$ are defined as 
$$
   \omega^{\Z,\pm}_s = \frac{\alpha^{\Z,\pm}_s}{\sum_{r=0}^2 \alpha^{\Z,\pm}_r},~\alpha^{\Z,\pm}_s = \gamma^{\pm}_s \left( 1 + \frac{\tau_5}{\beta_s + \epsilon} \right),~s=0,1,2,
$$
with $\tau_5 = \left| \beta_0 - \beta_2 \right|$.
The nonlinear positive and negative weights $\omega^{\ZR,\pm}_s$ are 
$$
   \omega^{\ZR,\pm}_s = \frac{\alpha^{\ZR,\pm}_s}{\sum_{r=0}^2 \alpha^{\ZR,\pm}_r},~\alpha^{\ZR,\pm}_s = \gamma^{\pm}_s \left( 1 + \left( \frac{\tau^\ZR_5}{\sqrt[p]{\beta_s} + \epsilon} \right)^p \right),~s=0,1,2,
$$
with $\tau^\ZR_5 = \left| \sqrt[p]{\beta_0} - \sqrt[p]{\beta_2} \right|$.
The nonlinear positive and negative weights $\omega^{\ZL,\pm}_s$ are 
$$
   \omega^{\ZL,\pm}_s = \frac{\alpha^{\ZL,\pm}_s}{\sum_{r=0}^2 \alpha^{\ZL,\pm}_r},~\alpha^{\ZL,\pm}_s = \gamma^{\pm}_s \left( 1 + \left( \frac{\tau^\ZL_5}{\beta_s + \epsilon} \right)^q \right),~s=0,1,2,
$$
with $\tau^\ZL_5 = \frac{1}{p} \left| \ln \frac{1+\beta_0}{1+\beta_2} \right|$.

Evaluating $p(x)$ and $p^{(s)}(x)$ at $x = x_{i+\xi_{\gamma}}$ gives
\begin{align*}
 v^{(S^5)}_{i+\xi_{\gamma}} &= \frac{-9+22\sqrt{15}}{2400} \vbar_{i-2} + \frac{29-41\sqrt{15}}{600} \vbar_{i-1} + \frac{1093}{1200} \vbar_i + \frac{29+41\sqrt{15}}{2400} \vbar_{i+1} - \frac{9+22\sqrt{15}}{2400} \vbar_{i+2}, \\
 v^{(0)}_{i+\xi_{\gamma}} &= \frac{2+3\sqrt{15}}{60} \vbar_{i-2} - \frac{1+3\sqrt{15}}{15} \vbar_{i-1} + \frac{62+9\sqrt{15}}{60} \vbar_i, \\
 v^{(1)}_{i+\xi_{\gamma}} &= \frac{2-3\sqrt{15}}{60} \vbar_{i-1} + \frac{14}{15} \vbar_i + \frac{2+3\sqrt{15}}{60} \vbar_{i+1}, \\
 v^{(2)}_{i+\xi_{\gamma}} &= \frac{62-9\sqrt{15}}{60} \vbar_i + \frac{-1+3\sqrt{15}}{15} \vbar_{i+1} + \frac{2-3\sqrt{15}}{60} \vbar_{i+2},
\end{align*}
There are positive linear weights $d_0 = \frac{1008-71\sqrt{15}}{5240}, d_1 = \frac{403}{655}$ and $d_2 = \frac{1008+71\sqrt{15}}{5240}$ with their sum equal one for consistency, such that
$$
   v^{(S^5)}_{i+\xi_{\gamma}} = \sum_{s=0}^2 d_s v^{(s)}_{i+\xi_{\gamma}}.
$$
The WENO reconstruction is 
$$
   v_{i+\xi_{\gamma}} = \sum_{s=0}^2 \omega^\N_s v^{(s)}_{i+\xi_{\gamma}},
$$
after applying the nonlinear weights $\omega^\JS_k$ \eqref{eq:weights_JS}, $\omega^\M_k$ \eqref{eq:weights_M}, $\omega^\Z_k$ \eqref{eq:weights_Z}, $\omega^\ZR_k$ \eqref{eq:weights_ZR}, and $\omega^\ZL_k$ \eqref{eq:weights_ZL}.


\end{document}